%%%%%%%%%%%% File brc9.tex with figures -- version june 29 1999 %%%%
\magnification=1200

\catcode`\@=11

\hsize=125 mm   \vsize =187mm
\hoffset=4mm    \voffset=10mm
\pretolerance=500 \tolerance=1000 \brokenpenalty=5000

\catcode`\;=\active
\def;{\relax\ifhmode\ifdim\lastskip>\z@
\unskip\fi\kern.2em\fi\string;}

\catcode`\:=\active
\def:{\relax\ifhmode\ifdim\lastskip>\z@\unskip\fi
\penalty\@M\ \fi\string:}

\catcode`\!=\active
\def!{\relax\ifhmode\ifdim\lastskip>\z@
\unskip\fi\kern.2em\fi\string!}

\catcode`\?=\active
\def?{\relax\ifhmode\ifdim\lastskip>\z@
\unskip\fi\kern.2em\fi\string?}

\def\^#1{\if#1i{\accent"5E\i}\else{\accent"5E #1}\fi}
\def\"#1{\if#1i{\accent"7F\i}\else{\accent"7F #1}\fi}

%\frenchspacing

\catcode`\@=12

\newif\ifpagetitre      \pagetitretrue
\newtoks\hautpagetitre  \hautpagetitre={\hfil}
\newtoks\baspagetitre   \baspagetitre={\hfil}

\newtoks\auteurcourant  \auteurcourant={\hfil}
\newtoks\titrecourant   \titrecourant={\hfil}

\newtoks\hautpagegauche \newtoks\hautpagedroite
\hautpagegauche={\hfil\the\auteurcourant\hfil}
\hautpagedroite={\hfil\the\titrecourant\hfil}

\newtoks\baspagegauche  \baspagegauche={\hfil\tenrm\folio\hfil}
\newtoks\baspagedroite  \baspagedroite={\hfil\tenrm\folio\hfil}

\headline={\ifpagetitre\the\hautpagetitre
\else\ifodd\pageno\the\hautpagedroite
\else\the\hautpagegauche\fi\fi}

\footline={\ifpagetitre\the\baspagetitre
\global\pagetitrefalse
\else\ifodd\pageno\the\baspagedroite
\else\the\baspagegauche\fi\fi}

%%%%%%%%
%\hautpagetitre={\hfill\tenbf preliminary version: Not for diffusion\hfill}
\hautpagetitre={\hfill}
\hautpagegauche={\tenrm\folio\hfill\the\auteurcourant}
\hautpagedroite={\tenrm\the\titrecourant\hfill\folio}
\baspagegauche={\hfil} \baspagedroite={\hfil}
\auteurcourant{Marmi, Moussa, Yoccoz}
\titrecourant{Complex Brjuno Functions}
%%%%%%%%%
\def\mois{\ifcase\month\or January\or February\or March\or April\or
May\or June\or July\or August\or September\or October\or November\or
December\fi}
\def\Date{\rightline{\mois\ /\ \the\day\ /\/ \the\year}}
\hfuzz=0.3pt
\font\tit=cmb10 scaled \magstep1
\def\dst{\displaystyle}

\def\IM{\mathop{\Im m}\nolimits}
\def\RE{\mathop{\Re e}\nolimits}
\def\H{\Bbb H}
\def\R{\Bbb R}
\def\T{\Bbb T}
\def\Z{\Bbb Z}
\def\Q{\Bbb Q}
\def\C{\Bbb C}
\def\N{{\Bbb N}}
\def\D{{\Bbb D}}

\def\S{{\Bbb S}}
\def\hH+{\hat{\H}^{+}}
\def\hHZ+{\widehat{\H^{+}/\Z}}

\def\cM{{\cal M}}

\def\be{\beta}

\def\eps{{\varepsilon}}

\def\abcd{\left(\matrix{ a & b\cr c & d \cr}\right)}
\def\O#1#2{{\cal O}^{#1}(\overline{\Bbb C}\setminus #2)}
\def\hC#1{\hat{\cal C}_{\Bbb R}(\overline{\Bbb C}\setminus #1)}
\def\Li{\hbox{Li}_{2}}
\def\cotg{\hbox{\rm cotg}\,}

\def\Lp#1#2{L^{#1}(#2)}
\def\Hp#1#2{H^{#1}(#2)}

\def\remark#1{\noindent{\it Remark\ }#1\ }
\def\proof{\noindent{\it Proof.\ }}
\def\qed{\hfill$\square$\par\smallbreak}
\def\Proc#1#2\par{\medbreak \noindent {\bf #1\enspace }{\sl #2}%
\par\ifdim \lastskip <\medskipamount \removelastskip \penalty 55\medskip \fi}%
\def\Def#1#2\par{\medbreak \noindent {\bf #1\enspace }{#2}%
\par\ifdim \lastskip <\medskipamount \removelastskip \penalty 55\medskip
\fi}
\def\qed{\hfill$\square$\par\smallbreak}
%\long\def\Proc#1#2\par{\medbreak \noindent {\bf #1\enspace }{\sl #2}%
%\par\ifdim \lastskip <\medskipamount \removelastskip \penalty 55\medskip
%\fi}
\def\hfl#1#2{\smash{\mathop{\hbox to 12mm{\rightarrowfill}}
\limits^{\scriptstyle#1}_{\scriptstyle#2}}}
\def\vfl#1#2{\llap{$\scriptstyle #1$}
\left\downarrow\vbox to 6mm{}\right.\rlap{$\scriptstyle #2$}}
\def\diagram#1{\def\normalbaselines{\baselineskip=0pt
\lineskip=10pt\lineskiplimit=1pt} \matrix{#1}}
\def\smk{\smallskip}
\def\mdk{\medskip}
\def\bgk{\bigskip}
%%% fin macros %%%%%%%%

%%%% new fonts %%%%%%%%

\font\pic=cmr9

%%%%%%%%%  Title  and Abstract  %%%%%%%%%%%%%%
\smallskip
\centerline{\tit Complex Brjuno Functions}
\medskip
\centerline{S. Marmi\footnote{$^1$}{\pic Dipartimento di Matematica ``U.
Dini'',
Universit\`a di Firenze, Viale Morgagni 67$/$A, 50134 Firenze, Italy},
P. Moussa\footnote{$^2$}{ \pic Service de Physique Th\'eorique, CEA/Saclay,
91191 Gif-Sur-Yvette, France}
and J.-C. Yoccoz\footnote{$^3$}{\pic Coll\`ege de France, 3, Rue d'Ulm,
75005 Paris and Universit\'e de Paris-Sud, Math\'ematiques. B\^at. 425, 
91405-Orsay, France} \footnote{}{\pic\hfill June 29, 1999}}
\bigskip
\centerline{{\bf  Abstract}}
\medskip
\input amssym.def
\input amssym.tex
The Brjuno function arises naturally in the study of one--dimensional
analytic small divisors problems. It belongs to $\hbox{BMO}({\Bbb
T}^{1})$ and it is stable under H\"older perturbations. It is related
to the size of Siegel disks by various rigorous and conjectural
results.

In this work we show how to extend the Brjuno function to a holomorphic
function on ${\Bbb H}/{\Bbb Z}$, the complex Brjuno function. This has
an explicit expression in terms of a series of transformed
dilogarithms under the action of the modular group.

The extension is obtained using a complex analogue of the
continued fraction expansion of a real number. Since our method
is based on the use of hyperfunctions it applies to less regular
functions than the Brjuno function and it is quite general.

We prove that the harmonic conjugate of the Brjuno function is
bounded. Its trace on ${\Bbb R}/{\Bbb Z}$ is continuous at all irrational
points and has a jump of $\pi/q$ at each rational point
$p/q\in {\Bbb Q}$.\par
%%%%%%%%%%%%%  1.  Introduction  %%%%%%%%%%%%%%
\beginsection 1. Introduction \par
%\medskip
\beginsection 1.1 The real Brjuno function\par
%\smallskip

Let $\alpha\in {\Bbb R}\setminus {\Bbb Q}$ and let
$(p_{n}/q_{n})_{n\ge 0}$ be the sequence of the convergents of
its continued fraction expansion. A {\it Brjuno number} is an
irrational number $\alpha$ such that $\sum_{n=0}^\infty {\log
q_{n+1}\over q_{n}}<+\infty$.

The importance of Brjuno numbers comes from the study of
one--dimensional analytic small divisors problems. In the case of
germs of holomorphic diffeomorphisms of one complex variable
with an indifferent fixed point, extending a previous result
of Siegel [Si], Brjuno proved [Br] that all germs with linear part
$\lambda = e^{2\pi i \alpha}$ are linearizable if $\alpha$ is
a Brjuno number. Conversely the third author proved that this
condition is also necessary [Yo1]. Similar results hold for the local
conjugacy of analytic diffeomorphisms of the circle [KH, Yo2, Yo3] and
for some area--preserving maps [Ma,Da1], including the standard family
[Da2, BG1, BG2].
The set of Brjuno numbers is invariant under the action of the modular
group $\hbox{PGL}\,(2,{\Bbb Z})$ and it can be characterized as
the set where the {\it Brjuno function} $B\, :{\Bbb R}\setminus {\Bbb Q}
\rightarrow {\Bbb R}\cup\{+\infty\}$ is finite.
\vfill\eject
This arithmetical
function is ${\Bbb Z}$--periodic and satisfies the remarkable
functional equation
$$
B(\alpha ) = -\log \alpha +\alpha B\left({1\over \alpha}\right)\;
\; , \; \; \;
\alpha \in (0,1)\; , \eqno(1.1)
$$
which allows $B$ to be interpreted as a {\it cocycle} under the
action of the modular group (see Appendix 5 for
details). In terms of the continued fraction expansion of $\alpha$
the Brjuno function is defined as follows:
$$
B(\alpha ) = \sum_{j=0}^{+\infty}\beta_{j-1}(\alpha
)\log\alpha_{j}^{-1}\; , \; \eqno(1.2)
$$
where
$
\beta_{-1}=1\; , \;\;\; \beta_{j}(\alpha ) = |p_{j}-q_{j}\alpha |
\; (j\ge 0)\; ,\;\;\; \alpha_{j}= -{\dst q_{j}\alpha-p_{j}\over\dst
q_{j-1}
\alpha - p_{j-1}}\; ,
$
(see Appendix 1 for a short summary of the relevant facts concerning
the continued fraction).

In a previous paper [MMY] we introduced the linear operator
$$
Tf(x) = xf\left({1\over x}\right) \; , \;\; \; x\in (0,1)
$$
acting in the space of
${\Bbb Z}$--periodic measurable functions and we studied the
equation
$$
(1-T)B_{f} = f\;  ,
$$
so that
$$
\eqalign{
&B_{f}(x+1) = B_{f}(x) \;\;\;\forall x\in \R\; , \cr
&B_{f}(x) = f(x)+xB_{f}(1/x)\;\;\;\forall x\in (0,1)\; . \cr}
$$
The choice $f(x)=-\log \{x\}$ (where $\{\cdot\}$ denotes fractional
part) leads to the Brjuno function $B$. For other choices of the
singular behaviour of $f$ at $0$ the condition $B_{f}<+\infty$ leads
to different diophantine conditions. On the other hand if $f$ is
H\"older continuous then $B_{f}$ is also H\"older continuous and this
fact could help to explain the numerical results of Buric, Percival
and Vivaldi [BPV].

Acting on $L^p([0,1])$ the operator $T$ has spectral radius bounded
above by ${\sqrt{5}-1\over 2}$ (thus ($1-T$) is invertible). A suitable
adaptation of this argument has led us to conclude that the Brjuno
function belongs to $\hbox{BMO}\,({\Bbb T}^{1})$ (bounded mean
oscillation, see references [Ga, GCRF] for its definition and
more information).

By Fefferman's duality theorem BMO is the dual of the Hardy space
$\hbox{H}^{1}$ thus one can add an $L^\infty$ function to $B$ so that
the harmonic conjugate of the sum will also be $L^\infty$. This
suggests to look for an holomorphic function ${\cal B}$ defined on
the upper half plane which is ${\Bbb Z}$--periodic and whose
trace on ${\Bbb R}$ has for imaginary part the Brjuno function $B$.
The function ${\cal B}$ will be called the {\it complex Brjuno
function}.

Another motivation for the introduction of the complex Brjuno function
comes from results concerning the linearization of the quadratic
polynomial $P_{\lambda}(z) = \lambda (z-z^{2})$
([Yo1], Chapter II). One has the following results:
\item{(1)} there exists a bounded holomorphic function $U\, : {\Bbb D}
\rightarrow {\Bbb C}$ such that $|U(\lambda )|$ is equal to the
radius of convergence of the normalized linearization of
$P_{\lambda}$;
\item{(2)} for all $\lambda_{0} \in {\Bbb S}^{1}$, $|U(\lambda )|$ has a
non--tangential limit in $\lambda_{0}$ (which is still equal to
the
radius of convergence of the normalized linearization of
$P_{\lambda_{0}}$);
\item{(3)} if $\lambda = e^{2\pi i\alpha}$, $\alpha\in {\Bbb
R}\setminus {\Bbb Q}$, $P_{\lambda}$ is linearizable if and only if
$\alpha$ is a Brjuno number. Moreover there exists a universal
constant $C_{1}>0$ and for all $\varepsilon >0$ there exists
$C_{\varepsilon}>0$ such that for all Brjuno numbers $\alpha$
one has
$$
(1-\varepsilon )B(\alpha ) -C_{\varepsilon} \le -\log |U(\lambda )|
\le B(\alpha )+C_{1}\; .
$$

\noindent
In [MMY] the authors proposed the following conjecture (see also [Ma]): the
function defined on the set of Brjuno numbers by
$\alpha \mapsto B(\alpha )+\log |U(e^{2\pi i \alpha })|$
extends to a $1/2$--H\"older continuous function as $\alpha $
varies in ${\Bbb R}$.
If this were true then the function $-i{\cal B}(z)+\log U(e^{2\pi i
z})$ would also extend to a H\"older continuous function
on $\overline{ \H}$.

%\vskip .5 truecm
\beginsection
1.2 The complex Brjuno function\par
%\vskip .3 truecm

A natural question now is how to extend the operator $T$ to complex
analytic functions. This is achieved as follows:
the operator $T$ extends to the space $A'([0,1])$ of
hyperfunctions $u$ with support contained in $[0,1]$ (see section 1.4
for a proof of this fact and
Appendix 2 for a very brief introduction to hyperfunctions). This
space is canonically isomorphic to the complex vector space
${\cal O}^{1}(\overline{\Bbb C}\setminus [0,1])$ of holomorphic
functions on $\overline{\Bbb C}\setminus [0,1]$ which vanish at
infinity. The connection between $u$ and the associated holomorphic function
$\varphi$ is commonly written as: $u(x)={1\over
2i}(\varphi(x+i0)-\varphi(x-i0))$,
which is also equal to $\IM \varphi(x+i0)$ when $\varphi$ is real (i.e.
$\varphi(\overline{z})=\overline{\varphi(z)}$).
On ${\cal O}^{1}(\overline{\Bbb C}\setminus [0,1])$
the formula for $T$ reads
$$
(T\varphi)(z) = -z\sum_{m=1}^\infty \left[\varphi\left({1\over z}-m\right)
-\varphi (-m)\right]+\sum_{m=1}^\infty \varphi'(-m)\; .
\eqno(1.3)
$$
Formally we have
$$
(1-T)^{-1} \varphi (z) = \sum_{r\ge 0}(T^r\varphi)(z)=
\sum_{g\in {\cal M}}(L_{g}\varphi)(z)\; , \; \eqno (1.4)
$$
where the monoid
$$
{\cal M} = \left\{ g=\left(\matrix{a & b\cr c & d\cr}\right)\in
\hbox{GL}\,(2,{\Bbb Z})\, ,
\; d\ge b \ge a\ge 0\, ,\; \hbox{and}\; d\ge c\ge a\;
\right\} \cup\left\{\left(\matrix{1 & 0\cr 0 &
1\cr}\right)\right\}\, ,
$$
acts on ${\cal O}^{1}(\overline{\Bbb C}\setminus [0,1])$
according to
$$
(L_{g}\varphi)(z) = (a-cz)\left[\varphi\left({dz-b\over a-cz}\right)-
\varphi\left(-{d\over c}\right)\right]-\hbox{det}(g)c^{-1}\varphi'
\left(-{d\over c}\right)\; . \eqno(1.5)
$$
The series (1.4) actually converges in
${\cal O}^{1}(\overline{\Bbb C}\setminus [0,1])$
to a function $\sum_{\cal M}\varphi$. To recover a holomorphic
periodic function on $\H$ one sums over integer translates:
$$
{\cal B}_{\varphi}(z)= \sum_{n\in {\Bbb Z}}\left(\sum_{\cal
M}\varphi\right)(z-n)\; . \eqno(1.6)
$$
To construct the complex Brjuno function one has to take
$\varphi_{0}(z)=-{1\over \pi}\hbox{Li}_{2}\left({1\over z}\right)$,
where $\hbox{Li}_{2}$ is the dilogarithm (Appendix 3, [O]).
Then the above formulas give
$$
\eqalign{
{\cal B}(z) &= -{1\over \pi}\sum_{p/q\in {\Bbb Q}}\left\{
(p'-q'z)\left[\hbox{Li}_{2}\left({p'-q'z\over qz-p}\right)-
\hbox{Li}_{2}\left(-{q'\over q}\right)\right]\right.\cr
&\left. +(p''-q''z)\left[
\hbox{Li}_{2}\left({p''-q''z\over qz-p}\right)-
\hbox{Li}_{2}\left(-{q''\over q}\right)\right]+{1\over q}
\log{q+q''\over q+q'}\right\}\; ,\cr} \eqno(1.7)
$$
where $\left[{p'\over q'},{p''\over q''}\right]$ is the Farey interval
such that ${p\over q}={p'+p''\over q'+q''}$ (with the convention
$p'=p-1$, $q'=1$, $p''=1$, $q''=0$ if $q=1$).

%\vskip .5 truecm
\beginsection
1.3 Main results: Properties of the complex Brjuno function\par
%\vskip .3 truecm

Our main result (Corollary 5.8)
is that {\it the real part of ${\cal B}$ is
bounded} on the upper half plane (note that this statement is
stronger than the result obtained by
the above mentioned general properties of the harmonic
conjugates of BMO functions); moreover {\it the trace of $\RE{\cal
B}$ on ${\Bbb R}$ is continuous at all irrational points and has a
jump of $\pi /q$ at each rational point $p/q\in {\Bbb Q}$}
(see Section 5.2.9).

A numerical study of the function $\arg U(e^{2\pi i \alpha})$ seems
to indicate a similar behaviour (see Figure 1).

Concerning the boundary behaviour of the imaginary part of ${\cal B}$
we prove the following (Theorem 5.19):
\item{(i)} if $\alpha$ is a Brjuno number then $\IM {\cal B}(\alpha
+w)$ converges to $B(\alpha )$ as $w\rightarrow 0$ in any domain with
a finite order of tangency to the real axis;
\item{(ii)} if $\alpha$ is diophantine one can allow domains with
infinite order of tangency (see (5.58))

The precise behaviour of at rational points is described by Theorem
5.10.

%\vskip .5 truecm
\beginsection 1.4 Hyperfunctions and operator $T$\par
%\vskip .3 truecm

Let $u,\psi\in \Lp{2}{[0,1]}$, $m\in\N$, $m\ge 1$. We consider
$$
T_{m}u(x) =\cases{
xu(1/x-m)\; & if $x\in\left[{1\over m+1},{1\over m}\right]$\ ,\cr
\hfill 0\hfill & \hfill otherwise.\hfill} \eqno(1.8)
$$
Note that $T_{m}u=(Tu)\mid_{\left[{1\over m+1},{1\over
m}\right]}$, thus $T=\sum_{m\ge 1}T_{m}$.
We define  the adjoint $T_{m}^{*}$ by
$$
\int_0^1 T_{m}u(x)\psi(x)dx = \int_{0}^{1}u(x)T_{m}^{*}\psi
(x) dx\; ,
$$
which gives
$$
T_{m}^{*}\psi (x) = {1\over (m+x)^3}
\psi\left({1\over m+x}\right)\; .
$$
The previous formula with $\psi$ analytic
in a neighborhood of $[0,1]$ allows to extend the domain of definition
of $T_{m}$ to the hyperfunctions $u\in A'([0,1])$
(see Appendix 5 for a very short summary of hyperfunctions)
and to obtain $T_{m}u\in A'\left(\left[{1\over m+1},{1\over
m}\right]\right)$. More generally, if $u\in
A'([\gamma_{0},\gamma_{1}])$, $\gamma_{0}>-1$, then
$T_{m}u\in A'\left(\left[{1\over m+\gamma_{1}},{1\over m+\gamma_{0}}
\right]\right)
\subset A'\left(\left[
0,{1\over 1+\gamma_{0}}\right]\right)$. One has
$$
\int_{0}^{1}Tu(x)\psi (x)dx=\sum_{m\ge 1}\int_{0}^{1}u(x)T_{m}^{*}
\psi (x) dx\; . \eqno(1.9)
$$
If $\psi$ is holomorphic in a neighborhood $V$ of $[0,1]$, then
also $T_{m}^{*}\psi$ is holomorphic in $V$ and one has
$$
\sup_{V}|T_{m}^{*}\psi| \le {1\over 2}m^{-3}\sup_{V}|\psi |\; ,
\eqno(1.10)
$$
(this follows immediately from the estimates of section 3.2 choosing
$V$ to be the complement of a neighborhood of $D_{\infty}$ with
respect to the Poincar\'e metric on $\overline{\Bbb C}\setminus
[0,1]$).
Therefore the series $\sum_{m\ge 1}T_{m}u$ converges in $A'([0,1])$ to
a hyperfunction which will be denoted $Tu$.

Let $u\in A'([\gamma_{0},\gamma_{1}])$, $\gamma_{0}>-1$,
and let $\varphi\in
{\cal O}^{1}(\overline{\Bbb C}\setminus [\gamma_{0},\gamma_{1}])$
be the associated holomorphic function, i.e. $\varphi$ is holomorphic
outside $[\gamma_0,\gamma_1]$ and vanishes at infinity.
For all $m\ge 1$ the holomorphic function
associated to $T_{m}u$ is $L_{g(m)}\varphi$ given by equation
(1.5) with $a=0$, $b=c=1$, $d=m$. Indeed if $z\notin
[\gamma_{0},\gamma_{1}]$
$$
(T_{m}u)(c_{z})= u(T_{m}^{*}c_{z}), \eqno(1.11)
$$
where $c_z(x)={\dst 1\over\dst \pi}{\dst 1\over\dst x-z}$, and
$$
\eqalign{
\pi T_{m}^{*}c_{z}(x) &= {1\over (m+x)^3}
{1\over {1\over m+x}-z}
%\cr &
= {1\over (m+x)^2}-z\left({1\over m+x-{1\over z}}-
{1\over m+x}\right)
\cr&
= \pi \left[ \left.{\partial\over\partial z} c_{z}(x)
\right\vert_{z=-m} -z\left(c_{{1\over
z}-m}(x)-c_{-m}(x)\right)\right]\; . \cr
}
$$
Thus we are led to define $T \varphi$,  according to (1.3),
as an element of the space ${\cal O}^{1}(\overline{\Bbb C}\setminus [
0,1/(1+\gamma_{0})])$.

To construct the complex analytic extension of the functions $B_{f}$
(defined in Section 1.1) our strategy is the following:
\item{1)} take the restriction of the periodic function $f$ to the
interval $[0,1]$;
\item{2)} consider its associated hyperfunction $u_{f}$ and its
holomorphic representative $\varphi\in\O{1}{[0,1]}$.

Then the series (1.6) converges (thanks to Corollary 3.6) to the
complex extension ${\cal B}_{f}$ of the function $B_{f}$. The main
difficulty (unless $f$ belongs to some $L^p$ space, see Section 4.3)
would be to recover  $B_{f}$ as non--tangential limit of
the imaginary part of ${\cal B}_{f}$ as $\IM z\rightarrow 0$.

%\vskip .5 truecm
\beginsection 1.5 Summary of the contents\par
%\vskip .3 truecm

Let us now briefly describe   the contents of this article.

In Section 2 we discuss the relation between the monoid ${\cal M}$ and
the full modular group $\hbox{GL}\,(2,{\Bbb Z})$.
We then describe various automorphic actions of ${\cal M}$.

In Section 3, the introduction of a complex analogue of the continued
fraction expansion of a real number allows us to prove the
convergence of the series (1.4) and (1.6) (Corollary 3.6). The main
feature of the complex continued fraction is that it reduces to the
real continued fraction when the number is real and it stops
after a finite number of iterations when the number is rational or
complex. In the latter case the absolute value of the imaginary part
of the iterates grows at least exponentially with the number of
iterations and when it reaches $1/2$ the iteration stops.

In Section 4 we use the complex continued fraction to study the
behaviour of the series (1.4) when $z$ is close to $[0,1]$.
This is interesting in itself and it will be very important
when applied to the complex Brjuno function in order to prove our
main results. Our study allows us to prove that the restriction of $T$
to the Hardy spaces $\Hp{p}{\overline{\Bbb C}\setminus [0,1]}\cap
\O{1}{[0,1]}$, $1\le p\le +\infty$, is continuous with spectral
radius bounded above by ${\sqrt{5}-1\over 2}$. The same results holds
also on the space of functions $\varphi\in \O{1}{[0,1]}$ with bounded
real part.

The complex Brjuno function is finally introduced in Section 5 where
we state and prove our main results.

In the Appendices we recall the results we need on the real continued
fraction, on the hyperfunctions and on the dilogarithm. Then we show
how to relate the complex Brjuno function with the even real Brjuno
function treated in [MMY]. Finally we describe how  the real Brjuno
function  can be viewed as a cocycle under the action of the modular
group.

\medskip
{\it Acknowledgements.} This work begun during a visit of the first
author at the S.Ph.T.--CEA Saclay and at the Dept. of Mathematics of
Orsay during the academic year 1993--1994. This research has been
supported by the  CNR, CNRS, INFN, MURST and a EEC grant.

%\vskip .5 truecm
%%%%%%%%%%%%%%  2. Modular group, monoid etc...%%%%%%%%%%%
\beginsection 2. Modular group, the monoid $\cM$ and its action  \par
%\vskip .5 truecm
According to (1.4) the inversion of $(1-T)$ leads to consider the
monoid ${\cal M}$ of matrices of $\hbox{GL}(2,{\Bbb Z})$. In this
section we study its algebraic properties and we describe various
actions of the modular group and of the monoid ${\cal M}$ on
meromorphic or holomorphic functions.

\mdk
\noindent
{\bf  2.1
Algebraic properties, notations, structure of the monoid ${\cal
M}$; relations of $\cM$ with the modular group and with Farey
intervals}

\smk
{\bf 2.1.1} Notations:
\item{$G$}  $= \hbox{GL}(2,{\Bbb Z})= \{ \left(\matrix{a & b\cr c &
d\cr}\right)
\, , \; a,b,c,d\in {\Bbb Z}\, , \; \varepsilon_{g}:= ad-bc =\pm 1\}$;
\item{$H$}  is the subgroup of order $8$ of matrices of the form
$\left(\matrix{\varepsilon & 0\cr 0 & \varepsilon'\cr}\right)$ or
$\left(\matrix{0 & \varepsilon \cr \varepsilon' & 0\cr}\right)$,
where $\varepsilon , \varepsilon'\in \{-1,+1\}$;
\item{${\cal M}$}  is the monoid with unit $\left(\matrix{1 & 0\cr 0 &
1\cr}\right)$ made of matrices $g=\left(\matrix{a & b\cr c &
d\cr}\right)\in G$ such that, if $g\not= \,\hbox{id}\,$, we have
$d\ge b \ge a\ge 0$ and $d\ge c\ge a$.
\item{$Z$}  is the subgroup of matrices of the form
$\left(\matrix{1 & n\cr 0 & 1\cr}\right)$, $n\in {\Bbb Z}$.

\smallskip
{\bf 2.1.2} Let $g(m)=\left(\matrix{0 & 1\cr 1 & m\cr}\right)$, where $m\ge 1$.
Clearly $g(m)\in {\cal M}$.  Moreover, ${\cal M}$ is the {\it free monoid}
generated by the elements $g(m)$, $m\ge 1$: each element $g$ of ${\cal
M}$ can be written as
$
g = g(m_1)
\cdots g(m_r)\; , \;\;r\ge 0\; , \; m_{i}\ge 1\; ,
$
and this decomposition is unique (see Proposition A1.2).

\smallskip
{\bf 2.1.3}
One has
$$
G = Z \cdot {\cal M}\cdot H\; ,
$$
i.e. the application
$
Z\times {\cal M}\times H\rightarrow G\; , \;\;
(z,m,h)\mapsto g=z\cdot m\cdot h
$
is a {\it bijection}.

\smallskip
{\bf 2.1.4} The subset $Z\cdot {\cal M}$ of $G$ is made of matrices
$g=\left(\matrix{a & b\cr c & d\cr}\right)$ such that
$d\ge c\ge 0$ with the following additional restrictions:
$a=1$ if $c=0$, and, $b=a+1$ if $d=c=1$.

We will also often use the following remark, which is an immediate
consequence of the structure of ${\cal M}$ and of the relation
$$
g(m) \left(\matrix{ 1 & 1\cr 0 & 1\cr}\right) = g(m+1)\; :
$$
one has the {\it partition}
$$
Z\cdot {\cal M} = Z\cdot {\cal M}\cdot \left(\matrix{ 1 & 1\cr 0 &
1\cr}\right) \bigsqcup Z\cdot {\cal M}\cdot \left(\matrix{
0 & 1\cr 1 & 1\cr}\right) \; .
$$

\smallskip
{\bf 2.1.5} Let us consider the usual action of $G$ on
$\overline{\Bbb C}= {\Bbb C}\cup\{\infty\}$ by homographies:
$\left(\matrix{a & b \cr c & d \cr}\right) \cdot z = {\dst az+b\over
\dst cz+d}$. The following facts are easy to check:
\item{1.} $g\cdot [0,1] = [0,1]$ if and only if $g$ belongs to the
subgroup of order $4$ of matrices of the form
$\pm \left(\matrix{ 1 & 0\cr 0 & 1}\right)$, $\pm \left(\matrix{
-1 & 1\cr 0 & 1\cr}\right)$.
\item{2.} The monoid of the elements $g$ such that $g\cdot
[0,1]\subset [0,1]$ admits the partition
$$
{\cal M}\bigsqcup {\cal M}\left(\matrix{ -1 & 0\cr 0 & -1\cr}\right)
\bigsqcup {\cal M}\left(\matrix{ -1 & 1\cr 0 & 1\cr}\right)
\bigsqcup {\cal M}\left(\matrix{ 1 & -1\cr 0 & -1\cr}\right)\; .
$$
Note that $\left(\matrix{ -1 & 1\cr 0 & 1\cr}\right)=
\left(\matrix{ 1 & 1\cr 0 & 1\cr}\right)
\left(\matrix{ -1 & 0\cr 0 & 1\cr}\right)$.
\item{3.} The application $g\mapsto g\cdot 1 = {\dst a+b\over\dst c+d}$
is a {\it bijection} of $Z{\cal M}$ over ${\Bbb Q}$ which
maps ${\cal M}$ onto ${\Bbb Q}\cap (0,1]$.
\item{4.} The application $g\mapsto g\cdot 0 = b/d$ maps
$Z{\cal M}$ onto ${\Bbb Q}$ and each rational number has exactly {\it
two inverse images}. The two elements which map $0$ on $1$ are
$\left(\matrix{ 1 & 1\cr 0 & 1\cr}\right)$ and
$\left(\matrix{ 0 & 1\cr 1 & 1\cr}\right)$. This makes the partition
of 2.1.4 less mysterious.
\item{5.} The application $g\mapsto g\cdot [0,+\infty]$ is a
{\it bijection} of $Z{\cal M}$ on the set of Farey intervals
(the convention we adopt here implies that
$[n,+\infty]$ is a Farey interval, but $[-\infty, n]$ is not).
For the
definition and properties of the Farey partition of $[0,1]$
we refer the reader to
[HW].
\item{6.} The application $g\mapsto g\cdot\infty = a/c\in
\overline{\Bbb Q} = {\Bbb Q}\cup\{\infty\}$ of $Z{\cal M}$ on
$\overline{\Bbb Q}$ is surjective. Moreover
\itemitem{$\bullet$} $g\cdot\infty = \infty$ if and only if $g\in Z$;
\itemitem{$\bullet$} $g\cdot\infty = n\in {\Bbb Z}$ if and only if
$$
g=\left(\matrix{n & 1+kn\cr 1 & k\cr}\right) \; , \;\; k\ge 1\;
\hbox{or}\; g = \left(\matrix{n & -1+kn\cr 1 & k\cr}\right) \; , \;\;
k\ge 2\;  ;
$$
\itemitem{$\bullet$} $g\cdot\infty = a/c$, $c>1$ if and only if
$$
g=\left(\matrix{a & a'+ka\cr c & c'+kc\cr}\right) \; , \;\; k\ge 1\;
\hbox{or}\; g = \left(\matrix{a & a''+ka\cr c & c''+kc\cr}\right) \; , \;\;
k\ge 1\;  ,
$$
where $\left[{a'\over c'},{a''\over c''}\right]$ is the Farey interval
which contains $a/c$.

\mdk
\beginsection  2.2 Actions of ${\cal M}$ on
some spaces of holomorphic functions\par
\smk

{\bf 2.2.1} Let $U$ be an open subset of $\overline{\Bbb C}$. We will
denote by ${\cal O}(U)$ the complex vector space of holomorphic
functions on $U$. Let $I\subset \R$ be a compact interval and
let $k\in \Z$.
If $k\ge 0$ ($k<0$ respectively), we will denote
with ${\cal O}^{k}(\overline{\Bbb C}\setminus I)$
the complex vector space of
functions holomorphic in ${\Bbb C} \setminus I$,
meromorphic in $\overline{\Bbb C}\setminus I$,
which have a zero  at infinity  of order at least $k$
(resp. a pole of order at most $|k|$).

Let $g=\left(\matrix{a & b\cr c & d\cr}\right)\in G$ and assume
that $\varphi$ is meromorphic in $U$. We define
$$
L_{g}^{(k)}\varphi (z) = (a-cz)^{-k} \varphi \left({dz-b\over
a-cz}\right)\; . \eqno(2.1)
$$
The function $L_{g}^{(k)}\varphi$ is meromorphic in $g\cdot U$. Note that
$$
\eqalign{
(L_{g}^{(k)}\varphi )' &= kcL_{g}^{(k+1)} \varphi +
\varepsilon_{g}L_{g}^{(k+2)} \varphi'\cr
(L_{g}^{(k)}\varphi )'' &= k(k+1) c^{2}L_{g}^{(k+2)}\varphi + \varepsilon_{g}
(2k+2) c L_{g}^{(k+3)}\varphi' + L_{g}^{(k+4)} \varphi''\cr
}
$$
thus, for $k=-1$
$$
(L_{g}^{(-1)}\varphi)'' = L_{g}^{(3)}\varphi''\; . \eqno(2.2)
$$
Note also that if $g\in Z$ then $L_{g}^{(k)}$ does not depend on $k$.

The formula (2.1) above clearly defines an action of $G$: if $g,g'\in
G$, $k\in {\Bbb Z}$ and $\varphi$ is meromorphic in $U$ then the
functions $L_{g}^{(k)}(L_{g'}^{(k)}\varphi )$ and $L_{gg'}^{(k)}\varphi$
(meromorphic on $gg'\cdot U$) coincide.

\smallskip
{\bf 2.2.2} Let $J$ denote a compact interval of ${\Bbb R}$ and let
$\varphi\in {\cal O}^{2}(\overline{\Bbb C}\setminus J)$. The series
$$
\sum_{n\in {\Bbb Z}}\varphi (z-n) = \sum_{g\in Z}L_{g}^{(k)}\varphi
\; \; \hbox{(for all }k)\; , \eqno(2.3)
$$
converges uniformly on compact subsets of ${\Bbb C}\setminus {\Bbb
R}$ and also on the domains $\{z\in {\Bbb C}\; , \;
|\IM z|\ge \delta >0\; , \; \; |\RE z|\le A\}$. The sum will be
denoted $\sum_{Z}\varphi$; it is a function holomorphic in ${\Bbb
C}\setminus {\Bbb R}$, periodic of period $1$ and vanishing at
$\pm i\infty$. Thus taking the quotient by ${\Bbb Z}$ it can be
represented by means of the variable $q=e^{\pm 2\pi i z}$.

If $\varphi\in {\cal O}^{1}(\overline{\Bbb C}\setminus J)$
with $[0,1]\subset J$ then one
can decompose in a unique way
$$
\varphi (z) = a_{0}\log {z\over z-1}+\varphi_{0}(z)\; , \eqno(2.4)
$$
where $a_{0}\in {\Bbb C}$, $\varphi_{0}\in {\cal O}^{2}(\overline{\Bbb C}
\setminus J)$ and we consider the main branch of the logarithm
in $\C\setminus\R^-$. We have
$$
\sum_{n=-N}^{N}\log {z-n\over z-n-1} = \log {z+N\over z-N-1}
$$
and this leads to the definition
$$
\sum_{Z}\varphi (z):= \sum_{Z}\varphi_{0}(z) + \cases{
-a_{0}\pi i & if $\IM z >0$\cr
+a_{0}\pi i & if $\IM z<0$\cr}
\; .
\eqno(2.5)
$$
Note that in order to insure the convergence of the series
$\sum_{g\in Z}L_{g}^{(k)} \varphi$ one must regroup together the terms
$\left(\matrix{ 1 & n \cr 0 & 1\cr}\right)$ and
$\left(\matrix{ 1 & -n \cr 0 & 1\cr}\right)$, $n\ge 1$
(symmetric summation or Eisenstein's summation).

\smallskip
{\bf 2.2.3} Let $k\in {\Bbb Z}$ and $g\in {\cal M}$. Since $g\cdot
[0,1]\subset [0,1]$ if $\varphi$ is meromorphic in $\overline{\Bbb
C}\setminus [0,1]$ then $L_{g}^{(k)}\varphi$ will still be meromorphic
in $\overline{\Bbb C}\setminus [0,1]$. Moreover if $\varphi\in
{\cal O}^{k}(\overline{\Bbb C}\setminus [0,1])$ then also
$L_{g}^{(k)}\varphi\in {\cal O}^{k}(\overline{\Bbb C}\setminus [0,1])$. Thus
one has an action of the monoid ${\cal M}$ on
${\cal O}^{k}(\overline{\Bbb C}\setminus [0,1])$.

\smallskip
{\bf 2.2.4} We will now define a new action $L_{g}$ of ${\cal M}$ on
${\cal O}^{1}(\overline{\Bbb C}\setminus [0,1])$ which differs from
the action
$L_{g}^{(-1)}$ on $\O{-1}{[0,1]}$ by an affine correction.

All functions $\varphi\in {\cal O}^{-1}(\overline{\Bbb C}\setminus [0,1])$
can be uniquely written as
$$
\varphi (z) = Az+B+p(\varphi )(z)\; , \;\;p(\varphi )\in
{\cal O}^{1}(\overline{\Bbb C}\setminus [0,1])\; . \eqno(2.6)
$$
Note that if $g=\left(\matrix{ a & b \cr c & d \cr}\right)\in G$ and
$\psi (z)= Az+b$ then
$$
L_{g}^{(-1)}\psi (z) = A(dz-b)+B(a-cz)
$$
is still an affine function. Thus the formula
$$
L_{g}\varphi := p(L_{g}^{(-1)}\varphi)\; , \eqno(2.7)
$$
where $g\in {\cal M}$, $\varphi\in {\cal O}^{1}
(\overline{\Bbb C}\setminus [0,1])$, defines an {\it action} of
${\cal M}$ on ${\cal O}^{1}(\overline{\Bbb C}\setminus [0,1])$
which makes the following diagram commute
$$
\diagram{
{\cal O}^{-1}(\overline{\Bbb C}\setminus [0,1]) &
\hfl{L_{g}^{(-1)}}{} & {\cal O}^{-1}(\overline{\Bbb C}\setminus [0,1])
    \cr
\vfl{p}{}&&\vfl{}{p}
    \cr
{\cal O}^{1}(\overline{\Bbb C}\setminus [0,1]) &
\hfl{}{L_{g}} & {\cal O}^{1}(\overline{\Bbb C}\setminus [0,1])
    \cr}
$$
More explicitly, if $g=\left(\matrix{ a & b \cr c & d \cr}\right)\in
{\cal M}$
$$
L_{g}\varphi (z) =
(a-cz)\left[\varphi\left({dz-b\over a-cz}\right)-\varphi\left(-{d\over
c}\right)\right] - \varepsilon_{g}c^{-1}\varphi'\left(-{d\over
c}\right)\; . \eqno(2.8)
$$
Since this definition differs from that of $L_{g}^{(-1)}$ only by an
affine correction one clearly has
$$
(L_{g}\varphi)'' = L_{g}^{(3)}\varphi''\; , \eqno(2.9)
$$
where $\varphi''\in {\cal O}^{3}(\overline{\Bbb C}\setminus [0,1])$.
Equivalently one can say that $L_{g}\varphi$ is obtained by taking the
double primitive of
$L_{g}^{(3)}\varphi''$ which vanishes at infinity.

Two other  formulas  will be used throughout
what follows: if $\varphi\in {\cal O}^{1}(\overline{\Bbb C}\setminus
[0,1])$, $g=\left(\matrix{ a & b \cr c & d \cr}\right)\in
{\cal M}$ and $z\notin \left[{b\over d},{a\over c}\right]=
g\cdot [0,+\infty]$ one has
$$
\eqalign{
L_{g}\varphi (z) &=
\varepsilon_{g}c^{-1} \left[\int_{0}^{1}\varphi'\left(
-{d\over c}+\varepsilon_{g}{t\over c(a-cz)}\right)dt -
\varphi' \left(-{d\over c}\right)\right] \cr
&= c^{-2}(a-cz)^{-1}\int_{0}^{1}\varphi''\left(-
{d\over c}+\varepsilon_{g}{t\over c(a-cz)}\right)(1-t)dt\; . \cr}
\eqno(2.10)
$$
Note that
the assumption about $z$ means that the segment whose extremities
are $-d/c$ and $(dz-b)/(a-cz)$ does not intersect the
interval [0,1]. The two formulas are thus nothing else than
Taylor's formulas of first and second order  with integral remainder.

\smallskip
{\bf 2.2.5}
If $\varphi\in {\cal O}(\overline{\C}\setminus I)
$ we denote by
$
\sigma\cdot\varphi (z) = \overline{\varphi(\bar{z})}\; .
$
$\sigma$ is an involution of ${\cal O}(\overline{\C}\setminus I)$ which
preserves all the subspaces
${\cal O}^k(\overline{\C}\setminus I)$, commutes with $p$ and
with the actions of ${\cal M}$ on ${\cal O}^k$. %%%%%%%%
If $\sigma\cdot \varphi=\varphi$,
then
$\varphi (x)\in {\Bbb R}$ for all $x\in {\Bbb R}\setminus I$,
and we say that $\varphi$ is real holomorphic.
\par
%\vskip .5 truecm
%%%%%%%%%%%%  3. Complex continued fraction %%%%%%%%%%%%%%%%%
\beginsection 3. Complex continued fractions  \par
%\vskip .5 truecm

Exactly as in the real case treated in [MMY], where the use of
the continued fraction expansion was important for the study of the
real Brjuno function,
the introduction of a
complex analogue of Gauss' algorithm of continued fraction expansion
of a real number will be essential for the study of the boundary
behaviour of $\sum_{\cal M}(L_{g}\varphi) (z)$ and the construction of the
complex Brjuno function. In this section we first define our complex
version of the continued fraction (section 3.1) then we use it
to estimate the spectral radius of $T$ (section 3.2) and to prove the
convergence of the series (1.4) and (1.6) (Corollary 3.6).

%\vskip .5 truecm
\beginsection 3.1 Definition of the complex continued fractions \par
%\vskip .3 truecm

{\bf 3.1.1} We consider the following  domains:
$$
\eqalign{
D_{0} &= \left\{ z\in {\Bbb C}\, , |z+1|\le 1\, , \,
\RE z \ge {\sqrt{3}\over 2} -1\right\}\; , \cr
D_{1} &= \left\{ z\in {\Bbb C}\, , |z|\ge 1\, , \,
\left|z-{1\over\sqrt{3}}\right|\le {1\over\sqrt{3}}\right\}\; , \cr
D &= \{ z\in {\Bbb C}\, , |z|\le 1\, , \,
|z-i|\ge 1\, ,
\, |z+i|\ge 1\, , \, \RE z > 0\}\; , \cr
H_{0} &= \{ z\in {\Bbb C}\, , |z-i|\le 1\, , \, |z+1|\ge 1\, , \,
\IM z \le 1/2\}\; , \cr
H_{0}' &= \{ z\in {\Bbb C}\, , \, \overline{z} \in H_{0}\}\cr
\Delta &= H_{0}\cup H_{0}'\cup D
= \{z\in \C\, , |z|\le 1\, , \, |z+1|\ge 1\, , \, |\IM z|\le
1/2\}\; , \cr
D_{\infty} &= \overline{\C}\setminus (D_{0}\cup \Delta\cup D_{1})\cr
&= \{|\IM z|>1/2\}\cup\{\RE z <{\sqrt{3}\over 2}-1\}\cup
\{\RE z>{\sqrt{3}\over 2}\, , \, |z-1/\sqrt{3}|>1/\sqrt{3}\}\; . \cr
}
$$
Figures 2 and 3 show these domains and their image under the inversion
$S(z)=1/z$. A fundamental property is the following
\item{$\bullet$} if $z\notin D\cup D_{1}$ (in particular if $z\in
D_{\infty}$) then $1/z-m\in D_{\infty}$ for all $m\ge 1$;
\item{$\bullet$} if $z\in D_{1}$, then $1/z-1\in D_{0}$ and
$1/z-m\in D_{\infty}$ for all $m\ge 2$.

Observe that
$
SD = \cup_{m\ge 1}(\Delta +m)\; ,
$
where the domains have disjoint interior. Thus, for $z\in D$, we
define
$$
A(z) = {1\over z}- m = (g(m))^{-1}\cdot z\; , \eqno(3.1)
$$
(we recall that $g(m) = \left(\matrix{ 0 & 1 \cr 1 & m\cr }\right)$,
$m\ge 1$) where $m\ge 1$ is the unique integer such that
$
A(z) \in \Delta\; , \;\;\; |A(z)|<1\; .
$
Iterating from $z_{0}\in D$, we define
$$
z_{i+1} = A(z_{i})=A^{i+1}(z_{0})\;  \eqno(3.2)
$$
as long as $z_{i}=A^i (z) \in D$. The iteration process
stops when one of the two following conditions is verified:
\item{$\bullet$} $z_{l}=0$ for some $l\ge 0$;
this happens if and only if  $z_{0}\in \Q$,
\item{$\bullet$} $z_{l}\notin (D\cup\{0\})$, for some $l\ge 0$;
this happens if and only if $z_{0}\notin \R$.

For all $0\le i<l$, we will denote $m_{i+1}$ the integer such that
$$
z_{i+1} = {1\over z_{i}} - m_{i+1}\, , \; m_{i+1}\ge 1\; . \eqno(3.3)
$$

\smallskip
{\bf 3.1.2} Let
$
\left(\matrix{p_{i-1} & p_{i}\cr q_{i-1} & q_{i}\cr}\right)
=
\left(\matrix{0 & 1\cr 1 & m_{1}\cr}\right)
\ldots
\left(\matrix{0 & 1\cr 1 & m_{i}\cr}\right)
\in {\cal M} \, , \; 0\le i\le l \; .
$
Then one has the same recurrence relations as for the real continued
fraction
$$
\eqalign{
p_{i+1} &= m_{i+1}p_{i}+p_{i-1}\; , \cr
q_{i+1} &= m_{i+1}q_{i}+q_{i-1}\; , \cr
} \eqno(3.4)
$$
with initial data $p_{-1}= q_{0}=1$ and $p_{0}=q_{-1}=0$.  Moreover
$$
z_{0} = {p_{i-1}z_{i}+p_{i}\over q_{i-1}z_{i}+q_{i}}\, , \;\;\;
z_{i} = {p_{i}-q_{i}z_{0}\over q_{i-1}z_{0}-p_{i-1}}\, , \eqno(3.5)
$$
and if one poses
$$
\beta_{i}(z_{0}) = \prod_{j=0}^iz_{j}= (-1)^i(q_{i}z_{0}-p_{i})\; ,
\eqno(3.6)
$$
then
$$
\beta_{i}(z_{0}) = {z_{i}\over q_{i}+q_{i-1}z_{i}}= {1\over q_{i+1}+
q_{i}z_{i+1}}
\; . \eqno(3.7)
$$
Finally one has
$$
\eqalign{
(-1)^i\IM z_{0}
&= |\beta_{i-1}(z_{0})|^{2}\IM z_{i} = |q_{i}+q_{i-1}z_{i}|^{-2}\IM
z_{i}\;  , \cr
{dz_{i}\over dz_{0}}
&= (-1)^i(\beta_{i-1}(z_{0}))^{-2} =
(-1)^i(q_{i}+q_{i-1}z_{i})^{2}\; . \cr
}
$$
Observe that, as $|z_{i+1}+1|\ge 1$ and $\RE z_{i+1}\ge {\sqrt{3}\over
2}-1$ for $i<l$, we have from (3.7)
$$
|\beta_{i}(z_{0})|\le q_{i+1}^{-1}[\cos\pi/12]^{-1} =
{2\sqrt{2}\over 1+\sqrt{3}}q_{i+1}^{-1} \eqno(3.8a)
$$
and, as $q_{i}\le q_{i+1}$, $|z_{i+1}|\le 1$,
$$
|\beta_{i}(z_{0})|\ge {1\over 2}q_{i+1}^{-1}\; . \eqno(3.8b)
$$

%\vskip .5 truecm
\beginsection 3.2 The operator $T$, its spectral radius
and the sum over the monoid \par
%\vskip .3 truecm

{\bf 3.2.1}
Let $\gamma_{1}>\gamma_{0}>-1$, $I=[\gamma_{0},\gamma_{1}]$, $\varphi
\in \O{1}{I}$.

For all $m\ge 1$ one has (compare with (2.10))
$$
\eqalign{
L_{g(m)} \varphi (z) &= - z \left( \varphi \left({1\over z} -
m\right) - \varphi (-m)\right) +\varphi' (-m)\cr
&= - z^{-1} \int_{0}^{1}\varphi'' \left( -m+{t\over z}\right)
(1-t) dt\; , \cr
}\eqno(3.9)
$$
provided that $z\notin [0,1/(\gamma_{0}+m)]$,
and an even simpler formula for the action at the level of second
derivatives
$$
L^{(3)}_{g(m)} \psi (z) = -z^{-3}\psi\left( {1\over z}-m\right) \; .
\eqno(3.10)
$$

\smallskip
{\bf 3.2.2}
Let $\varepsilon >0$,
$$
U_{\varepsilon} = \{ z\in \C\, , \RE z\le \gamma_{0}-\varepsilon\,
\hbox{or}\, \RE z\ge \gamma_{1}+\varepsilon\,\hbox{or}\, |\IM z|\ge
\varepsilon\}\; . \eqno(3.11)
$$
We have the following
\smallskip
\Proc{Proposition 3.1}{Let $\gamma_{0},\gamma_{1}$ as above,
$I=[\gamma_{0},\gamma_{1}]$, $J=[0,1/(1+\gamma_{0})]$.
\item{1.} For all $\varphi\in {\cal O}^{1}(\overline{\Bbb C}\setminus
I)$, the series $\sum_{m\ge 1}L_{g(m)}\varphi $ converges
uniformly on compact subsets $K$ of $\overline{\Bbb C}\setminus
J$ to a function $T\varphi \in {\cal O}^{1}(\overline{\Bbb C}\setminus
J)$ and there exists $\varepsilon > 0$  and $C_{K}>0$ such that :
$
\sup_{K}|T\varphi |\le
C_{K}\sup_{U_{\varepsilon}}|\varphi |\; .
$
\item{2.} For all $\psi \in {\cal O}^{3}(\overline{\Bbb C}\setminus
I)$ the series $\sum_{m\ge 1}L_{g(m)}^{(3)}\psi $ converges
uniformly on compact subsets $K$ of $\overline{\Bbb C}\setminus
J$ to a function $T^{(3)}\psi \in {\cal O}^{1}(\overline{\Bbb C}\setminus
J)$ and there exists $\varepsilon > 0$  and $C_{K}>0$ such that :
$
\sup_{K}|T^{(3)}\psi |\le
C_{K}\sup_{U_{\varepsilon}}|\psi |\; .
$
\item{3.} For all $\varphi \in {\cal O}^{1}(\overline{\Bbb C}\setminus
I)$  one has }
$
T^{(3)} \varphi''  =  (T\varphi )''\; .
$

\smallskip
\proof
Let $\varepsilon >0$.
There exist $c_{\varepsilon}>0$ such that
for $\psi\in\O{3}{I}$, $z\in U_{\varepsilon}$ one has
$$
|\psi (z)| \le c_{\varepsilon}|z|^{-3}\sup_{U_{\varepsilon}}|\psi |\; .
$$
If $K$ is a compact subset of $\overline{\C}\setminus J$ there exists
$\varepsilon=\varepsilon_{K}$ such that $1/z-m\in U_{\varepsilon}$
for all $m\ge 1$ and $z\in K$. Moreover there exists $c_{K}>0$
and $M=M(\gamma_0)$, such that $|1/z-m|^{-1}\le c_{K}m^{-1}$
for all $z\in K$, $m\ge M$.
Consequently, for $z\in K$ and $\psi \in \O{3}{I}$ we have
$$
\sum_{m\ge 1} \left| \psi \left({1\over z}-m\right)\right|
\le c'_K \sup_{U_{\varepsilon}}|\psi | \quad,\quad\hbox{with}\quad
c'_K=M+c_{\varepsilon}c_{K}^{3}\sum_{m\ge M}m^{-3}\ ,
$$
which proves the second part of the proposition.

By integrating twice from $\infty$ one deduces that for
$\varphi\in \O{1}{I}$ the series
$\sum_{m\ge 1}L_{g(m)}\varphi$ converges uniformly on compact subsets
of $\overline{\C}\setminus J$ to a function $T\varphi\in \O{1}{J}$.
Moreover for any compact $K\subset \overline{\C}\setminus J$ there
exist $c_{\varepsilon}>0$ and $c_{K}>0$ such that for $z\in K$,
$
|(T\varphi )''(z) | = |(T^{(3)}\varphi'')(z)|\le
c_{K}|z|^{-3}\sup_{U_{\varepsilon}}|\varphi'' |\; .
$
On the other hand there exists $c''_{\varepsilon}>0$ such that,
$
\sup_{U_{\varepsilon}}|\varphi'' |\le
c''_{\varepsilon}\sup_{U_{\varepsilon /2}}|\varphi |
$
(by Cauchy's formula), hence we get for $z\in K$,
$
|T\varphi (z)|\le \tilde c_{K}|z|^{-1}\sup_{U_{\varepsilon /2}}|\varphi |
\; .
$
The third part of the
Proposition is immediate.
\qed
\smallskip

\smallskip
{\bf 3.2.3} The open set
$\overline{\Bbb C}\setminus [0,1]$
is an hyperbolic Riemann surface which is naturally
equipped with a Poincar\'e metric. The following well-known fact will be
crucial for the proof of Lemma 3.2: given two hyperbolic Riemann
surfaces $M,N$ and an analytic map $f\,:\, M\rightarrow N$ either its
differential $df$ contracts the hyperbolic metric or $f$ is a
surjective local isometry. In what follows we will denote $d_{hyper}$
the Poincar\'e metric on the Riemann surface under consideration.

Given $\rho >0$ we denote
$$
V_{\rho}(D_{\infty}) = \{z\in\overline{\Bbb C}\setminus [0,1]\, , \,
d_{hyper}(z,D_{\infty})<\rho\}\; . \eqno(3.12)
$$
the $\rho$--neighborhood of $D_{\infty}$ in $\overline{\Bbb C}\setminus
[0,1]$.

\smallskip
\Proc{Lemma 3.2}{Let $\rho \ge 0$. For all $m\ge 1$ and $z\in
V_{\rho}(D_{\infty})$ one has ${1\over z}-m\in V_{\rho}(D_{\infty})$}

\smallskip
\proof
The M\"obius transformation
$z\mapsto {1\over z}-m$ maps $D_{\infty}$ into itself and
$\overline{\Bbb C}\setminus [0,1]$ onto $\overline{\Bbb C}\setminus
[1-m,+\infty ]$ which is contained in $\overline{\Bbb C}\setminus [0,1]$.
Thus it decreases the hyperbolic distance and the
$\rho$--neighborhood of the image of $D_{\infty}$ w.r.t.\ the
Poincar\'e metric of $\overline{\Bbb C}\setminus [1-m,+\infty]$ is
contained into the $\rho$--neighborhood of $D_{\infty}$ w.r.t.\ the
Poincar\'e metric of $\overline{\Bbb C}\setminus [0,1]$. \qed
\smallskip

\smallskip
{\bf 3.2.4}
Once the existence of the operator $T$ on holomorphic functions is
established (Proposition 3.1)
one can ask for more information on Hardy spaces. The
results we will prove are completely analogous to those obtained for
the real Brjuno operator in [MMY].

\Proc{Proposition 3.3}{Let $\rho\ge 0$. There exists $c_{\rho}'>0$
such that for all $r\ge 0$ and $\psi \in {\cal O}^{3}(\overline{\Bbb
C}\setminus
[0,1])$ one has }
$$
\sup_{V_{\rho}(D_{\infty})}|((T^{(3)})^r\psi )(z)|\le
c_{\rho}' \left({\sqrt{5}-1\over 2}\right)^r\sup_{V_{\rho}(D_{\infty})}
|\psi (z)|\; .
$$

\smallskip
\proof
First of all note that if $z\in D_{\infty}$ and $g\in \cM$ then
$g^{-1}\cdot z\in D_{\infty}$, thus taking into account  Lemma 3.2,
if $z\in V_{\rho}(D_{\infty})$ then $g^{-1}z\in
V_{\rho}(D_{\infty})$.

Given an integer  $r\ge 0$ we denote ${\cal M}^{(r)}$ the set of
elements $g$ of ${\cal M}$ of the form $g(m_{1})\ldots g(m_{r})$,
$m_{i}\ge 1$ for all $1\le i\le r$. Then one has
$$
(T^{(3)})^r\psi = \sum_{g\in {\cal M}^{(r)}}L_{g}^{(3)}\psi
= \sum_{g'\in {\cal M}^{(r-1)}}L_{g'}^{(3)} (T^{(3)}\psi )\; .
$$
Let
$g' = \left(\matrix{ a' & b' \cr c' & d' \cr}\right) \in {\cal
M}^{(r-1)}$ and $z\in V_{\rho}(D_{\infty})$; let
$z' = {d'z-b'\over a'-c'z}$. For all $m\ge 1$ one has
$$
\left(\matrix{ a' & b' \cr c' & d' \cr}\right)
\left(\matrix{ 0 & 1 \cr 1 & m \cr}\right) =
\left(\matrix{ b' & d' \cr a'+mb' & c'+md' \cr}\right)\; ,
$$
from which it follows that
$$
L_{g'}^{(3)} (T^{(3)}\psi ) (z) =
(b'-d'z)^{-3} \sum_{m\ge 1}\psi \left({1\over z'}-m\right) \;  .
$$
Since $z'\in V_{\rho}(D_{\infty})$, as we have seen during the proof
of  Proposition 3.1 one has
$$
\left| \sum_{m\ge 1}\psi \left({1\over z'}-m\right)\right| \le
c_{\rho }\sup_{V_{\rho}(D_{\infty})}|\psi |\; .
$$
On the other hand one has
$|z-b'/d'|^{-1}\le c_{\rho}$ for all $z\in V_{\rho}(D_{\infty})$
and $b'/d'\in [0,1]$. Thus we get
$
|L_{g'}^{(3)} (T^{(3)}\psi ) (z) |\le \tilde{c}_{\rho}(d')^{-3}\; .
$
But now it is enough to recall that (see appendix A1)
$
\min_{{\cal M}^{(r-1)}} d' \ge C \left({\sqrt{5}+1\over
2}\right)^{r-1}\; ,
$
and
$
\sum_{{\cal M}^{(r-1)}} d'^{-2} \le C \; ,
$
to obtain the desired estimate.
\qed
\smallskip
\remark{3.4} In a completely analogous way we may prove that
for $\rho\ge 0$ and all  $\varphi\in {\cal O}^1(\overline{{\Bbb C}}\setminus
[0,1])$,
$$
\sup_{V_{\rho}(D_{\infty})}|T^r\varphi (z)| \le c_{\rho}'
\left({\sqrt{5}-1\over
2}\right)^{r} \sup_{V_{\rho}(D_{\infty})} |\varphi (z)| \; .
$$

\smallskip
\remark{3.5} We may also consider the Hardy space $H^p
(D_{\infty})$, $1\le p<+\infty$ of analytic functions
$\varphi \, : \, D_{\infty}\rightarrow {\Bbb C} $ such that
the subharmonic function $|\varphi |^p$ has a harmonic majorant. It
is an immediate consequence of the Riemann mapping theorem that
this space is isomorphic to $H^p({\Bbb D})$. Indeed if $h$
maps $D_{\infty}$ conformally onto $\D$ one can use the
norm
$$
\Vert \varphi\Vert_{H^p(D_{\infty})}=\Vert \varphi\circ h^{-1}
\Vert_{H^p(\D)} = \left(
\int_{\partial D_{\infty}} |\varphi (z)|^p |h'(z)||dz|\right)^{1/p}
\; .
$$
Note that since $\partial D_{\infty}$ is a
rectifiable Jordan curve $h$ extends to a homeomorphism of $\partial
D_{\infty}$ onto $\S^{1}$ which is conformal almost everywhere.
It is immediate to check that
the proof of Proposition 3.3 can be easily adapted so as to show that
$T$ is bounded linear operator on $H^p(D_{\infty})$ with spectral
radius $\le {\sqrt{5}-1\over
2}$.

\smallskip
{\bf 3.2.4}
We have now the following important Corollary, which establishes the
convergence of the series (1.4).

\smallskip
\Proc{Corollary 3.6}{}
\item{1.} Let $\psi \in {\cal O}^{3}(\overline{\Bbb C}\setminus
[0,1])$.
The family $(L_{g}^{(3)}\psi )_{g\in {\cal M}}$ is summable,
uniformly on compact subsets of $\overline{\Bbb C}\setminus
[0,1]$.
Its sum is equal to $\sum_{r\ge 0}(T^{(3)})^r\psi $, will
be denoted $\sum^{(3)}_{\cal M}\psi$ and for all compact subset $K$
of $\overline{\Bbb C}\setminus
[0,1]$ there
exists $\varepsilon >0$ such that
$$
\sup_{K}|\sum\nolimits^{(3)}_{\cal M}\psi |\le C_{K}
\sup_{U_{\varepsilon}} |\psi |\; .
$$
The family $(L_{g}^{(3)}\psi )_{g\in Z\cdot {\cal M}}$
is summable, uniformly on all domains of the form
$\{|\RE z|<A\, , \; |\IM z|\ge \delta\} $ (where $A$ and $\delta$ are
positive). Its sum is equal to $\sum_{\Bbb Z}
\sum^{(3)}_{\cal M}\psi $, will
be denoted $\sum^{(3)}_{Z\cdot \cal M}\psi$. It is holomorphic in
${\Bbb C}\setminus {\Bbb R}$, periodic of period $1$ and bounded in a
neighborhood of $\pm i \infty$.
\item{2.} Let $\varphi \in {\cal O}^{1}(\overline{\Bbb C}\setminus
[0,1])$. The family $(L_{g}\varphi )_{g\in {\cal M}}$ is summable,
uniformly on compact subsets of $\overline{\Bbb C}\setminus
[0,1]$.
Its sum is equal to $\sum_{r\ge 0}T^r\varphi $, will
be denoted $\sum_{\cal M}\varphi $. The function
$\sum_{\Bbb Z}(\sum_{\cal M}\varphi )$ will be  denoted
$\sum_{Z\cdot \cal M}\varphi$. It is holomorphic in
${\Bbb C}\setminus {\Bbb R}$, periodic of period $1$ and vanishes at
$\pm i \infty$. One has
$$
%\eqalign{
 (\sum_{\cal M}\varphi )'' = \sum\nolimits^{(3)}_{\cal M} \varphi''
\quad ,\quad
 (\sum_{Z\cdot \cal M}\varphi )'' =
 \sum\nolimits^{(3)}_{Z\cdot \cal M} \varphi''
\; .
%}
$$
\item{3.} $\sum_{\cM}$ (resp. $\sum_{\cM}^{(3)}$) and
$(1-T)$ (resp. $(1-T^{(3)})$) acting on
$\O{1}{[0,1]}$ (resp.
$\O{3}{[0,1]}$) are the inverses of one another:
$$
\eqalign{
(1-T)\sum_{\cM} &= \sum_{\cM}(1-T)= \hbox{id.}\cr
(1-T^{(3)})\sum\nolimits^{(3)}_{\cal M} &=
\sum\nolimits^{(3)}_{\cal M} (1-T^{(3)}) = \hbox{id.}\cr
}
$$\par
\proof
\item {1)} The only non trivial assertions are the summability
of the families $(L^{(3)}_g\psi)_{g\in {\cal M}}$
and $(L^{(3)}_g\psi)_{g\in Z{\cal M}}$. Writing $g=g'g(m)$,
$(m\ge1, g'\in {\cal M})$ for $g\in{\cal M}$, $g\neq {\rm id}$,
we have
$$
L^{(3)}_g\psi(z) = (b'-d'z)^{-3}\psi\left({1\over z'}-m\right)\ ,
$$
with
$$ g'=\left(\matrix{a'&b'\cr c'&d'\cr}\right)\quad,\quad
z'={d'z-b'\over a'-c'z}\ .
$$
Now for $z\in V_{\rho}(D_{\infty})$,
$
\left|L^{(3)}_g\psi(z)\right| = c_{\rho}
\left(md' {\rm dist\big(z,[0,1]\big)
}\right)^{-3}\sup_{V_{\rho}(D_{\infty})}|\psi|\ ,
$
and the summability assertions follow (see proof of
Proposition 3.3).
\par
\item {2)} Again, the only non trivial assertion is the summability one,
which is obtained from the first part by integrating twice.\par
\item {3)} The third part of the corollary is immediate.\par
\qed
%\vskip .5 truecm
%%%%%%%%%%%%%%%%%%%%   4. Boundary behaviour %%%%%%%%%%%%%%%
\beginsection 4. Boundary behaviour of $\sum_{\cal M}\varphi$  \par
%\vskip .5 truecm

In this section we will study the behaviour of $\sum_{\cal M}\varphi
(z)$ when $z$ is close to $[0,1]$. Our main tool for this study will
be the complex continued fraction introduced in the previous section.
In this section and in the next one we will for shortness
often denote $c$ or $C$ various
positive universal constants.

%\vskip .5 truecm
\beginsection 4.1 Decomposition into principal and residual terms \par
%\vskip .3 truecm

{\bf 4.1.1}
We begin our study of the boundary behaviour of
$\sum_{\cal M}\varphi
(z)$ considering the case $z$ is close to $0$.

\smallskip
\Proc{Proposition 4.1}{
\item{1.} Let $I=[\gamma_{0},\gamma_{1}]$, $\gamma_{0}>
-1$. There exists $c=c_{I}>0$ such that for all
$\varphi \in {\cal O}^{1}(\overline{\Bbb C}\setminus
I)$ and for all $z\in D_{0}\cup H_{0}\cup H_{0}'$ one has
$$
|T\varphi (z) -\sum_{m\ge 1} \varphi'(-m)|\le c_{I}|z|\log
(1+|z|^{-1})\sup_{U}|\varphi |\; , \eqno(4.1)
$$
where $U=\{z\in\C\, , |\IM z|\ge 1/2\, \hbox{or}\,
\RE z\le \gamma_{0}-1\,\hbox{or}\,\RE z\ge \gamma_{1}+1\}$.
\item{2.} There exists $c>0$ such that for all
$\varphi\in\O{1}{[0,1]}$ and
for all $z\in D_{0}\cup H_{0}\cup H_{0}'$ one has}
$$
|\sum_{\cal M}\varphi (z) - \varphi (z) - \sum_{m\ge 1} (\sum_{\cal
M}\varphi )'(-m) | \le  c |z|(1+\log
|z|^{-1})\sup_{D_{\infty}}|\varphi |\; . \eqno(4.2)
$$

\smallskip
\proof
We will only prove (4.2), since the proof of (4.1)
is essentially the same.

Let $\tilde{\varphi} = \sum_{\cal M}\varphi $. One has
$$
\tilde{\varphi} = \varphi + T\tilde{\varphi} =
\varphi + \sum_{m\ge 1}L_{g(m)}\tilde{\varphi}\; .
$$
Let $z\in D_{0}\sqcup H_{0}\sqcup H_{0}'$. We will consider the cases
$m\ge 3|z|^{-1}$ and $m< 3|z|^{-1}$ separately.

If $m \ge 3|z|^{-1}$ the segment $[-m,-m+{1\over z}]$ is contained in
the closed half--plane $\{\RE w \le -2/3\}$ and one has
$$
L_{g(m)}\tilde{\varphi} (z) = -{1\over z} \int_{0}^{1}(1-t)\tilde{\varphi} ''
(-m+{t\over z})dt\; .
$$
Applying Cauchy's estimate, it follows that
$
|L_{g(m)}\tilde{\varphi} (z) | \le c|z|^{-1}m^{-3}\sup_{D_{\infty}}
|\tilde{\varphi }|\; .
$
Since by Remark 3.4 $\sup_{D_{\infty}}
|\tilde{\varphi }|\le c \sup_{D_{\infty}}
|\varphi |$ one gets
$$
\left|\sum_{m\ge 3|z|^{-1}}L_{g(m)}\tilde{\varphi} (z)\right| \le
c |z| \sup_{D_{\infty}}|\varphi |\; .
$$
In the case $m< 3|z|^{-1}$ we separate the three terms constituting
$L_{g(m)}\tilde{\varphi }(z)$ and we obtain
$$
\eqalign{
|\tilde{\varphi}'(-m)| \le c m^{-2}
\sup_{D_{\infty}}|\tilde{\varphi}|\quad  &,\quad
|z \tilde{\varphi}(-m)| \le c m^{-1} |z|
\sup_{D_{\infty}}|\tilde{\varphi}|\;  , \cr
\left|z \tilde{\varphi}\left(-m+{1\over z}\right)\right| &\le c |-m
+{1\over z}|^{-1} \sup_{D_{\infty}}|\tilde{\varphi}|\;  , \cr
}
$$
thus
$$
\left|\sum_{m< 3|z|^{-1}}L_{g(m)}\tilde{\varphi} (z) - \sum_{m\ge 1}
\tilde{\varphi}'(-m)\right| \le
c |z| (1+\log |z|^{-1}) \sup_{D_{\infty}}|\varphi |\; ,
$$
and the assertion is proved.
\qed

\smallskip
{\bf 4.1.2}
Next we consider the behaviour near $z=1$:

\Proc{Proposition 4.2}{There exists $c>0$ such that for all
$\varphi\in {\cal O}^{1}(\overline{\Bbb C}\setminus
[0,1])$ and $z\in D_{1}$ one has}
$$
\eqalignno{
\left| T\varphi (z) +z\varphi\left({1\over z}-1\right)-
\sum_{m\ge 1}\varphi' (-m)\right| &\le c|z-1|\sup_{D_{\infty}}
|\varphi |\; , &(4.3) \cr
\left|\sum_{\cal M}\varphi (z) - \varphi (z) + z\varphi\left({1\over
z}-1\right)\right| &\le c |z-1| (1+\log |z-1|^{-1})\sup_{D_{\infty}}|\varphi|
\; . \qquad &(4.4) \cr
}
$$

\smallskip
\proof
For $z\in D_{1}$ we have ${1\over z}-1\in D_{0}$.  Moreover,
there exists $c>0$ such that for $m\ge 2$ and $w$ on the segment
with endpoints $1-m$ and $1/z-m$
we have
$$
|\varphi' (-w)|\le cm^{-2}\sup_{D_{\infty}}|\varphi |\; . \eqno(4.5)
$$
Since
$$
T\varphi (z) +z\varphi\left({1\over z}-1\right)-
\sum_{m\ge 1}\varphi' (-m) = z\sum_{m\ge 2}\left[
\varphi (1-m)-\varphi\left({1\over z}-m\right)\right]
$$
from (4.5) one deduces easily (4.3).
\smk
Let $\tilde{\varphi} = \sum_{\cal M}\varphi $. We have
$
\tilde{\varphi} (z) = \varphi (z) +T\tilde{\varphi}(z)\; ,
$
and
$$
\eqalign{
\tilde{\varphi}(z)  &=  \varphi (z) - z\varphi\left({1\over
z}-1\right) +T\tilde{\varphi}(z) +z\tilde{\varphi}\left({1\over
z}-1\right)-\sum_{m\ge 1}\tilde{\varphi}'(-m)\cr
&+z\left[\varphi\left({1\over
z}-1\right)-\tilde{\varphi}\left({1\over
z}-1\right)+\sum_{m\ge 1}\tilde{\varphi}'(-m)\right]
+(1-z)\sum_{m\ge 1}\tilde{\varphi}'(-m)\; . \cr}
$$
By (4.3) one has
$$
\left| T\tilde{\varphi} (z)+z\tilde{\varphi}\left({1\over z}-1\right)
-\sum_{m\ge 1}\tilde{\varphi}'(-m)\right|\le c|z-1|\sup_{D_{\infty}}
|\tilde{\varphi}|\; .
$$
We also have :
$
\sup_{D_{\infty}}|\tilde{\varphi}|\le c\sup_{D_{\infty}}|\varphi |
$
and by (4.2)
$$
\left|{\varphi} ({1\over z}-1) - \tilde\varphi ({1\over z}-1)
+ \sum_{m\ge 1} \tilde{\varphi}' (-m)\right| \le c |z-1|(1+\log |z-1|^{-1})
\sup_{D_{\infty}}|\varphi |\; .
$$
As
$
\left|\sum_{m\ge 1} \tilde{\varphi}' (-m)\right|  \le c
\sup_{D_{\infty}}|\tilde{\varphi}|\le c
\sup_{D_{\infty}}|\varphi |\ ,
$
we get the second inequality. \qed
\smallskip

\remark{4.3} It is easy to check that the estimates in Proposition 4.2
are valid if $z$ is such that $z-1\in \Delta$.

%\vskip .5 truecm

\beginsection 4.2 Boundary behaviour and continued fraction\par
%\vskip .3 truecm

Let $k\ge 1$ and $m_{1},\ldots ,m_{k}$ be integers $\ge 1$. We denote
by $D(m_{1},\ldots ,m_{k})$ the set of $z_{0}\in D$ for which the
complex continued fraction is
$$
z_{i}^{-1} = m_{i+1}+z_{i+1}\; , \;\;\; 0\le i<k\; ,
$$
with $z_{i}\in D$ for $0\le i<k$ and $z_{k}\in \Delta$. In the
following we set for $0<i\le k$
$$
\varepsilon_{i} = \cases{
0 & if $m_{i}=1\; ,$ \cr
1 & if $m_{i}>1\; .$\cr}
$$

\smallskip
\Proc{Proposition 4.4}{For $\varphi\in\O{1}{[0,1]}$ we have in
$D(m_{1},\ldots ,m_{k})$
$$
\eqalign{
T^k \varphi (z_{0}) &= (p_{k-1}-q_{k-1}z_{0})[\varphi (z_{k})+
\varphi (z_{k}-1)+\varepsilon_{k}\varphi (z_{k}+1)]\cr
&- (p_{k-2}-q_{k-2}z_{0})(1+z_{k-1})\varepsilon_{k-1}\varphi\left(
-{z_{k-1}\over 1+z_{k-1}}\right)+R^{[k]}(\varphi )(z_{0})\; .\cr
}\eqno(4.6)
$$
The remainder term $R^{[k]}(\varphi )$ is holomorphic in
$\hbox{int}\,D(m_{1},\ldots ,m_{k})$, continuous in
$D(m_{1},\ldots ,m_{k})$ and satisfies there }
$$
|R^{[k]}(\varphi )(z_{0})|\le ck\left({\sqrt{5}+1\over 2}\right)^{-k}
\sup_{D_{\infty}}|\varphi |\; . \eqno(4.7)
$$

\smallskip
\proof
One has
$$
T\varphi (z_{0}) = -z_{0}[\varphi (z_{1})+
\varphi (z_{1}-1)+\varepsilon_{1}\varphi (z_{1}+1)]+
R^{(m_{1})}(\varphi )(z_{0})\; , \eqno(4.8)
$$
with
$$
\eqalign{
R^{(m_{1})}(\varphi )(z_{0}) &=
\sum_{m\ge 1}\varphi' (-m)-z_{0}\sum_{m\ge 1\, , \, |m-m_{1}|> 1}
[\varphi (z_{1}+m_{1}-m)-\varphi (-m)]\cr
&+ z_{0}\sum_{m\ge 1\, , \, |m-m_{1}|\le  1}
\varphi (-m)\; . \cr
}
$$
For $z_{0}\in D(m_{1})$, one easily checks that :
$
|R^{(m_{1})}(\varphi )(z_{0})|\le c\sup_{D_{\infty}}|\varphi |
$
and that $R^{(m_{1})}(\varphi )$ is holomorphic in the neighborhood
of $D(m_{1})$.

Iterating (4.8) $k$ times we have
$$
\eqalign{
T^k\varphi (z_{0}) &= \prod_{i=0}^{k-1}(-z_{i})[\varphi (z_{k})+
\varphi (z_{k}-1)+\varepsilon_{k}\varphi (z_{k}+1)]\cr
&+\sum_{j=1}^{k-1}[\prod_{i=0}^{j-1}(-z_{i})](T^{k-j}\varphi
(z_{j}-1)+\varepsilon_{j}T^{k-j}\varphi (z_{j}+1))\cr
&+\sum_{j=1}^{k}[\prod_{i=0}^{j-2}(-z_{i})]
R^{(m_{j})}(T^{k-j}\varphi )(z_{j-1})\; . \cr
}\eqno(4.9)
$$
We have here :
$
\prod_{0}^{l-1}(-z_{i})=p_{l-1}-q_{l-1}z_{0}\; ,
$
and,
$
|p_{l-1}-q_{l-1}z_{0}|\le cq_{l}^{-1}\; .
$
The function
$$
R_{0}^{[k]}(z_{0}) = \sum_{j=1}^k [\prod_{i=0}^{j-2}(-z_{i})]
R^{(m_{j})}(T^{k-j}\varphi )(z_{j-1}) \eqno(4.10)
$$
is holomorphic in a neighborhood of $D(m_{1},\ldots ,m_{k})$ and
satisfies there
$$
|R_{0}^{[k]}(z_{0})|\le c\sum_{j=1}^k q_{j-1}^{-1}\sup_{D_{\infty}}
|T^{k-j}\varphi |\le ck\left({\sqrt{5}+1\over 2}\right)^{-k}
\sup_{D_{\infty}}|\varphi |\; .
$$
The function
$$
R_{1}^{[k]}(z_{0})= \sum_{j=1}^{k-1}(p_{j-1}-q_{j-1}z_{0})T^{k-j}\varphi
(z_{j}-1) \eqno(4.11)
$$
is holomorphic in the interior of $D(m_{1},\ldots ,m_{k})$
(and even in a neighborhood of $D(m_{1},\ldots ,m_{k})$ if $m_{k}>1$).
Applying Proposition 4.1 $(k-j)$ times to each term of the sum
we see that it is continuous in
$D(m_{1},\ldots ,m_{k})$ and satisfies there
$$
|R_{1}^{[k]}(z_{0})|\le c\sum_{j=1}^{k-1}q_{j}^{-1}\sup_{D_{\infty}}
|T^{k-j-1}\varphi |\le ck\left({\sqrt{5}+1\over 2}\right)^{-k}
\sup_{D_{\infty}}|\varphi |\; .
$$
The function
$$
R_{2}^{[k]}(z_{0})= \sum_{j=1}^{k-2}(p_{j-1}-q_{j-1}z_{0})\varepsilon_{j}
T^{k-j}\varphi
(z_{j}+1) \eqno(4.12)
$$
is holomorphic in a neighborhood of $D(m_{1},\ldots ,m_{k})$
and satisfies, according to Propositions 4.1 and 4.2
$$
|R_{2}^{[k]}(z_{0})|\le c\sum_{j=1}^{k-2}q_{j}^{-1}\sup_{D_{\infty}}
|T^{k-j-2}\varphi |\le ck\left({\sqrt{5}+1\over 2}\right)^{-k}
\sup_{D_{\infty}}|\varphi |\; .
$$
Finally, we apply Proposition 4.2 to
$
(p_{k-2}-q_{k-2}z_{0})\varepsilon_{k-1}T\varphi (z_{k-1}+1)
$
and summing up the different contributions to (4.9)
we get the desired result. \qed
%\vskip .5 truecm

\beginsection 4.3 $H^p$--estimates \par
%\vskip .3 truecm
{\bf 4.3.1} Let $1<p<+\infty$. We consider the space
$H^p(\H^\pm)$ of functions  $F\in\O{1}{[0,1]}$
whose restrictions to both $\H^{+}$ and $\H^-$ belong to
the Hardy space $H^p$ of these half--spaces, endowed with the norm
$$
\Vert F\Vert_{H^p(\H^\pm)}= \Vert F\mid_{\H^+}\Vert_{H^p}+
\Vert F\mid_{\H^-}\Vert_{H^p}\; . \eqno(4.13)
$$
See references [Du, St] for details. Then, it is a classical result
[Ga] that $F\in H^p(\H^\pm)$ if and
only if the associated hyperfunction $u$ belongs to $L^p([0,1])$ and
that the correspondence is an isomorphism of
Banach spaces.

But we know [MMY] that $T$ acting on $L^p([0,1])$ has spectral radius
$\le {\sqrt{5}-1\over 2}$. Consequently the same is true for
$T$ acting on $H^p (\H^\pm )$.

\smallskip
{\bf 4.3.2} For the case $p=\infty$ and on the larger
domain $\overline{\C}\setminus [0,1]$ we have the following

\smallskip
\Proc{Proposition 4.5}{The restriction of $T$ to
$H^\infty (\overline{\C}\setminus [0,1])\cap \O{1}{[0,1]}$ is a
bounded operator on this space with spectral radius $\le {\sqrt{5}-1\over 2}$.}

\smallskip
\proof
Let $k\ge 1$, $\varphi\in H^\infty (\overline{\C}\setminus [0,1])
\cap \O{1}{[0,1]}$ and $z_{0}\in \overline{\C}\setminus [0,1]$. We
estimate $|T^k\varphi (z_{0})|$ in various cases.
\item{(i)} If $z_{0}\in D_{\infty}$ we have
$$
|T^k\varphi (z_{0})|\le c\left({\sqrt{5}+1\over 2}\right)^{-k}
\sup_{D_{\infty}}|\varphi |\; ,
$$
according to Remark 3.4.
\item{(ii)} If $z_{0}\in D_{0}\cup H_{0}\cup H_{0}'$ we have,
according to Proposition 4.1 and to Remark 3.4
$$
|T^k\varphi (z_{0})|\le c\sup_{D_{\infty}}|T^{k-1}\varphi|\le
c\left({\sqrt{5}+1\over 2}\right)^{-k}
\sup_{D_{\infty}}|\varphi |\; .
$$
\item{(iii)} if $z_{0}\in D_{1}$ we have, according to Proposition 4.2
$$
|T^k\varphi (z_{0})| \le c\left[\sup_{D_{\infty}}|T^{k-1}\varphi|+
\left|T^{k-1}\varphi\left({1\over z_{0}}-1\right)\right|\right]
$$
which gives, for $k=1$ :
$
|T\varphi (z_{0})|\le c\Vert\varphi\Vert_{H^\infty}
$
and for $k>1$, according to Proposition 4.1 :
$$
|T^k\varphi (z_{0})| \le c\sup_{D_{\infty}}|T^{k-2}\varphi|\le
c\left({\sqrt{5}+1\over 2}\right)^{-k}
\sup_{D_{\infty}}|\varphi |\; .
$$
\item{(iv)} If $z_{0}\in D$ has continued fraction
$z_{i}^{-1} = m_{i+1}+z_{i+1}$ with $0\le i<l$, $z_{l}\in
\Delta\setminus D$, $l<k$, we apply Proposition 4.4 to write
$$
\eqalign{
T^k \varphi &(z_{0}) = (p_{l-1}-q_{l-1}z_{0})[T^{k-l}\varphi (z_{l})
+T^{k-l}\varphi (z_{l}-1)+\varepsilon_{l}T^{k-l}\varphi (z_{l}+1)]\cr
&-(p_{l-2}-q_{l-2}z_{0})(1+z_{l-1})\varepsilon_{l-1}T^{k-l}\varphi
\left({-z_{l-1}\over
1+z_{l-1}}\right)
+R^{[l]}(T^{k-l}\varphi)(z_{0})\; . \cr
}
$$
Here, we have $z_{l}$, $z_{l}-1$ and ${\dst -z_{l-1}\over\dst
1+z_{l-1}}$ in $D_{\infty}\cup D_{0}\cup H_{0}\cup H_{0}'$ hence the
value of $T^{k-l}\varphi$ at these points is in absolute value less
than $c\sup_{D_{\infty}}|T^{k-l-1}\varphi |$ (Proposition 4.1).
We also have, according to Proposition 4.2 and Proposition 4.1
$$
| T^{k-l}\varphi(z_{l}+1)|\le \cases{
c\Vert\varphi\Vert_{H^\infty} &if $k=l+1\; , $\cr
c\sup_{D_{\infty}}|T^{k-l-2}\varphi| &if $k>l+1\; . $\cr
}
$$
On the other hand we have
$$
\eqalign{
|p_{l-2}-q_{l-2}z_{0}| &\le cq_{l-1}^{-1}\le c
\left({\sqrt{5}+1\over 2}\right)^{-l}\; , \cr
|p_{l-1}-q_{l-1}z_{0}| &\le cq_{l}^{-1}\le c
\left({\sqrt{5}+1\over 2}\right)^{-l}\; , \cr
}
$$
and thus from the estimate of $R^{[l]}$ in Proposition 4.4
and from Remark 3.4  we have
$$
|T^{k}\varphi(z_{0})|\le \cases{
c\,l\left({\sqrt{5}+1\over 2}\right)^{-l}
\Vert\varphi\Vert_{H^\infty} &if $k=l+1\; , $\cr
c\,k\left({\sqrt{5}+1\over 2}\right)^{-k}
\sup_{D_{\infty}}|\varphi| &if $k>l+1\; . $\cr
}
$$
\item{(v)} If $z_{0}\in D$ has continued fraction
$z_{i}^{-1} = m_{i+1}+z_{i+1}$ with $0\le i<k$, $z_{k}\in
\Delta$, we again apply Proposition 4.4 to get
$$
\eqalign{
T^k \varphi (z_{0}) &= (p_{k-1}-q_{k-1}z_{0})[\varphi (z_{k})
+\varphi (z_{k}-1)+\varepsilon_{k}\varphi (z_{k}+1)]\cr
&-(p_{k-2}-q_{k-2}z_{0})(1+z_{k-1})\varepsilon_{k-1}\varphi
\left(-{z_{k-1}\over
1+z_{k-1}}\right)+R^{[k]}(\varphi)(z_{0})\; ,\cr
}
$$
which gives
$$
|T^k\varphi (z_{0})|\le cq_{k-1}^{-1}\Vert\varphi\Vert_{H^\infty}
+ck\left({\sqrt{5}+1\over 2}\right)^{-k}
\sup_{D_{\infty}}|\varphi|\; .
$$

Collecting the estimates obtained in (i)--(v)
we get the desired result.  \qed
\smallskip

\smallskip
{\bf 4.3.3} Let us now consider the conformally invariant version of
$H^p$, $1<p<+\infty$.
The application $w\mapsto z={(w+1)^{2}\over 4w}$ is a conformal
representation of $\D$ on $\C\setminus [0,1]$; letting
$w=e^{2\pi i \theta}$, one finds $z=\cos^{2}\theta$ thus
$|dz|=2|\sin\theta||\cos\theta|d\theta$, i.e.
$d\theta = {dx\over 2\sqrt{x(1-x)}}$. We can therefore identify
$\Lp{p}{\T^{1}}$ with the direct sum of two copies of
$\Lp{p}{[0,1],{dx\over 2\sqrt{x(1-x)}}}$ (one for each side of
$[0,1]$). Let now $\varphi\in\O{1}{[0,1]}\cap
\Hp{p}{\overline{\C}\setminus [0,1]}$ which we
will assume to be {\it real}. The associated hyperfunction $u$ is
$u = \IM \varphi (x+i0)$, and $u\in \Lp{p}{[0,1],{dx\over
2\sqrt{x(1-x)}}}$. We have the following

\smallskip
\Proc{Proposition 4.6}{For $p\ge 1$
the operator $T$ defines a bounded operator on
$\Lp{p}{[0,1],{dx\over 2\sqrt{x(1-x)}}}$ and its spectral radius is
$\le {\sqrt{5}-1\over 2}$. }

\smallskip
\proof
One has $Tu(x)=xu\left({1\over x}-m\right)$ if
${1\over m+1}\le x\le {1\over m}$, thus
$$
\eqalign{
\int_{1/(m+1)}^{1/m}|Tu(x)|^p{dx\over 2\sqrt{x(1-x)}} &\le
m^{1/2-p}\int_{1/(m+1)}^{1/m}\left|u\left({1\over x}-m\right)\right|^p
{dx\over 2\sqrt{(1-x)}} \cr
&\le m^{1/2-p}\int_{0}^{1}|u(s)|^p\left(1-{1\over s+m}\right)^{-1/2}
{ds\over (s+m)^{2}}\cr
&\le m^{1/2-p}\int_{0}^{1}|u(s)|^p (s+m)^{-3/2} {ds\over\sqrt{s}}\; .
\cr
}
$$
Summing over $m\ge 1$ we have
$$
\Vert Tu\Vert_{\Lp{p}{[0,1],{dx/ 2\sqrt{x(1-x)}}}}\le
\left[\sum_{m\ge 1}m^{-p-1}\right]^{1/p}\Vert u
\Vert_{\Lp{p}{[0,1],{dx/\sqrt{x}}}}\; ,
$$
which proves the first assertion.

As far as the spectral radius is concerned, le us consider an integer
$k>1$ and the set ${\cal I}(k)$ of the intervals of definition of the
branches of $A^k$. For $I\in {\cal I}(k)$ one has
(we denote with $x_{I}$ the center of $I$)
$$
\eqalign{
\int_{I} |T^ku(x)|^p{dx\over 2\sqrt{x(1-x)}} &=
\int_{I}\beta_{k-1}(x)^p |u(A^k x)|^p {dx\over 2\sqrt{x(1-x)}}\cr
&\le c|I|^p\left[x_{I}(1-x_{I})\right]^{-1/2}\int_{I}
|u(A^k x)|^p dx\cr
}
$$
since $c^{-1}\left[x_{I}(1-x_{I})\right]^{-1/2}\le
\left[x(1-x)\right]^{-1/2}\le c \left[x_{I}(1-x_{I})\right]^{-1/2}$
for all $x\in I$, $I\in{\cal I}(k)$ and $k>1$; also
$c^{-1}|I|\le \beta_{k-1}(x)\le c|I|$ for all $x\in I$. From this we
get
$$
\int_{I} |T^ku(x)|^p{dx\over 2\sqrt{x(1-x)}} \le
c|I|^{p+1} \left[x_{I}(1-x_{I})\right]^{-1/2} \Vert u
\Vert_{\Lp{p}{[0,1]}}^p\; ,
$$
(since $A^k$ has bounded distorsion on $I$).

Taking the sum over all the intervals $I$ (i.e. the branches of
$A^k$), since the $I\in {\cal I}(k)$ form a partition $(\hbox{mod}\, 0)$
of $[0,1]$, one obtains
$$
\Vert T^ku \Vert_{\Lp{p}{[0,1],{dx\over 2\sqrt{x(1-x)}}}}^p\le
c\int_{0}^{1}{dx\over\sqrt{x(1-x)}} \Vert u\Vert_{\Lp{p}{[0,1]}}^p
\left[\max_{{\cal I}(k)}|I|\right]^p\; ,
$$
which proves the second part of the proposition. \qed
\smallskip

\smallskip
\Proc{Corollary 4.7}{For $p> 1$, the operator $T$ maps $\O{1}{[0,1]}\cap
\Hp{p}{\overline{\C}\setminus [0,1]}$ into itself, is bounded and its
spectral radius is $\le {\sqrt{5}-1\over 2}$}

\smallskip
\proof
Let $\varphi\in \O{1}{[0,1]}\cap
\Hp{p}{\overline{\C}\setminus [0,1]}$ be real, $u=\IM \varphi (x+i0)$
the associated hyperfunction. For $k\ge 0$ the function of
$\O{1}{[0,1]}$ associated to $T^k u$ is $T^k\varphi$. Using the
conformal representation of $\overline{\C}\setminus [0,1]$ onto
$\D$ and the fact that the Hilbert transform is bounded on
$L^p (d\theta)$ one obtains the desired result. \qed
\smallskip

\smallskip
{\bf 4.3.4} Here we consider the operator $\sum_{\cM}$ acting on the
space
$\Hp{1}{\overline{\C}\setminus [0,1]}\cap \O{1}{[0,1]}$.

\smallskip
\Proc{Lemma 4.8}{For all $g\in \cM$ the restriction to
$\Hp{1}{\overline{\C}\setminus [0,1]}\cap \O{1}{[0,1]}$ of
$L_{g}$ is a bounded operator of this space into itself.}

\smallskip
\proof
It is sufficient to consider the case $g=g(m)$, $m\ge 1$. One has then
$$
L_{g(m)}=
\sigma\circ\tau^{-1}\circ\chi_{1}\circ\chi_{0}\circ\tau\circ\sigma\circ
\iota_{m}\circ\tau_{m}
$$
where
\item{$\tau_{m}$} is the isomorphism $\varphi (z) \mapsto
\varphi \left({1\over z}-m\right)$ of $\Hp{1}{\overline{\C}\setminus [0,1]}$
onto $\Hp{1}{\overline{\C}\setminus [1/(m+1),1/m]}$;
\item{$\iota_{m}$} is the canonical injection of
$\Hp{1}{\overline{\C}\setminus [1/(m+1),1/m]}$ into
$\Hp{1}{\overline{\C}\setminus [0,1]}$;
\item{$\sigma$} is the bounded operator $\varphi\mapsto \varphi
-\varphi (\infty )$ of $\Hp{1}{\overline{\C}\setminus [0,1]}$
into $\Hp{1}{\overline{\C}\setminus [0,1]}\cap \O{1}{[0,1]}$;
\item{$\tau$} is the isomorphism $\varphi\mapsto \varphi_{1}(w)
=\varphi\left({(w+1)^{2}\over 4w}\right)$ of
$\Hp{1}{\overline{\C}\setminus [0,1]}$ onto $\Hp{1}{\D}$, whose
restriction to $\Hp{1}{\overline{\C}\setminus [0,1]}\cap \O{1}{[0,1]}$
is an isomorphism onto $\{\varphi_{1}\in \Hp{1}{\D}\; , \;
\varphi_{1}(0)=0\}$;
\item{$\chi_{0}$} is the isomorphism $\varphi_{1}(w)\mapsto {1\over
w}\varphi_{1}(w)$ of $\{\varphi_{1}\in \Hp{1}{\D}\; , \;
\varphi_{1}(0)=0\}$ onto $\Hp{1}{\D}$;
\item{$\chi_{1}$} is the multiplication operator by the function
$-{(w+1)^{2}\over 4}\in\Hp{\infty}{\D}$ into
$\Hp{1}{\D}$. \qed
\smallskip
\noindent
We want now to estimate the norm of $L_{g}$ acting on
$\Hp{1}{\overline{\C}\setminus [0,1]}\cap \O{1}{[0,1]}$.

\Proc{Proposition 4.9}{There exists a constant $K>0$ such that,
if $g\in\cM$ is different from
$\left(\matrix{1 & 0\cr 0 & 1\cr}\right)$,
$\left(\matrix{0 & 1\cr 1 & 1\cr}\right)$ (thus $d>b>0$) one has}
$$
\Vert L_{g}\Vert_{H^{1}}\le K d^{-5/2}[\min (b,d-b)]^{-1/2}
\log (1+d)\; .
$$

\smallskip
\proof
A function $\varphi\in \Hp{1}{\overline{\C}\setminus [0,1]}$ has non
tangential limits at almost all points of the boundary and
$$
\Vert\varphi\Vert_{H^{1}} = \int_{0}^{1}
{|\varphi (x+i0)|dx \over
2\sqrt{x(1-x)}} + \int_{0}^{1}
{|\varphi (x-i0)|dx \over
2\sqrt{x(1-x)}} \; .
$$
On the other hand a function $\varphi\in \Hp{1}{\overline{\C}\setminus [0,1]}$
verifies
$$
\eqalign{
|\varphi (t)| &\le C |t|^{-1/2}\Vert\varphi\Vert_{H^{1}}\;\;  \forall
t\in (-1,0)\; , \cr
|\varphi (t)| &\le C (t-1)^{-1/2}\Vert\varphi\Vert_{H^{1}}\;\;  \forall
t\in (1,2)\; , \cr
}\eqno(4.14)
$$
as one can check directly applying Poisson's integral formula.
Moreover,  if it belongs to
$\Hp{1}{\overline{\C}\setminus [0,1]}\cap \O{1}{[0,1]}$ the same
argument leads to the estimates
$$
%\eqalign{
|\varphi (z)| \le C|z|^{-1}\Vert\varphi\Vert_{H^{1}} \; ,\quad
|\varphi' (z)| \le C|z|^{-2}\Vert\varphi\Vert_{H^{1}}\; ,\quad
|\varphi'' (z)| \le C|z|^{-3}\Vert\varphi\Vert_{H^{1}}\; ,
\eqno(4.15)
$$
for all $z$ such that $|z-1/2|>1$.

Given these preliminary elementary estimates, let $\varphi\in
\Hp{1}{\overline{\C}\setminus [0,1]}\cap \O{1}{[0,1]}$,
$g=\abcd\in\cM$, with $d>b>0$. Let
$\tilde{\varphi} = L_{g}\varphi$; $\tilde{\varphi}$ is holomorphic
outside  the interval with end points $b/d$ and $(a+b)/(c+d)$.
We must estimate $\int_{0}^{1}
{|\tilde{\varphi} (x\pm i0)|dx \over
\sqrt{x(1-x)}}$ and by symmetry it is enough to consider
$\tilde{\varphi} (x+ i0)$. We define an interval $I_{g}$
of the following form:
\item{(i)} if $g=\left(\matrix{0 & 1\cr 1 & m\cr}\right)$, $m\ge 2$,
$I_{g}=[0,3/2m]$;
\item{(ii)} if $g= \left(\matrix{1 & m-1\cr 1 & m\cr}\right)$,
$m\ge 2$, $I_{g}=[1-3/2m, 1]$;
\item{(iii)} in all the other cases one has $c>a>0$,
${a\over c}={b\over d}+{\eps_{g}\over cd}$; $I_{g}$
is the interval whose end points ${b\over d}+{3\eps_{g}\over
2cd}$ and ${b\over d}-{\eps_{g}\over 2cd}$; clearly
$I_{g}\subset [0,1]$.

\noindent
We will now directly estimate the integral of
${\tilde{\varphi}(x+ i0)dx \over
\sqrt{x(1-x)}}$ on $I_{g}$ and on $[0,1]\setminus I_{g}$. We begin
with the latter.

If $x\in [0,1]\setminus I_{g}$, $\tilde{\varphi}$ is holomorphic in a
neighborhood of $x$ and one has
$$
\tilde{\varphi}(x) = c^{-2}(a-cx)^{-1}\int_{0}^{1}(1-t)\varphi''
\left(-{d\over c}+\eps_{g}{t\over c(a-cx)}\right)dt\; .
$$
We have
$$
x\notin I_{g}\Rightarrow {\eps_{g}\over c(a-cx)}\in
\left[-2{d\over c},{2\over
3}{d\over c}\right]\eqno(4.16)
$$
thus the values of the second derivative in the integral are $\le C
c^{3}d^{-3}\Vert \varphi\Vert_{H^{1}}$,  by the third estimate a
of (4.15), and  one has
$$
|\tilde{\varphi}(x)|\le C
d^{-3}|x-b/d|^{-1}\Vert\varphi\Vert_{H^{1}}\; ,
$$
since $x-a/c$ and $x-b/d$ are comparable outside $I_{g}$. We thus
obtain
$$
\int_{[0,1]\setminus I_{g}} {|\tilde{\varphi}(x\pm i0)|dx\over
\sqrt{x(1-x)}}\le C d^{-3}\Vert\varphi\Vert_{H^{1}}
\int_{[0,1]\setminus I_{g}}{dx\over \sqrt{x(1-x)}\left|x-{b\over
d}\right|}\; .
$$
An easy estimate of this last integral gives
$$
\int_{[0,1]\setminus I_{g}}{|\tilde{\varphi}(x\pm i0)|dx\over
\sqrt{x(1-x)}}\le Cd^{-5/2}[\min (b,d-b)]^{-1/2}\log
(1+\min (b,d-b))\; .
$$
For the integral inside $I_{g}$ we distinguish the three different
contributions to $\tilde{\varphi}=L(g)\varphi$.
First of all one has (applying the second estimate of (4.15))
$$
|c^{-1}\varphi' (-d/c)|\le C cd^{-2}\Vert\varphi\Vert_{H^{1}}\; ,
$$
and
$$
\int_{I_{g}}{dx\over\sqrt{x(1-x)}}\le C c^{-1}d^{-1/2}[\min
(b,d-b)]^{-1/2}
\Vert\varphi\Vert_{H^{1}}\; ,
$$
thus
$$
\int_{I_{g}}{|c^{-1}\varphi' (-d/c)|dx\over \sqrt{x(1-x)}}\le
Cd^{-5/2}[\min (b,d-b)]^{-1/2}\Vert\varphi\Vert_{H^{1}}\; .
$$
For $x\in I_{g}$ one has applying (4.15) and (4.16)
$$
|(a-cx)\varphi (-d/c)|\le Ccd^{-2}\Vert\varphi\Vert_{H^{1}}\; ,
$$
from which follows the same estimate above for the second term.
We are left with
$$
I= \int_{I_{g}}{|a-cx||\varphi (y\pm i0)|\over \sqrt{x(1-x)}}dx\; , \;\;
y={dx-b\over a-cx}\; .
$$
Inside $I_{g}$ one has
$
{\dst |a-cx|\over\dst \sqrt{x(1-x)}}\le Cd^{-1/2}[\min (b,d-b)]^{-1/2}
\; ;
$
on the other hand one has
$
\int_{I_{g}}|\varphi (y\pm i0)|dx = c^{-2}\int
{\dst |\varphi (y\pm i0)|dy\over\dst  (y+d/c)^{2}}
\; ,
$
where this integral is taken on the complement in $\overline{\R}$ of
the interval $\left[-{\dst 3d\over\dst c},-{\dst d\over\dst 3c}\right]$
(if $c>1$; when $g=\left(\matrix{0 & 1\cr 1 & m\cr}\right)$ or
$\left(\matrix{1 & m-1\cr 1 & m\cr}\right)$, $m\ge 2$, the integral
is taken from $-m/3$ to $+\infty$).

One has then
$$
\eqalign{
\int_{-1/2}^{3/2}{|\varphi (y\pm i0)|dy\over (y+d/c)^{2}}
& \le Cc^{2}d^{-2}\Vert\varphi\Vert_{H^{1}}\; , \cr
\int_{-d/3c}^{-1/2}{|\varphi (y\pm i0)|dy\over (y+d/c)^{2}}
& \le Cc^{2}d^{-2}\Vert\varphi\Vert_{H^{1}}\log \left(1+{d\over c}\right)
\; , \cr
}
$$
and the same for $\int_{3/2}^{3d/c}$.  Finally one has
$$
\int_{|y|>3d/c}{|\varphi (y\pm i0)|dy\over (y+d/c)^{2}}
\le Cc^{2}d^{-2}\Vert\varphi\Vert_{H^{1}}\; .
$$
Putting the various bounds together one gets
$$
\int_{I_{g}} {|\tilde{\varphi}(x\pm i0)|dx\over \sqrt{x(1-x)}}\le
C d^{-5/2}[\min(b,d-b)]^{-1/2}\log (1+d/c)\Vert\varphi\Vert_{H^{1}}\; .
$$
To conclude the proof it is now enough to note that
$\log (1+d/c)\le \log (1+d)$ and
$\log (1+c\min (b,d-b))\le \log (1+d^{2})\le 2\log (1+d)$. \qed
\smallskip

\smallskip
\Proc{Corollary 4.10}{The series $\sum_{\cM}\Vert L_{g}\Vert_{H^{1}}$
is convergent. $\sum_{\cM}$ defines a bounded operator of
$\Hp{1}{\overline{\C}\setminus [0,1]}\cap \O{1}{[0,1]}$ into itself;
moreover $T$ is a contraction of this space, its spectral radius being
bounded above by ${\sqrt{5}-1\over 2}$. }

\smallskip
\proof
Let $k\ge 1$. We will take the sum of $\Vert L_{g}\Vert_{H^{1}}$
on the elements $g\in \cM$ which are the product of exactly $k$
generators $g(m)$. The branch of $A^k$ associated to $g$ has domain
$\left[{a+b\over c+d},{b\over d}\right]=I(g)$ and on this interval
the $\max$ and the $\min $ of $x(1-x)$ are comparable. The length
$[d(c+d)]^{-1}$ of $I(g)$ is bounded below by ${1\over 2}d^{-2}$, thus
$$
\Vert L_{g}\Vert_{H^{1}}\le {Kd^{-3}\log (1+d)\over
\sqrt{{b\over d}\left(1-{b\over d}\right)}}\le
K' d^{-1}\log (1+d)\int_{I(g)}{dx\over\sqrt{x(1-x)}}\; .
$$
When we sum on the elements $g$ considered, the intervals
$I(g)$ form a partition of $[0,1]$ $\hbox{mod}\, 0$
thus we get
$$
\sum_{m_{1},\ldots ,m_{k}\, , g=g(m_{1})\cdots
g(m_{k})}\Vert L_{g}\Vert_{H^{1}}\le K'' \max_{m_{1},\ldots ,m_{k}}
[d^{-1}\log (1+d)]
$$
which gives the desired result. \qed
%\vskip .5 truecm
\beginsection 4.4 Real holomorphic functions with bounded real part
\par
% \vskip .3 truecm
We denote by $E$ the space of functions $\varphi\in\O{1}{[0,1]}$
whose real part is bounded, endowed with the norm
$
\Vert \varphi\Vert_{E} = \sup\nolimits_{\overline{\C}\setminus [0,1]}
|\RE \varphi |\; .
$
We have then  for $|\IM z|\le 1/2$
$$
|\IM\varphi (z) |\le {2\over \pi}\log [(2\sqrt{2}-2)|\IM z|]^{-1}
\Vert \varphi\Vert_{E}\; ,\eqno(*)
$$
as one can prove from the analogue estimate for functions in the unit
disk $\D$
$$
|\IM\Phi (w)|\le {2\over \pi}\sup_{\vartheta\in [0,1]}|
\RE\Phi (e^{2\pi i \vartheta})|\log{1+|w|\over 1-|w|}
$$
applying the conformal representation of $\overline{\C}\setminus
[0,1]\rightarrow \D$, $w=(\sqrt{z}-\sqrt{z-1})^{2}$.
\smallskip
\Proc{Proposition 4.11}{The restriction of $T$ to $E$ is a bounded
operator with spectral radius $\le {\sqrt{5}-1\over 2}$. }

\smallskip
\proof
Let $\varphi\in E$, $k\ge 1$, $z_{0}\in
\overline{\C}\setminus [0,1]$. We estimate
$\RE T^k\varphi (z_{0})$ in various cases.
\item{(i)} If $z_{0}\in D_{\infty}$ we have
$$
|\RE T^k\varphi (z_{0})|\le |T^k\varphi (z_{0})|\le
c\left({\sqrt{5}+1\over 2}\right)^{-k}\sup_{D_{\infty}}|\varphi |\; ,
$$
and, on the other hand, for all $\varphi\in E$,
$
\sup_{D_{\infty}}|\varphi |\le C\Vert \varphi\Vert_{E}\; .
$
\item{(ii)} If $z_{0}\in D_{0}\cup H_{0}\cup H_{0}'$, or if
$z_{0}\in D_{1}$, $k>1$, or if $z_{0}\in D$ has continued fraction
$z_{i}^{-1}=m_{i+1}+z_{i+1}$, $0\le i<l$, with $z_{l}\in
\Delta\setminus D$, $l<k-1$, we have obtained in the proof of
Propostion 4.5 the estimate
$$
|T^k\varphi (z_{0})|\le
c\,k\left({\sqrt{5}+1\over 2}\right)^{-k}\sup_{D_{\infty}}|\varphi |\; .
$$
\item{(iii)} If $z_{0}\in D_{1}$, $k=1$, we have from Proposition 4.2
$$
\left|T\varphi (z_{0})+z_{0}\varphi\left({1\over
z_{0}}-1\right)\right|\le
c\sup_{D_{\infty}}|\varphi |\; ,
$$
we have
$$
\left|\RE z_{0}\RE \varphi\left({1\over
z_{0}}-1\right)\right|\le C\Vert \varphi\Vert_{E}\; ,
$$
and
$$
\left|\IM \left({1\over
z_{0}}-1\right)\right|\ge  C^{-1}|\IM z_{0}|
$$
hence from $(*)$
$$
\left|\IM z_{0}\IM \varphi\left({1\over
z_{0}}-1\right)\right|\le C\Vert \varphi\Vert_{E}\; .
$$
\item{(iv)} If $z_{0}\in D$ has continued fraction
$z_{i}^{-1}=m_{i+1}+z_{i+1}$, $0\le i<k-1$, with $z_{k-1}\in
\Delta\setminus D$, the only term  in the proof of
Propostion 4.5 which gives some trouble is
$
(p_{k-2}-q_{k-2}z_{0})\varepsilon_{k-2}T\varphi (z_{k-1}+1)
$
(the others are once again dominated by $Ck
\left({\sqrt{5}+1\over 2}\right)^{-k}\sup_{D_{\infty}}|\varphi |$).
In fact, from Proposition 4.2, we are even left with
$$
(p_{k-2}-q_{k-2}z_{0})\varepsilon_{k-2}(1+z_{k-1})\varphi\left(
{1\over 1+z_{k-1}}-1\right)\; .
$$
We have here
$
(p_{k-2}-q_{k-2}z_{0})(1+z_{k-1})=
(p_{k-2}-q_{k-2}z_{0})-(p_{k-1}-q_{k-1}z_{0})
$
hence
$$
\eqalign{
\left|\RE [(p_{k-2}-q_{k-2}z_{0})(1+z_{k-1})]\right|
&\le Cq_{k-1}^{-1}\cr
\left|\IM [(p_{k-2}-q_{k-2}z_{0})(1+z_{k-1})]\right|
&\le Cq_{k-1}|\IM z_{0}|
\le Cq_{k-1}^{-1}|\IM z_{k-1}|\cr
}
$$
and
$$
\left|\IM\left({1\over 1+z_{k-1}}-1\right)\right|\ge C^{-1}|\IM
z_{k-1}|\; .
$$
Thus, from $(*)$ we get
$$
\left|\RE (p_{k-2}-q_{k-2}z_{0})(1+z_{k-1})\varphi\left(
{1\over 1+z_{k-1}}-1\right)\right|\le
Cq_{k-1}^{-1}\Vert\varphi\Vert_{E}\; ,
$$
and finally
$
|\RE T^k\varphi (z_{0})|\le Ck\left({\dst\sqrt{5}+1\over\dst
2}\right)^{-k} \Vert\varphi\Vert_{E}\; .
$
\item{(v)} If $z_{0}\in D$ has continued fraction
$z_{i}^{-1}=m_{i+1}+z_{i+1}$, $0\le i<k$, with $z_{k}\in
\Delta$, we apply Proposition 4.4. We have
$$
\eqalign{
|\RE (p_{k-1}-q_{k-1}z_{0})| &\le Cq_{k}^{-1}\cr
|\IM (p_{k-1}-q_{k-1}z_{0})| &\le q_{k-1}|\IM z_{0}|
\le Cq_{k-1}q_{k}^{-2}|\IM z_{k}|\cr
}
$$
hence, for $\varepsilon =-1,0,1$, from $(*)$
$$
|\RE [(p_{k-1}-q_{k-1}z_{0})\varphi (z_{k}+\varepsilon )]|
\le Cq_{k}^{-1}\Vert\varphi\Vert_{E}\; .
$$
We deal similarly with
$$
\RE \left[(p_{k-2}-q_{k-2}z_{0})(1+z_{k-1})\varphi\left(-
{z_{k-1}\over 1+z_{k-1}}\right)\right]
$$
(see (iv) above) and conclude that
$$
|\RE T^k\varphi (z_{0})|\le Ck\left({\sqrt{5}+1\over 2}\right)^{-k}
\Vert\varphi\Vert_{E}\; .
$$
\qed
%\vskip .5 truecm

%%%%%%%%%%%%%%%% 5. The complex Brjuno function %%%%%%%%%%%%%%
\beginsection 5. The complex Brjuno function  \par
\vskip .5 truecm
In this Section we introduce and study the complex Brjuno function.
Preliminarly we need some further results on the monoid ${\cal M}$ and
on the algebraic properties of $\sum_{\cal M}$ and $\sum_{Z{\cal M}}$.

\smallskip
{\bf 5.0.1}
Let us recall that one has
$$
Z\cM = Z\cM\left(\matrix{1 & 1\cr 0 & 1\cr}\right)\sqcup
Z\cM\left(\matrix{0 & 1\cr 1 & 1 \cr}\right)\; .
$$
More precisely, if one denotes $\cM^{*}=\cM\setminus\{1\}$, one has
$$
\cM\left(\matrix{1 & 1\cr 0 & 1\cr}\right)\sqcup
\cM\left(\matrix{0 & 1\cr 1 & 1 \cr}\right)=\cM^{*}
\sqcup\left\{\left(\matrix{1 & 1\cr 0 & 1\cr}\right)\right\} \; .
$$

\smallskip
{\bf 5.0.2}
Let $\varphi\in\O{1}{[0,1]}$. Let
$$
\varphi_{1} = \left(L_{\left(\matrix{0 & 1\cr 1 & 1 \cr}\right)}+
L_{\left(\matrix{1 & 1\cr 0 & 1\cr}\right)}\right)\varphi\in
\O{1}{[1/2,2]}\; .\eqno(5.1)
$$
The family $(L_{g}\varphi)_{g\in \cM}$ is uniformly summable
on compact subsets of $\overline{\C}\setminus [0,1]$, thus
$$
\sum_{\cM}\varphi (z)-\varphi (z) = \sum_{\cM}\varphi_{1}(z)-
\varphi (z-1)\; .
$$
On the other hand
$
\sum\limits_{\cM}\varphi_{1} = \varphi_{1}+
\sum\limits_{\cM}(T\varphi_{1} )\; ,
$
and $T\varphi_{1}\in \O{1}{[0,1]}$. We can therefore conclude that
$$
\sum_{Z\cM}\varphi (z) = \sum_{Z}[\varphi_{1}+\sum_{\cM}
(T\varphi_{1} )]\; , \eqno(5.2)
$$
where the definition of $\sum_{Z}\varphi_{1} $,
$\varphi_{1}\in \O{1}{[1/2,2]}$, is obtained extending the one
given for functions in  $\O{1}{[0,1]}$.

\smallskip
\remark{5.1} One has
$$
\varphi_{1}(z)=-z\left[\varphi\left({1\over z}-1\right)-\varphi
(-1)\right]+\varphi'(-1)+\varphi (z-1)\; , \eqno(5.3)
$$
from which it follows that
$$
z\varphi_{1}\left({1\over z}\right) = - [\varphi (z-1)-\varphi (-1)]
+z\varphi' (-1)+z\varphi\left({1\over z}-1\right)\; ,
$$
and
$$
\varphi_{1}(z)+z\varphi_{1}\left({\dst 1\over\dst z}\right)
=(1+z)(\varphi (-1) +\varphi' (-1))\; .
$$

%\vskip .5 truecm
\beginsection 5.1 The dilogarithm \par
%\vskip .3 truecm
\indent
{\bf 5.1.1} Let us define
$$
\varphi_0 (z) =-{1\over \pi}\Li \left({1\over z}\right)\; , \eqno(5.4)
$$
where the dilogarithm is taken with its principal branch in
${\Bbb C}\setminus [1,+\infty ]$ (see Appendix 3 for a short summary
of the properties of the dilogarithm and [O] and references therein
for more details).
The function $\varphi_{0}$
belongs to $\O{1}{[0,1]}$. It is real on the real axis
outside $[0,1]$ and its only singular points are $0$ and $1$.
It is bounded outside of any neighborhood of $0$ and
$$
\IM \varphi_{0} (x\pm i0) = \pm \log {1\over x}\; , \;\; 0<x\le 1\; ,
\eqno(5.5)
$$
thus the relation with the real Brjuno function $B$ is clear:
$B(x) = [(1-T)^{-1}f](x)$ with $f(x)=\sum_{n\in \Z}
\IM \varphi_{0} (x+ i0-n)$.

\smallskip
{\bf 5.1.2} Let us now consider
$$
\varphi_{1}= \left[ L_{\left(\matrix{ 1 & 1 \cr 0 & 1}\right)}+
L_{\left(\matrix{ 0 & 1 \cr 1 & 1}\right)}\right]\varphi_{0}
\in \O{1}{[1/2,2]}\; . \eqno(5.6)
$$
Since $\varphi_{0}(-1) =\pi/12$ and $\varphi_{0}'(-1)={1\over\pi}\log
2$ one has
$$
\varphi_{1}(z) = {1\over\pi}\left[ z\Li \left({z\over 1-z}\right)
-\Li\left({1\over z-1}\right)\right]+{\pi\over 12}z+{1\over\pi}\log 2.
\eqno(5.7)
$$
The function $\varphi_{1}$ is real. It admits as unique singularities
the points $1/2$, $1$ and $2$ and has two cuts along $(1/2,1)$ and
$(1,2)$. It can be continuously extended to $1/2$ and $2$ and it is
bounded outside any neighborhood of $1$. Moreover
$$
\IM \varphi_{1}(x\pm i0) = \cases{ \pm x\log {x\over 1-x} & if $1/2\le
x<1\; ,$\cr
\pm \log {1\over x-1} & if $1<x\le 2\; .$\cr} \eqno(5.8)
$$
Note that if $1/2\le x <1$
$$
x\log {x\over 1-x} = \log {1\over 1-x} + x\log x + (1-x) \log
(1-x)\; .
$$
One also has
$$
\varphi_{1}(z)+z\varphi_{1}\left({1\over z}\right) =
(1+z)\left({\pi\over 12}+{1\over\pi}\log 2\right)\; . \eqno(5.9)
$$

\smallskip
\Proc{Lemma 5.2}{The function $\varphi_{1}(z) +i\log (1-z)$ is
continuous on $\overline{\H^{+}}$ and its value at $1$ is ${1\over
\pi}\log 2+{7\pi\over 12}$. }

\smallskip
\proof
Applying (A3.7) to (5.7) twice one gets
$$
\eqalign{
\varphi_{1}(z) &= {1\over \pi}\log 2-{\pi\over 6}+{\pi\over 4}z
+{1\over\pi}\left[\Li (z-1)-z\Li \left({1-z\over z}\right)\right]\cr
& + {1\over 2\pi} \left[\log^{2}(1-z)-z\log^{2}{1-z\over -z}\right]\cr
}
$$
In this expression the function $\Li (z-1) - z\Li \left({1-z\over
z}\right)$  is regular and vanishing at $z=1$. Moreover
$$
\log^{2}(1-z)-z\log^{2}{1-z\over -z} =
\log^{2}(1-z)-\log^{2}{1-z\over -z} + (1-z)\log^{2}{1-z\over -z}\; ,
$$
where $(1-z)\log^{2}{1-z\over -z}$ vanishes at $z=1$, and
$$
\log^{2}(1-z)-\log^{2}{1-z\over -z} = -\log^{2} (-z) + 2\log (-z)
\log (1-z )\; .
$$
In a neighborhood of $1$ in $\overline{\H^{+}}$ one has
$\log (-z)+i\pi= \hbox{O}\, (|z-1|)$ thus
$$
\log^{2}(1-z)-\log^{2}{1-z\over -z} = \pi^{2} - 2i\pi \log (1-z)
+ \hbox{O}\, \left(|z-1|\log {1\over |z-1|}\right) \; .
$$
\qed

\smallskip
This Lemma leads to the following important

\smallskip
\Proc{Corollary 5.3}{The real part of $\varphi_{1}$ is bounded in
$\overline{\C}\setminus [1/2 ,2]$. It has an extension to a
continuous function on $\overline{\C}\setminus\{1\}$ and}
$$
\lim_{x\rightarrow 1^\pm}\RE \varphi_{1}(x) =
{1\over \pi}\log 2 +{\pi \over 12} \mp {\pi\over 2}\; . \eqno(5.10)
$$

\smallskip
This Corollary  is the motivation for using $\varphi_{1}$ instead of
$\varphi_{0}$ as the starting point of the construction of the
complex Brjuno function. Equation (5.2) shows that this leads to the
same result.

\vskip .5 truecm
\beginsection 5.2 A natural compactification of $\H^{+}$ \par
\vskip .3 truecm
By Lemma 5.2 above, $\RE \varphi_{1}$ extends continuously to
$\overline{\H^{+}}\setminus\{1\}$ with limits at $1$
along rays. This means that $\RE \varphi_{1}$ extends
continuously to the compactification of $\H^{+}$ obtained from
$\overline{\H^{+}}$ by blowing out $1$ into a semi--circle
(corresponding to all rays in $\H^{+}$ which end in $1$). If we
want to obtain a similar result for the complex Brjuno function
$\sum_{Z{\cal M}}\varphi_{0}$
we have to do the same thing at every point of $\overline{\Q}$.

\smallskip
{\bf 5.2.1} We will consider
$$
\hat{\H}^{+} = \H^{+}\sqcup (\R\setminus\Q )\sqcup
\left(\overline{\Q}\times \left[-{\pi\over 2},+{\pi\over
2}\right]\right) \eqno(5.6)
$$
(where $\overline{\Q}=\Q\cup\{\infty\}$)
equipped with the topology defined by the following
fundamental system of neighborhoods at any point $z\in\hH+$:
\item{a)} if $z_{0}\in \H^{+}$ a fundamental system of neighborhoods
is given by $\{|z-z_{0}|<\varepsilon\}$, $0<\varepsilon <\IM z_{0}$;
\item{b)} if $\alpha_{0}\in \R\setminus\Q$, a fundamental system of
neighborhoods
is given by the sets ($\varepsilon >0$)
$$
\eqalign{
V_{\varepsilon}(\alpha_{0}) &= \{ z\in \H^{+}\; , \; |z-\alpha_{0}|
<\varepsilon \} \cr
& \cup \{\alpha\in\R\setminus\Q\; ,\; |\alpha
-\alpha_{0}|<\varepsilon\} \cr
&\cup \left\{(\alpha ,\theta )\in \Q\times \left[-{\pi\over 2},+{\pi\over
2}\right]\; , \; |\alpha -\alpha_{0}|<\varepsilon\right\}\; ;  \cr}
$$
\item{c)} if $\alpha_{0}\in \Q$ a fundamental system of neighborhoods
of $(\alpha_{0}, \pi/2)$ (resp. $(\alpha_{0}, -\pi/2)$) is given by
($0<\varepsilon <\pi$)
$$
\eqalign{
V_{\varepsilon} (\alpha_{0}, \pi/2) &=
\{z\in\H^{+}\; , \; |z-\alpha_{0}|<\varepsilon\; , \;
0<\arg (z-\alpha_{0})<\varepsilon\}\cr
&\cup \{\alpha\in\R\setminus\Q\; , \; 0<\alpha
-\alpha_{0}<\varepsilon\} \cr
&\cup \left\{(\alpha ,\theta )\in \Q\times \left[-{\pi\over 2},+{\pi\over
2}\right]\; , \; 0<\alpha -\alpha_{0}<\varepsilon\right\}\cr
&\cup \left\{(\alpha_{0},\theta )\; , \; \pi/2-\varepsilon <\theta \le
\pi/2\right\}\; , \cr
V_{\varepsilon}(\alpha_{0}, -\pi/2) &=
\{z\in\H^{+}\; , \; |z-\alpha_{0}|<\varepsilon\; , \;
\pi>\arg (z-\alpha_{0})>\pi -\varepsilon\}\cr
&\cup \{\alpha\in\R\setminus\Q\; , \; 0<\alpha_{0}
-\alpha <\varepsilon\} \cr
&\cup \left\{(\alpha ,\theta )\in
\Q\times \left[-{\pi\over 2},+{\pi\over
2}\right]\; , \; 0<\alpha_{0} -\alpha <\varepsilon\right\}\cr
&\cup \left\{(\alpha_{0} ,\theta )\; , \; -\pi/2\le \theta <
-\pi/2+\varepsilon\right\}\; ; \cr
}
$$
\item{d)} If $\alpha_{0}\in\Q$, $-\pi/2<\theta_{0}<\pi/2$
a fundamental system of neighborhoods of $(\alpha_{0},\theta_{0})$
is given by ($0<\varepsilon <\pi/2-|\theta_{0}|$)
$$
\eqalign{
V_{\varepsilon} (\alpha_{0}, \theta_{0}) &=
\{z\in\H^{+}\; , \; |z-\alpha_{0}|<\varepsilon\; , \;
\theta_{0}-\varepsilon <\pi/2-\arg (z-\alpha_{0})<\theta_{0}+\varepsilon\}\cr
&\cup \left\{(\alpha_{0} ,\theta )\in \Q\times \left[-{\pi\over 2},+{\pi\over
2}\right]\; , \; \theta_{0}-\varepsilon <\theta <\theta_{0}+
\varepsilon \right\}\; ; \cr
}
$$
\item{e)} A fundamental system of neighborhoods of $(\infty ,\pi/2)$
(resp. $(\infty, -\pi/2)$ is given by ($0<\varepsilon <\pi $)
$$
\eqalign{
V_{\varepsilon}(\infty ,\pi/2) &=
\{z\in\H^{+}\; , \; |z|>\varepsilon^{-1}\; , \;
\pi>\arg z>\pi -\varepsilon\}\cr
&\cup \{\alpha\in\R\setminus\Q\; , \; \alpha
<-\varepsilon^{-1}\} \cr
&\cup \left\{(\alpha ,\theta )\in
\Q\times \left[-{\pi\over 2},+{\pi\over
2}\right]\; , \; \alpha
<-\varepsilon^{-1}\right\}\cr
&\cup \left\{(\infty ,\theta )\; , \; \pi/2\ge \theta >
\pi/2-\varepsilon\right\}\; , \cr
V_{\varepsilon}(\infty ,-\pi/2) &=
\{z\in\H^{+}\; , \; |z|>\varepsilon^{-1}\; , \;
0<\arg z<\varepsilon\}\cr
&\cup \{\alpha\in\R\setminus\Q\; , \; \alpha
>\varepsilon^{-1}\} \cr
&\cup \left\{(\alpha ,\theta )\in
\Q\times \left[-{\pi\over 2},+{\pi\over
2}\right]\; , \; \alpha
>\varepsilon^{-1}\right\}\cr
&\cup \left\{(\infty ,\theta )\; , \; -\pi/2\le \theta <
-\pi/2+\varepsilon\right\}\; ; \cr
}
$$
\item{f)} If $-\pi/2<\theta_{0}<\pi/2$ a fundamental system of
neighborhoods of $(\infty ,\theta_{0})$ is given by
($0<\varepsilon <\pi/2-|\theta_{0}|$)
$$
\eqalign{
V_{\varepsilon}(\infty ,\theta_{0}) &=
\{z\in\H^{+}\; , \; |z|>\varepsilon^{-1}\; , \;
\theta_{0}-\varepsilon<\arg z-\pi/2<\theta_{0}+\varepsilon\}\cr
&\cup \left\{(\infty ,\theta )\; , \; |\theta -\theta_{0}|<
\varepsilon\right\}\; .\cr
}
$$

\smallskip
{\bf 5.2.2} One can check that the axioms for a system of neighborhoods are
verified. The following is the only non trivial property: if
$\zeta\in \hH+$ and $V$ is a neighborhood of $\zeta$ then there exists a
neighborhood $W$ of $\zeta$ such that $V$ is a neighborhood of each point
of $W$. This must be checked directly for each of the above listed
cases.

\smallskip
{\bf 5.2.3} It is clear that the topology of $\hH+$ induces on
$\H^{+}$ the usual topology and that $\H^{+}$ is an open dense
subset of $\hH+$.

\smallskip
{\bf 5.2.4} The space $\hH+$ is compact and Hausdorff. More precisely
there exists a homeomorphism of $\overline{\D}$ onto $\hH+$
whose restriction to $\D$ is a homeomorphism onto $\H^{+}$. One can
therefore give to $\hH+$ the structure of a topological manifold with
boundary, the boundary is $\partial\hH+ = \hH+\setminus \H^{+}$
and is homeomorphic to $\S^{1}$.

\smallskip
{\bf 5.2.5} The action of $\hbox{PSL}(2,\Z)$ on $\H^{+}$
by homographies has a {\it continuous} extension to an action
on $\hH+$: just  define
$$
\eqalign{
\left(\matrix{ a & b\cr c & d \cr}\right)\cdot\alpha &=
{a\alpha +b\over c\alpha +d}\; , \;\;
\forall\alpha\in\R\setminus\Q\; , \cr
\abcd\cdot (\alpha ,\theta ) &= \left({a\alpha +b\over c\alpha +d},
\theta \right)\; , \;\; \forall\alpha\in\overline{\Q}\; , \;
\forall\theta\in \left[-{\pi\over 2},+{\pi\over
2}\right]\; . \cr
}\eqno(5.11)
$$
Let $I$ denote a compact non trivial interval in $\R$
and $g=\abcd\in\hbox{PGL}\,(2,\Z)$.
\smallskip
\Proc{Lemma 5.4}{
Assume that $-d/c\notin I$ (so
that $g\cdot I = I'$ is also a compact interval of $\R$).
If $\varphi\in\O{1}{I}$  has the following properties:
\item{(i)} $\varphi$ is real;
\item{(ii)} the harmonic function $\RE\varphi$ on $\H^{+}$
has a continuous extension to $\hH+$.}

{\sl Then the function $L_{g}\varphi\in \O{1}{I'}$
also has these two properties. }

\smallskip
\proof
It is enough to distinguish three cases:
\item{(a)} if $\alpha$ is real and irrational one has
$$
\RE L_{g}\varphi (\alpha ) = (a-c\alpha )\left[\RE\varphi \left(
{d\alpha-b\over a-c\alpha}\right)-\varphi\left(-{d\over
c}\right)\right] - \eps_{g}c^{-1}\varphi'\left(-{d\over c}\right)\; .
$$
Note that from the assumption $-d/c\notin I$ follows that $\varphi
(-d/c)$ and $\varphi' (-d/c)$ are both real.
\item{(b)} If $\alpha = \infty$ and $\theta\in
\left[-{\pi\over 2},+{\pi\over
2}\right]$ is arbitrary then
$\RE L_{g}\varphi (\infty,\theta ) = 0$.
\item{(c)} Finally, if $\alpha$ is a rational number and
$\theta\in
\left[-{\pi\over 2},+{\pi\over
2}\right]$ is arbitrary one has
$$
\RE  L_{g}\varphi (\alpha,\theta ) = (a-c\alpha )\left[\RE\varphi \left(
{d\alpha-b\over a-c\alpha},\eps_{g}\theta \right)-\varphi\left(-{d\over
c}\right)\right]
- \eps_{g}c^{-1}\varphi'\left(-{d\over c}\right)\; .
$$
\qed

\smallskip
{\bf 5.2.6} If $I$ is a non trivial compact interval of $\R$,
we denote $\hC{I}$ the space of holomorphic functions
$\varphi\in\O{1}{I}$ which are {\it real} and whose real part on
$\H^{+}$ extends to a continuous function on $\hH+$.
This is a Banach space with the norm
$$
\eqalign{
\Vert\varphi\Vert_{\hat{\cal C}} &:= \sup\{|\RE\varphi (z)|\; , \;
z\in\H^{+}\}\cr
&= \sup\{|\RE\varphi (z)|\; , \;
z\in\hH+\}\cr
&= \sup\{|\RE\varphi (z)|\; , \;
z\in\partial\hH+\}\cr
}\eqno(5.12)
$$
(the equality of all these norms is a
trivial consequence of the maximum principle).
$\hC{I}$ is a real closed vector subspace of the Banach space
$E(I)$ of holomorphic functions in $\O{1}{I}$ with bounded real part
with respect to the norm
$$
\Vert\varphi\Vert_{E(I)} = \sup_{z\in \overline{\C}\setminus I}
|\RE\varphi (z)|\; \eqno(5.13)
$$
(generalizing the definition used in Section 4.4 for the
particular case $I=[0,1]$).

If $g=\abcd\in\hbox{PGL}\, (2,\Z)$ verifies $-d/c\notin I$, then, by
Lemma 5.4,
$L_{g}$ defines an operator of $\hC{I}$ into $\hC{g\cdot I}$. From
its definition  it is easy to see that this operator is bounded.

Let us assume now that $I\subset (-1,+\infty )$,
$I=[\gamma_{0},\gamma_{1}]$. For all $m\ge 1$, $L_{g(m)}$ defines a
bounded operator of $\hC{I}$ into
$\hC{[1/(m+\gamma_{1}),1/(m+\gamma_{0}]}$. If $\varphi\in\hC{I}$ then
$$
T\varphi = \sum_{m\ge 1}L_{g(m)}\varphi\;
$$
is holomorphic, real and belongs to $\O{1}{[0,1/(1+\gamma_{0})]}$.

\smallskip
\Proc{Proposition 5.5}{The function $T\varphi$ belongs to
$\hC{[0,1/(1+\gamma_{0})]}$ and the operator $T$ from $\hC{I}$ into
this space is bounded. More precisely, for $k\ge 0$, $2^k\le j<j'\le
2^{k+1}$ one has
$$
\Vert\sum_{j\le m<j'}L_{g(m)}\varphi\Vert_{C}\le C_{I}(1+k)2^{-k}
\Vert\varphi\Vert_{\hat{\cal C}}\; ,
$$
thus the series $\sum\RE L_{g(m)}\varphi$ is uniformly convergent in
$\hH+$.}

\smallskip
\proof
As in the proof of Proposition 4.9, by Poisson's formula one has
the following estimates
$$
\eqalign{
|\varphi (z)| &\le C_{I}|z|^{-1}\Vert\varphi\Vert_{\hat{\cal C}}\; , \cr
|D\varphi (z)| &\le C_{I}|z|^{-2}\Vert\varphi\Vert_{\hat{\cal C}}\; , \cr
|D^{2}\varphi (z)| &\le C_{I}|z|^{-3}\Vert\varphi\Vert_{\hat{\cal
C}}\; , \cr
}\eqno(5.14)
$$
provided that  $\hbox{dist}\,(z,I)\ge 1$ and  $\varphi\in\hC{I}$.

Now let $k\ge 0$,  and $j,j'$ such that $2^k\le j<j'\le 2^{k+1}$. Let
also denote $\varphi_{j,j'}=\sum_{j\le m<j'}L_{g(m)}\varphi$.
Clearly  we have $\varphi_{j,j'}\in\hC{I_{k}}$, where
$I_{k}=\left[{1\over 2^{k+1}-1+\gamma_{1}},{1\over
2^k+\gamma_{0}}\right]$. By the maximum principle the supremum of
$|\RE\varphi_{j,j'}|$ on $\hH+$ is attained at a point of the boundary
$\partial\hH+$ of the form $\alpha\in\R\setminus\Q$,
or  $(\alpha ,\vartheta )$, $\alpha\in\Q$, $-\pi/2\le\theta\le\pi/2$,
such that $\alpha\in I_{k}$.
Note that
$
\sum\limits_{j\le m<2^{j'}}|\varphi' (-m)|
\le C2^{-k}\Vert\varphi\Vert_{\hat{\cal C}}
$
and
$
\sum\limits_{j\le m<2^{j'}}|\varphi (-m)|
\le  C\Vert\varphi\Vert_{\hat{\cal C}}\; .
$
If $\alpha$ is irrational and contained in $I_{k}$ one has
$$
\sum_{j\le m<2^{j'}}
\left|\RE\varphi\left({1\over\alpha}-m\right)\right|\le
C(1+k)\Vert\varphi\Vert_{\hat{\cal C}}
$$
and the same estimate holds if $\alpha\in\Q\cap I_{k}$ for all
$\theta\in [-\pi/2,\pi/2]$. Since $|\alpha |<C2^{-k}$ for
$\alpha\in I_{k}$ one obtains the desired inequality which also
implies all the other properties. \qed

\smallskip
Proposition 4.11 leads to the following

\smallskip
\Proc{Corollary 5.6}{The spectral radius of $T$ on $\hC{[0,1]}$
is less or equal to ${\sqrt{5}-1\over 2}$. $\sum_{\cM}$
defines a bounded operator on this space. }

\smallskip
\proof
$\hC{[0,1]}$ is a closed subspace of $E([0,1])$.Ú\qed

\smallskip
{\bf 5.2.7} One constructs a compactification $\hHZ+$ of
$\H^{+}/\Z$ adding the point $i\infty$ and glueing
$(\R\setminus\Q)/\Z\sqcup (\Q/\Z\times [\pi/2,\pi/2])$
in the same way we proceeded for $\hH+$.  One obtains a
topological manifold with boundary homeomorphic to
$\overline{\D}$. The restriction of this homeomorphism to
$\H^{+}/\Z$ is onto $\D^{*}$. The boundary
$\partial\hHZ+ = \hHZ+\setminus (\H^{+}/\Z\cup\{i\infty\})$
is once again homeomorphic to $\S^{1}$.

\smallskip
{\bf 5.2.8} If $I$ is a compact non--trivial interval of
$\R$ and $\varphi\in\O{1}{I}$, we defined in (2.5)
the $1$--periodic holomorphic function
$\sum_{Z}\varphi$ on $\H^{+}$ which extends continuously to
$i\infty$.

\smallskip
\Proc{Proposition 5.7}{Assume that $\varphi\in\hC{I}$. Then
$\sum_{Z}\varphi$ has the following properties:
\item{(i)} $\RE (\sum_{Z}\varphi)$ is bounded on
$\H^{+}$ and the function which it defines on
$\H^{+}/\Z$ has a continuous extension to $\hHZ+$;
\item{(ii)} One has }
$$
\sup_{\H^{+}} |\RE (\sum_{Z}\varphi) |\le C\sup_{\H^{+}}
|\RE\varphi |\; .
$$

\smallskip
We postpone the proof of this Proposition after the statement of the
two following  consequences:

\smallskip
\Proc{Corollary 5.8}{If $\varphi\in\hC{[0,1]}$ then
$\sum_{Z\cM} \varphi $ has a bounded real part which extends
continuously to $\hHZ+$ and verifies}
$$
\sup_{\H^{+}}|\RE \sum_{Z\cM} \varphi |\le C\sup_{\H^{+}}
|\RE\varphi |\; .
$$

\smallskip
\Proc{Corollary 5.9}{ Let $\varphi_{0}(z)=-{1\over \pi}\Li \left(
{1\over z}\right)$; the complex Brjuno function ${\cal B} =
\sum_{Z\cM} \varphi_{0}$ has a
bounded real part which extends continuously to $\hHZ+$.}

\smallskip\noindent
{\it Proof of Corollary 5.9.}
Let $\varphi_{1}$ be defined as in (5.6). Then $\varphi_{1}\in
\hC{[1/2,2]}$, $T\varphi_{1}\in \hC{[0,2/3]}$ and
$\sum_{\cM}T\varphi_{1}\in \hC{[0,1]}$. Thus
$\varphi_{1}+\sum_{\cM}T\varphi_{1}\in \hC{[0,2]}$.
Applying (5.2) to $\varphi_{0}$ we get the desired result.
\qed

\smallskip\noindent
{\it Proof of Proposition 5.7.}
It is not restrictive to assume $I=[0,1]$.

Let us consider the holomorphic function $F_{0}$ associated to
the hyperfunction
$$
u(x) = \cases{ x & if $0\le x\le 1/2\; , $\cr
               1-x & if $1/2\le x\le 1\; . $\cr}
$$
One has
$$
F_{0}(z) = {1\over\pi}\left[z\log{z-1/2\over z}+(1-z)\log{z-1\over
z-1/2}\right]
$$
with the principal branch of the logarithm: if $z\notin [0,1]$,
$(z-1/2)/z$ and $(z-1)/(z-1/2)$ do not belong to $[-\infty ,0]$.

One can easily check that $F_{0}\mid_{\H^\pm}$ extends continuously
to $\overline{\H^\pm}$. At infinity one has (with
$\tilde{z}=z-1/2$):
$$
\eqalign{
z\log{z-1/2\over z} &= \left({1\over 2}+\tilde{z}\right)\log{\tilde{z}\over
\tilde{z}+1/2} \cr
&= -{1\over 2}-{1\over 8\tilde{z}}+\hbox{O}\,(\tilde{z}^{-2})\; , \cr
(1-z)\log{z-1\over z-1/2} &= \left({1\over 2}-\tilde{z}\right)
\log{\tilde{z}-1/2\over\tilde{z}} \cr
&= {1\over 2} - {1\over 8\tilde{z}}+\hbox{O}\,(\tilde{z}^{-2})\; , \cr
}
$$
thus
$$
F_{0}(z) = -{1\over 4\pi z} +\hbox{O}\,(z^{-2})\;
$$
(this could also be checked  directly by
observing that $\int_{0}^{1}u(x)dx=1/4$).

It is easy to verify that:
\item{$\bullet$} the real part of $F_{0}$ is bounded in
$\overline{\C}\setminus [0,1]$;
\item{$\bullet$} the real part of $\sum_{Z}F_{0}$ is
bounded in $\H$.

Let $\varphi\in\hC{I}$; consider the unique decomposition
$$
\varphi = c(\varphi ) F_{0}+\varphi^\sharp
$$
where $\varphi^\sharp\in\O{2}{[0,1]}\cap \hC{[0,1]}$,
$c(\varphi )\in \R$. One has
$
\Vert\varphi^\sharp\Vert_{\hat {\cal C}}\le C
\Vert\varphi\Vert_{\hat{\cal C}}
$
and
$
|c(\varphi )|\le C\Vert\varphi\Vert_{\hat{\cal C}}
$
i.e. this decomposition is continuous.

Thus we are lead to consider only the case $\varphi\in
\O{2}{[0,1]}\cap \hC{[0,1]}$. But in this case
$\sum_{Z}\varphi$ converges uniformly on all domains
$\{|\RE z|\le A\}$  and this fact immediately
leads to the assertions of the Proposition. \qed

\smallskip
{\bf 5.2.9} Note that the topology induced by
$\hHZ+$ on $\H^{+}/\Z\sqcup (\R\setminus\Q )/\Z$ is the same as the
topology induced by $\C/\Z$. Therefore the
 continuity of $\RE \sum_{Z\cM} \varphi_{0}$
on $\hHZ+$ implies that the real part $\RE \sum_{Z\cM} \varphi_{0}$
of the complex Brjuno function
is continuous on $\H^{+}/\Z\sqcup (\R\setminus\Q )/\Z$
in the usual sense. The value $\RE \sum_{Z\cM} \varphi_{0}
(\alpha_{0}, \pi/2)$ (resp. $(\alpha_{0},-\pi/2 )$), with
$\alpha_{0}\in\Q/\Z$, is the right (resp. left) limit of
$\RE \sum_{Z\cM} \varphi_{0} (\alpha )$, as $\alpha\in (\R\setminus\Q
)/\Z$ tends to $\alpha_{0}$

Recalling Lemma 5.2, one has
$$
\RE\varphi_{1} (1,\pi/2 )-\RE\varphi_{1}(1,-\pi/2 ) =-\pi
$$
and more precisely
$$
\RE\varphi_{1} (1,\theta ) = \RE\varphi_{1}(1,0)-\theta\; .
$$
If $\alpha_{0}\in \Q$, $\alpha_{0}\not= 1$ then
$$
\RE\varphi_{1} (\alpha_{0},\theta ) = \RE\varphi_{1}(\alpha_{0},0)
$$
for all $\theta\in [-\pi/2,\pi/2]$. Thus by (5.2) one
obtains that for all $p/q\in\Q$ ($p\wedge q =1$)
$$
\RE \sum_{Z\cM} \varphi_{0} (p/q,\theta ) =
\RE\sum_{Z\cM} \varphi_{0} (p/q, 0)-\theta /q \eqno(5.15)
$$
Thus {\it the  real part  $\RE\sum_{Z\cM} \varphi_{0} $ of the complex
Brjuno function has
at each rational $p/q\in\Q/\Z$ a decreasing jump of $\pi/q$}.

\vskip .5 truecm

\beginsection 5.3 Boundary behaviour of the imaginary part
 of the complex Brjuno function\par
\vskip .3 truecm
{\bf 5.3.0. Notations.}
We simply denote by $\Vert\varphi\Vert$ the norm
(5.13) in the Banach space $E(I)$.

We set $\varphi_{0}(z) = -{1\over\pi}\Li\left( {1\over z}\right)$ as
in (5.4) and $\varphi_{1}$ as in (5.6), thus by (5.2)
$$
{\cal B} = \sum_{Z}(\varphi_{1}+\sum_{\cM}T\varphi_{1})
= \sum_{Z\cM}\varphi_{0}\; . \eqno(5.16)
$$
We have $\varphi_{1}\in E([1/2,2])$, $T\varphi_{1}\in E([0,2/3])$
and $T^k\varphi_{1}\in E([0,1])$ for all $k\ge 2$.

In this Section we want to estimate the imaginary part of the $1$--periodic
function
${\cal B}$ near the real axis. We have
$$
\{|\IM z|\le 1/2\} = \cup_{n\in \Z}(\Delta +n)\; . \eqno(5.17)
$$
For $r\ge 1$, $m_{1},\ldots ,m_{r}\ge 1$, we recall the definition of
$D (m_{1},\ldots ,m_{r})$ given in Section 4.2, namely  the set of $z_0\in D$
such that the continued fraction is
$$
z_{i}^{-1} = z_{i+1}+m_{i+1}\; , \;\; 0\le i<r\; , \;\;
z_{i}\in D\; , \;\; z_{r}\in \Delta \; . \eqno(5.18)
$$
We also set
$$
\eqalign{
H &= H_{0}\cup H_{0}' = \Delta \setminus \hbox{int}\, D\cr
H(m_{1},\ldots ,m_{r}) &= D (m_{1},\ldots ,m_{r})
\setminus \hbox{int}\,\left( \cup_{m_{r+1}\ge 1}D (m_{1},\ldots
,m_{r+1})\right)\; . \cr
}\eqno(5.19)
$$
Then we have
$$
\{|\IM z|\le 1/2\} = \cup_{n\in \Z}\cup_{r\ge 0}
\cup_{m_{1},\ldots ,m_{r}\ge  1}\left[
H(m_{1},\ldots ,m_{r})+n\right]\sqcup \R\setminus\Q\; ,\eqno(5.20)
$$
where the sets in the right--hand term have disjoint interiors.

\smallskip
{\bf 5.3.1}
A set $H(m_{1},\ldots ,m_{r})+n$ meets $\R$ in a unique point, which
belongs to $\Q$. Conversely, any rational belongs to exactly two
such sets: if $p/q\in \Q$ has continued fraction
$p/q= n+1/m_{1}+1/m_{2}+\ldots +1/m_{r}$ with $m_{r}\ge 2$ when
$q>1$ (i.e. $r>0$), these two sets are
$H+n$ and $H(1)+n-1$,  if $p/q=n$ and $r=0$,
$H(m_{1},\ldots ,m_{r})+n$ and $H(m_{1},\ldots ,m_{r}-1,1)+n$
if $r>0$.
The union of these two sets will be denoted by $V(p/q)$; the boundary
of $V(p/q)\cap\H^{+}$ is formed by parts of the three horocycles,
attached to $p/q$, $p'/q'$ and $p''/q''$
(where $[p'/q',p''/q'']$ is the Farey interval with ``center''
$p/q$; when $p/q=n\in\Z$ we have $p'/q'=n-1$, $p''/q''=\infty$), which
are deduced from $\IM z=1/2$ by the action of $\hbox{SL}\, (2,\Z)$
(see Figure 4.).

\smallskip
{\bf 5.3.2}
We plan to compare, when $z\in V(p/q)$, the imaginary part of
${\cal B} (z)$ to the truncated real Brjuno function
$$
B_{finite}(p/q) = \sum_{j=0}^{r-1}\beta_{j-1}(p/q)\log\left[
A^j(p/q-n)\right]^{-1}\; , \eqno(5.21)
$$
where $p/q$ is as above and $A$ is the Gauss map (A1.1). The point is that
we want the dependence on $p/q$ to be explicit in this comparison
(i.e. all the constants $c$ are {\it independent} of $p/q$).
This will be achieved by the following

\smallskip
\Proc{Theorem 5.10}{ For $k\ge 0$, $m_{1},\ldots ,m_{k}\ge 1$,
$z_{0}\in H(m_{1},\ldots ,m_{k})$ one has }
$$
\IM {\cal B}(z_{0}) = B_{finite}(p_{k}/q_{k})+(p_{k-1}-q_{k-1}\RE
z_{0})\IM\varphi_{1}(z_{k}+1)+\hat{r}(z_{0})
$$
{\sl with}
$$
|\hat{r}(z_{0})|\le cq_{k}^{-1}|z_{k}|\log (1+|z_{k}|^{-1})
$$
{\sl for some positive constant $c$ independent of $k,m_{1},\ldots ,m_{k}$.}

\smallskip
The strategy of the proof of the theorem is to start from the formula
$$
\IM B(z)=\sum_{n\in {\Bbb Z}} \sum_{k\ge 0}
\IM T^K\varphi_1(z+n)\;,
$$
and to consider each term separately, putting most of them in the
remainder part, and replacing the others by simpler expressions.
What makes the proof a little lengthy
is that while $T^k\varphi_{1}
\in E([0,1])$ for $k>1$, we have
$\varphi_{1}\in E([1/2,2])$ and
$T\varphi_{1}\in E([0,2/3])$, which leads to distinguish several cases
(Sections 5.3.3 through 5.3.7) before giving its complete
proof (Section 5.3.9).

{\bf 5.3.3} Recall that $T^k\varphi_{1}$ is real,
$T^k\varphi_{1} (\overline{z})= \overline{T^k\varphi_{1}(z)}$,
hence $\IM T^k\varphi_{1} (z)$ vanishes on $\R$ outside of
$$
\cases{
[1/2,2] & for $k=0\; ,$ \cr
[0,2/3] & for $k=1\; ,$ \cr
[0,1] & for $k>1\; .$ \cr
}
$$
As we have, on the other hand,
$$
| \varphi' (z) |\le C_{K}|z|^{-2}\Vert\varphi\Vert\eqno(5.22)
$$
for all $\varphi\in E(I)$, $K\subset\overline{\C}\setminus I$
compact, $z\in K$, we obtain for $z\in \Delta$
$$
\eqalign{
|\IM \varphi_{1}(z+n )| &\le cn^{-2}|\IM z|\Vert\varphi_{1}\Vert \; ,
\; \; \hbox{if}\, n\not= 0,1,2\; , \; n\in \Z\; , \cr
|\IM T\varphi_{1}(z+n )| &\le cn^{-2}|\IM z|\Vert T\varphi_{1}\Vert \; ,
\; \; \hbox{if}\, n\not= -1,0\; , \; n\in \Z\; , \cr
|\IM T^k\varphi_{1}(z+n )| &\le cn^{-2}|\IM z|\Vert T^k\varphi_{1}\Vert \; ,
\; \; \hbox{if}\, k>1\; , \; n\not= -1,0,1\; , \; n\in \Z\; . \cr
}\eqno(5.23)
$$
Since, by Proposition 4.11,
we have $\Vert T^k\varphi_{1}\Vert \le CkG^{-k}\Vert
\varphi_{1}\Vert $, with $G=(\sqrt{5}+1)/2$, we obtain

\smallskip
\Proc{Lemma 5.11}{For $z\in \Delta \setminus [0,1]$, if we write
$$
\eqalignno{
\IM {\cal B}& (z) = \IM \varphi_{1}(z) +\IM \varphi_{1}(z+1) + \IM
\varphi_{1}(z+2 )\cr
& +\IM T\varphi_{1}(z-1) +\IM T\varphi_{1}(z)\cr
 +\sum_{k>1}&\left[
\IM T^k\varphi_{1}(z-1)+ \IM T^k \varphi_{1}(z)+\IM T^k
\varphi_{1}(z+1)\right]
+r_{0}(z)\;  ,&(5.24) \cr
}
$$
then we have}
$$
|r_{0}(z)| \le C |\IM z|\Vert\varphi_{1}\Vert\; . \eqno(5.25)
$$

\smallskip
{\bf 5.3.4} One has
$$
| \varphi_{1}'(z)|\le C\log (1+|z-1/2|^{-1})
$$
in the neighborhood of $1/2$ and
$$
| \varphi_{1}'(z)|\le C\log (1+|z-2|^{-1})
$$
in the neighborhood of $2$. Therefore, for $z\in\Delta \setminus\{0\}$
we have
$$
|\IM \varphi_{1}(z+2)|\le C|\IM z|\log (1+|z|^{-1})\; , \eqno(5.26)
$$
and, for $z\in \Delta\setminus D(1)$
$$
|\IM \varphi_{1}(z)|\le C|\IM z|\log (1+|z-1/2|^{-1})\; .
\eqno(5.27)
$$
Next we have

\smallskip
\Proc{Lemma 5.12}{For $k\ge 1$ and $z\in \Delta \setminus\{1\}$ we
have}
$$
|\IM T^k \varphi_{1}(z-1)|\le C|\IM z|\log (1+|z-1|^{-1})\Vert
T^{k-1}\varphi_{1}\Vert \; .
$$

\smallskip
\proof
We have
$$
|\IM T^k\varphi_{1}(z) | = \left|\IM\left[ z\sum_{m\ge 1}\left(
T^{k-1}\varphi_{1}\left({1\over
z}-m\right)-T^{k-1}\varphi_{1}(-m)\right)\right]
\right|\; .
$$
We distinguish two cases:
\item{(a)} If $|\IM z^{-1}|=|z|^{-2}|\IM z|\ge 1$, we choose
$m_{0}\ge 1$ such that
$$
\left| {1\over z}+m_{0}\right|\le C|z|^{-2} |\IM z|\;\;  \hbox{and}\;
C^{-1}|z|^{-1} \le m_{0}\le C|z|^{-1}\; .
$$
One has then
$$
\left| T^{k-1}\varphi_{1}\left({1\over z}-m\right)-T^{k-1}\varphi_{1}(-m-m_{0})
\right| \le C {|\IM z| \Vert T^{k-1}\varphi_{1}\Vert \over
|z|^{2} (m+m_{0})^{2}}
$$
from which it follows that
$$
\left|\left[ z\sum_{m\ge 1}\left(
T^{k-1}\varphi_{1}\left({1\over z}-m\right)-
T^{k-1}\varphi_{1}(-m-m_{0})\right)\right]
\right|\le C|\IM z|\Vert T^{k-1}\varphi_{1}\Vert\; .
$$
On the other hand
$$
|\IM z\sum_{1}^{m_{0}}T^{k-1}\varphi_{1}(-m)|\le C|\IM z|\log
(1+|z|^{-1})\Vert T^{k-1}\varphi_{1}\Vert\; ,
$$
from which the lemma follows in this case.
\item{(b)} If $|\IM z^{-1}|\le 1$ we choose $m_{0}\ge 1$ such that
$-m_{0}\le \RE 1/z\le -m_{0}+1$ thus
$$
\eqalign{
\left|\IM T^{k-1}\varphi_{1}\left({1\over z}-m\right)\right|
\le C|\IM z| |z|^{-2} &(m+m_{0})^{-2}\Vert T^{k-1}\varphi_{1}\Vert\; , \cr
\left|\RE\left[ T^{k-1}\varphi_{1}\left({1\over z}-m\right) -
T^{k-1}\varphi_{1}
(-m-m_{0})\right]\right|
&\le C(m+m_{0})^{-2}\Vert T^{k-1}\varphi_{1}\Vert\; . \cr
}
$$
We thus obtain
$$
\eqalign{
|\RE z|\left|\IM\sum_{m\ge 1}T^{k-1}\varphi_{1}\left({1\over z}-m\right)\right|
\le C|\IM z|&\Vert T^{k-1}\varphi_{1}\Vert\; , \cr
\left| \RE\sum_{m\ge 1}
\left[ T^{k-1}\varphi_{1}\left({1\over z}-m\right) - T^{k-1}\varphi_{1}
(-m-m_{0})\right]\right|
\le C|z|&\Vert T^{k-1}\varphi_{1}\Vert\; , \cr
%|\IM z|
|\sum_{1}^{m_{0}}T^{k-1}\varphi_{1}(-m)|
\le C
%|\IM z|
\log (1+|z|^{-1})&\Vert T^{k-1}\varphi_{1}\Vert\; ,\cr
}
$$
which give the desired result. \qed

\smallskip
\Proc{Lemma 5.13}{}

\item{1.} {\sl If $z\in D(1) \setminus D(1,1)$ we have}
$$
|\IM T\varphi_{1}(z) |\le C|\IM z|\log (1+|z-2/3|^{-1})\; .
$$
\item{2.} {\sl For $k>1$, $z\in \Delta\setminus\{0\}$ we have}
$$
|\IM T^k \varphi_{1}(z+1)|\le C|\IM z|\log (1+|z|^{-1})\Vert
T^{k-2}\varphi_{1}\Vert\; .
$$

\smallskip
\proof
In the first case,
in the domain considered, we have
$ 1/z -1\in \Delta\setminus D(1)$ hence, by (5.27),
$$
\left|\IM \varphi_{1}\left({1\over z}-1\right)\right|\le
C|\IM z|\log (1+|z-2/3|^{-1})\; .
$$
On the other hand we have
$$
\left|\RE \left[\varphi_{1}\left({1\over z}-1\right)- \varphi_{1}(-1)
\right]\right|\le 2\Vert\varphi_{1}\Vert\; ,
$$
and for $m\ge 2$
$$
\eqalign{
\left|\varphi_{1}\left({1\over z}-m\right) - \varphi_{1}(-m)\right| &\le
Cm^{-2}\Vert\varphi_{1}\Vert\; , \cr
\left|\IM \varphi_{1}\left({1\over z}-m\right)\right| &\le
Cm^{-2}|\IM z|
\Vert\varphi_{1}\Vert\; , \cr
}
$$
from which the first inequality of the Lemma follows.

In the second case we use Lemma 5.12 to get
$$
\left| \IM T^{k-1}\varphi_{1}\left({1\over z+1}-1\right) \right|
\le C |\IM z| \log (1+|z|^{-1}) \Vert T^{k-2}\varphi_{1}\Vert\; ,
$$
and deduce the second inequality as above.
\qed
\smallskip

{\bf 5.3.5} For $z\in \Delta$, $k\ge 1$, we set (as in Section 4.2)
$\varepsilon_{k}=
\varepsilon_{k} (z) =0$ if $z$ belongs to some $D (m_{1},\ldots ,
m_{k})$ with $m_{k}=1$, $\varepsilon_{k} =1$ otherwise.

Starting from Lemma 5.11, we use (5.26) to deal with
$\IM \varphi_1(z+2)$, (5.27) to deal with $\IM\varphi_1(z)$
when $\varepsilon=1$, Lemma 5.12 to deal
with $\IM T^k\varphi_1(z-1)$, Lemma 5.13 to deal with
$\IM T^k\varphi_1(z+1)$, and also with $\IM T\varphi_1(z)$ when
$\varepsilon_1=0$, $\varepsilon_2=1$.
This gives

\smallskip
\Proc{Lemma 5.14}{ For $z\in\Delta\setminus [0,1]$ we have
$$
\eqalign{
\IM {\cal B}(z) &= (1-\varepsilon_{1})\IM \varphi_{1}(z)+\IM \varphi_{1}
(z+1)+[1-\varepsilon_{2}(1-\varepsilon_{1})]\IM T\varphi_{1}(z)\cr
&+ \sum_{k>1}\IM T^k\varphi_{1}(z)+r_{0}(z)+r_{1}(z)\cr
}
$$
and }
$$
|r_{1}(z)|\le C|\IM z|\log (1+|\IM z|^{-1})\Vert\varphi_{1}\Vert\; .
$$

\smallskip
We next recall that, by Proposition 4.1, we have, when $z\in H$,
$k\ge 1$
$$
|\IM T^k\varphi_{1}(z)|\le C|z|\log (1+|z|^{-1})\sup_{D_{\infty}}
|T^{k-1}\varphi_{1}|\; . \eqno(5.28)
$$

\smallskip
{\bf 5.3.6} In the following two steps, we
obtain for $\IM T^k\varphi_1(z)$ an approximation
similar to Proposition 4.4, first for large $k$ (Proposition 5.15),
and then for small $k$ (Proposition (5.16). Here we assume that $k\ge 1$,
$m_{1},\ldots ,m_{k}\ge 1$ and $z_{0}\in D (m_{1},\ldots
,m_{k})$ and let
$l>0 $ and $\varphi = T^l\varphi_{1}$.

\smallskip
\Proc{Proposition 5.15}{ For $k\ge 1$, $m_{1},\ldots ,m_{k}\ge 1$,
$z_{0}\in D (m_{1},\ldots ,m_{k})$, $l>0$ we have }
$$
\eqalign{
|\IM & T^{k+l}\varphi_{1}(z_{0})-[\IM T^l\varphi_{1}(z_{k})]
(p_{k-1}-q_{k-1}\RE z_{0})| \cr
&\le Cq_{k}^{-1}|\IM z_{k}|\log (1+|\IM z_{k}|^{-1})\times\cases{
\Vert T^{l-2}\varphi_{1}\Vert &if $l>1$, \cr
\Vert \varphi_{1}\Vert &if $l=1$. \cr
}\cr
}
$$

\smallskip
\proof
As we did in the proof of Proposition 4.4 we write
$$
\eqalign{
T^k\varphi (z_{0}) &=
(p_{k-1}-q_{k-1}z_{0})\varphi (z_{k})+\sum_{j=1}^k
(p_{j-1}-q_{j-1}z_{0})T^{k-j}\varphi (z_{j}-1)\cr
&+ \sum_{j=1}^k
(p_{j-1}-q_{j-1}z_{0})\varepsilon_{j}T^{k-j}\varphi (z_{j}+1)\cr
&+\sum_{j=1}^k
(p_{j-2}-q_{j-2}z_{0})R^{(m_{j})}(T^{k-j}\varphi )(z_{j-1})\cr
&= (p_{k-1}-q_{k-1}z_{0})\varphi (z_{k})+ \hat{R}_{1}^{[k]} (\varphi )
(z_{0})+\hat{R}_{2}^{[k]} (\varphi )
(z_{0})+\hat{R}_{0}^{[k]} (\varphi )
(z_{0})\cr
}
\eqno(5.29)
$$
where
$$
\eqalign{
R^{(m_{j})}(\psi)(z_{j-1}) &= \sum_{m\ge 1} \psi'
(-m)+z_{j-1}\sum_{m\ge 1\, , \,|m-m_{j}|\le 1}\psi (-m)\cr
&-z_{j-1}\sum_{m\ge 1\, , \,|m-m_{j}|> 1}[\psi (z_{j}+m_{j}-m)
-\psi (-m)]\; . \cr
}
\eqno(5.30)
$$
We will use repeatedly the following inequalities
$$
\eqalignno{
|\IM (w_{1}w_{2})| &\le |w_{1}||\IM w_{2}|+|w_{2}||\IM w_{1}|\; ,
&(5.31)\cr
|p_{j-1}-q_{j-1}z_{0}| &\le cq_{j}^{-1}\; , &(5.32)\cr
|\IM (p_{j-1}-q_{j-1}z_{0})| &= q_{j-1}|\IM z_{0}|\le Cq_{j-1}
q_{k}^{-2}|\IM z_{k}|\; . &(5.33)\cr
}
$$
By Proposition 4.1, as $l>0$, for $1\le j\le k$ we have
$$
|T^{k-j}\varphi (z_{j}-1)|\le
C\sup_{D_{\infty}}|T^{k-j+l-1}\varphi_{1}|\; , \eqno(5.34)
$$
and, by Lemma 5.12,
$$
|\IM T^{k-j}\varphi (z_{j}-1)|\le C|\IM z_{j}|\log (1+|z_{j}-1|^{-1})
\Vert T^{k-j+l-1}\varphi_{1}\Vert \; . \eqno(5.35)
$$
We observe that for $j<k-1$
$$
\log (1+|z_{j}-1|^{-1})\le C\log (1+|z_{j+1}|^{-1})\le C\log
(1+m_{j+2})\; ,
$$
and one also has
$$
\eqalign{
\log (1+|z_{k-1}-1|^{-1})&\le C\log (1+|z_{k}|^{-1})\le C\log
(1+|\IM z_{k}|^{-1})\; , \cr
\log (1+|z_{k}-1|^{-1})&\le \log (1+|\IM z_{k}|^{-1})\; . \cr
}
$$
Thus, from (5.31)--(5.35) we get
$$
|\IM \hat{R}_{1}^{[k]}(\varphi )(z_{0})|\le
Cq_{k}^{-1}|\IM z_{k}| \log (1+|\IM z_{k}|^{-1})\Vert T^{l-1}\varphi_{1}
\Vert\; .
$$
Similarly, by Lemma 5.13, we have, for $j<k$, $l\ge 1$
or $j\le k$, $l>1$
$$
|\IM T^{k-j}\varphi (z_{j}+1)|\le C|\IM z_{j}|\log (1+|z_{j}|^{-1})
\Vert T^{k-j+l-2}\varphi_{1}\Vert \; , \eqno(5.36)
$$
and, by Proposition 4.2,
$$
|T^{k-j}\varphi (z_{j}+1)|\le
C\sup_{D_\infty}|T^{k-j+l-2}\varphi_{1}|\; . \eqno(5.37)
$$
On the other hand,  when $j=k$, $l=1$ (i.e. $T^{k-j}\varphi =T\varphi_{1}$)
one has
$$
\eqalign{
|\IM T\varphi_{1}(z_{k}+1)| &\le C|\IM z_{k}|\Vert
T\varphi_{1}\Vert\; , \cr
|T\varphi_{1}(z_{k}+1)| &\le C\Vert T\varphi_{1}\Vert\; , \cr
}
$$
(note that as $T\varphi_{1}\in E([0,2/3])$, $z_{k}+1$ is well separated from
the boundary). This gives
$$
|\IM \hat{R}_{2}^{[k]}(\varphi )(z_{0})|\le Cq_{k}^{-1}|\IM z_{k}|
\log (1+|z_{k}|^{-1})\times \cases{\Vert T^{l-2}\varphi_{1}\Vert &
$(l>1)$\cr
\Vert\varphi_{1}\Vert & $(l=1)$\cr}\; . \eqno(5.38)
$$
Finally one has
$$
\eqalignno{
|R^{(m_{j})}(\psi)(z_{j-1})| &\le C\Vert\psi\Vert\; , &(5.39)\cr
|\IM R^{(m_{j})}(\psi)(z_{j-1})| &\le C\Vert\psi\Vert
[|\IM z_{j-1}|\log (1+m_{j})+\RE z_{j-1}|\IM z_{j}|]\cr
&\le cm_{j}^{-1}|\IM z_{j}|\Vert\psi\Vert\; , &(5.40)\cr
}
$$
giving
$$
|\IM R_{0}^{[k]}(\varphi )(z_{0})|\le cq_{k}^{-1}|\IM z_{k}|\Vert T^l
\varphi_{1}\Vert \; . \eqno(5.41)
$$
To conclude our proof we only need to observe that
$$
|\IM (p_{k-1}-q_{k-1}z_{0})\RE \varphi (z_{k})|\le C q_{k}^{-1}|\IM z_{k}|
\Vert T^l\varphi_{1}\Vert \; , \eqno(5.42)
$$
to get the desired result. \qed

\smallskip
{\bf 5.3.7} For $k\ge 1$, $m_{1},\ldots ,m_{k}\ge 1$, $z_{0}\in
D(m_{1},\ldots ,m_{k})$ we now consider $\IM T^k \varphi_{1}(z_{0})$
(i.e. the case $l=0$ left out from Proposition 5.15):

\smallskip
\Proc{Proposition 5.16}{For $k\ge 1$, $m_{1},\ldots ,m_{k}\ge 1$,
$z_{0}\in D(m_{1},\ldots ,m_{k})$ we have }
$$
\eqalign{
|\IM &T^k\varphi_{1}(z_{0})-(p_{k-1}-q_{k-1}\RE z_{0})[
(1-\varepsilon_{k+1})\IM\varphi_{1}(z_{k})+\varepsilon_{k}\IM
\varphi_{1}(z_{k}+1)]|\cr
&\le Cq_{k}^{-1}|\IM z_{k}|\log (1+|\IM
z_{k}|^{-1})\Vert\varphi_{1}\Vert \; . \cr
}
$$

\smallskip
\proof
We now write (with $\varepsilon_{k}'=0$ if $m_{k}\le 2$,
$\varepsilon_{k}'=1$ otherwise)
$$
\eqalign{
T^k \varphi_{1}(z_{0}) &= (p_{k-1}-q_{k-1}z_{0})(\varphi_{1}(z_{k})+
\varepsilon_{k}\varphi_{1}(z_{k}+1)+\varepsilon_{k}'\varphi_{1}
(z_{k}+2))\cr
&+ \sum_{j=1}^{k-1}(p_{j-1}-q_{j-1}z_{0})T^{k-j}\varphi_{1}(z_{j}-1)\cr
&+ \sum_{j=1}^{k-1}(p_{j-1}-q_{j-1}z_{0})\varepsilon_{j}T^{k-j}\varphi_{1}
(z_{j}+1)\cr
&+ \sum_{j=1}^{k}(p_{j-2}-q_{j-2}z_{0})
\tilde{R}^{(m_{j})}(T^{k-j}\varphi_{1})(z_{j-1})\cr
}
\eqno(5.43)
$$
where $\tilde{R}^{(m_{j})}(T^{k-j}\varphi_{1}) =
R^{(m_{j})}(T^{k-j}\varphi_{1})$ for $j<k$ but
$$
\eqalign{
\tilde{R}^{(m_{k})}(\varphi_{1})(z_{k-1}) &= \sum_{m\ge 1} \varphi_{1}'
(-m)+z_{k-1}\sum_{m\ge 1\, , \, |m-m_{k}+1|\le 1}\varphi_{1}(-m) \cr
&-z_{k-1}\sum_{m\ge 1\, , \, |m-m_{k}+1|>
1}[\varphi_{1}(z_{k}+m_{k}-m)-\varphi_{1}(-m)]\; . \cr
}
$$
For the last sum $\tilde{R}_{0}^{[k]} (\varphi_{1})$ in (5.43),
we have as in (5.39)--(5.41)
$$
|\IM \tilde{R}_{0}^{[k]} (\varphi_{1})(z_{0})|\le Cq_{k}^{-1}|\IM
z_{k}|\Vert\varphi_{1}\Vert \; . \eqno(5.44)
$$
The first sum $R_{1}^{[k]}(\varphi_{1})$ is dealt with as in
(5.31)--(5.35) to get
$$
|\IM R_{1}^{[k]} (\varphi_{1})(z_{0})|\le Cq_{k}^{-1} |\IM z_{k}|\log
(1+|\IM z_{k}|^{-1})\Vert\varphi_{1}\Vert\; . \eqno(5.45)
$$
The middle sum $R_{2}^{[k]}(\varphi_{1})$ satisfies the same
estimate, proved as in (5.36)--(5.38).

By (5.26) we have
$$
|\IM [(p_{k-1}-q_{k-1}z_{0})\varphi_{1}(z_{k}+2)]|\le
q_{k}^{-1}
[C|\IM z_{k}|\log (1+|z_{k}|^{-1})+1]\Vert\varphi_{1}\Vert
\; . \eqno(5.46)
$$
Similarly, if $z_{k}\notin D(1)$ (i.e. if $\varepsilon_{k+1}=1$) we
have
$$
|\IM [(p_{k-1}-q_{k-1}z_{0})\varphi_{1}(z_{k})]|\le
Cq_{k}^{-1}|\IM z_{k}|\log (1+|z_{k}-1/2|^{-1})\Vert\varphi_{1}\Vert
\; , \eqno(5.47)
$$
according to (5.27). Finally, for $\varepsilon =0,1$
$$
|\IM (p_{k-1}-q_{k-1}z_{0})\RE \varphi_{1}(z_{k}+\varepsilon)|\le
Cq_{k}^{-1}|\IM z_{k}|\Vert\varphi_{1}\Vert
\; , \eqno(5.48)
$$
and from (5.43)--(5.48) we get our result. \qed

\smallskip
{\bf 5.3.8} Starting from Lemma 5.14, we now make use of
Propositions 5.15 and 5.16 to obtain a simpler approximation
for $\IM {\cal B}(z)$.

\smallskip
\Proc{Proposition 5.17}{For $k\ge 0$, $m_{1},\ldots ,m_{k}\ge 1$,
$z_0\in H(m_{1},\ldots ,m_{k})$ we have
$$
\eqalign{
\IM {\cal B}(z_{0}) &=
\sum_{l=0}^k\IM \varphi_{1}(z_{l}+1)(p_{l-1}-q_{l-1}\RE
z_{0})\varepsilon_{l}\cr
&+\sum_{l=0}^{k-1}\IM \varphi_{1}(z_{l})(p_{l-1}-q_{l-1}\RE
z_{0})(1-\varepsilon_{l+1})+r(z_{0})\cr
}
$$
with $\varepsilon_0=1$, and
$$
|r(z_{0})|\le Cq_{k}^{-1}|z_{k}|\log (1+|z_{k}|^{-1})\; .
$$}

\smallskip
\proof
First, assume $z_{0}\in H$. Then, from (5.28) and Lemma 5.14 we get
$$
\IM {\cal B}(z) = (1-\varepsilon_{1})\IM\varphi_{1}(z_{0})+
\IM\varphi_{1}(z_{0}+1) +r(z_{0})\; , \eqno(5.49)
$$
with
$$
|r(z_{0})|\le C|z_{0}|\log (1+|z_{0}|^{-1})\; . \eqno(5.50)
$$
As we also have, for $z_{0}\in H$,
$$
|\IM \varphi_{1}(z_{0})|\le C|\IM z_{0}|\; ,
$$
we obtain
$$
|\IM {\cal B}(z)-\IM \varphi_{1}(z_{0}+1)|\le C|z_{0}|\log (1+|z_{0}|^{-1})\; .
\eqno(5.51)
$$
Next consider $z_{0}\in H(m_{1},\ldots ,m_{k})$ with $k\ge 1$,
$m_{1},\ldots ,m_{k}\ge 1$. Observe that
$$
[1-\varepsilon_{2}(1-\varepsilon_{1})](1-\varepsilon_{2})=
1-\varepsilon_{2}\; , \;\;\; [1-\varepsilon_{2}(1-\varepsilon_{1})]
\varepsilon_{1}=\varepsilon_{1}\; .
$$
In the terms which appear in the right hand term of Lemma 5.14, we
use Proposition 5.15 to deal with $\IM T^l\varphi_{1}(z_{0})$,
$l>k$ and Proposition 5.16 to deal with $\IM T^l\varphi_{1}(z_{0})$,
$0<l\le k$. We have
$$
\eqalign{
|\sum_{l>k}\IM T^l\varphi_{1}(z_{0}) &- (p_{k-1}-q_{k-1}z_{0})
\sum_{l>0}\IM T^l\varphi_{1}(z_{k})|\cr
&\le Cq_{k}^{-1}|\IM z_{k}|\log (1+|\IM z_{k}|^{-1})\; , \cr
}
\eqno(5.52)
$$
but using (5.28) (since $z_{k}\in H$) we obtain
$$
|\sum_{l>k}\IM T^l\varphi_{1}(z_{0})|\le
Cq_{k}^{-1}|z_{k}|\log (1+| z_{k}|^{-1})\; . \eqno(5.53)
$$
For $0\le l\le k$, we have in Proposition 5.16 and Lemma 5.14
$$
q_{l}^{-1}|\IM z_{l}|\log (1+|\IM z_{l}|^{-1})\le
Cq_{l}q_{k}^{-2}|\IM z_{k}|\log (1+|\IM z_{l}|^{-1})\; ,
$$
and
$$
\eqalign{
\sum_{0}^k q_{l}q_{k}^{-1}\log (1+|\IM z_{l}|^{-1})&\le
C\sum_{0}^k q_{l}q_{k}^{-1}\log (1+Cq_{l}^{-2}q_{k}^{2}|\IM
z_{k}|^{-1})\cr
&\le  C \log (1+|\IM z_{k}|^{-1})\; . \cr
}
$$
\qed

\smallskip
{\bf 5.3.9} We can finally complete the {\sl proof of Theorem 5.10}.

Let $k\ge 0$, $m_{1},\ldots ,m_{k}\ge 1$. As we already know, the
domain $H(m_{1},\ldots ,m_{k})$ meets $\R$ in a unique point which is
$p_{k}/q_{k}$. Let us denote $x_{0}=p_{k}/q_{k}$ and consider its
continued fraction
$$
x_{i}^{-1}=m_{i+1}+x_{i+1}\; , \;\; 0\le i<k\; , \;\; x_{k}=0
$$
(of course $x_{i}$ is the point of intersection of $H(m_{i+1},\ldots ,m_{k})$
with $\R$).

Let $z_{0}\in H(m_{1},\ldots ,m_{k})$, $\IM z_{0}>0$.
We will then have
$(-1)^l\IM z_{l}>0$, $0\le l\le k$. On the other hand one has, for
$0\le l<k$
$$
\varepsilon_{l}\IM \varphi_{1}(x_{l}+1+(-1)^li0) =
(-1)^l\varepsilon_{l}\log {1\over x_{l}}\; , \eqno(5.54)
$$
and for $0\le l<k-1$
$$
(1-\varepsilon_{l+1})\IM\varphi_{1}(x_{l}+(-1)^li0) =
(1-\varepsilon_{l+1})(-1)^lx_{l}\log {1\over x_{l+1}} \eqno(5.55)
$$
since ${x_{l}\over 1-x_{l}}={1\over x_{l+1}}$ when
$1-\varepsilon_{l+1}\not= 0$. (5.54) and (5.55) imply that
$$
\eqalign{
\sum_{0}^{k-1} & \IM\varphi_{1}(x_{l}+1+(-1)^li0)(p_{l-1}-q_{l-1}x_{0})
\varepsilon_{l}\cr
&+\sum_{0}^{k-2}\IM \varphi_{1}(x_{l}+(-1)^li0)(p_{l-1}-
q_{l-1}x_{0})(1-\varepsilon_{l+1})\cr
&= \sum_{0}^{k-1} \beta_{l-1}(x_{0})\log {1\over x_{l}} :=
B_{finite}(p_{k}/q_{k})\; . \cr
}\eqno(5.56)
$$
When we set this expression in Theorem 5.10, and compare the result
with expression for $\IM B(z_0)$ given by Proposion 5.17, we see
that we have to deal with the following expressions.
For $0\le l\le k-1$, and $\varepsilon_l=1$:
$$
{\bf A}_l= \IM\varphi_1(z_l+1)(p_{l-1}-q_{l-1}\RE z_0)-
\IM\varphi_1(x_l+1 +(-1)^li0)(p_{l-1}-q_{l-1}x_0)\; .
$$
For $0\le l\le k-2$, and $\varepsilon_{l+1}=0$:
$$
{\bf B}_l=
\IM\varphi_1(z_l)(p_{l-1}-q_{l-1}\RE z_0)-
\IM\varphi_1(x_l +(-1)^li0)(p_{l-1}-q_{l-1}x_0)\; .
$$
And when $\varepsilon_k=0$ :
$$
{\bf C}_l=
\IM\varphi_1(z_k+1)(p_{k-1}-q_{k-1}\RE z_0)-
\IM\varphi_1(z_{k-1})(p_{k-2}-q_{k-2}\RE z_0)\; .
$$

First of all one has
$$
|\RE z_{0}-x_{0}|\le Cq_{k}^{-2}|z_{k}|\; .
$$
On the other hand, near $1$ one has
$$
| \varphi_{1}' (z)|\le C|z-1|^{-1}\; .
$$
For $0\le l<k$ (resp. $0\le l< k-1$) the distances of $x_{l}$
and $z_{l}$ from $0$ (resp. $1$) are comparable. Thus, for $0\le l<k$
$$
|\IM \varphi_{1}(z_{l}+1)-\IM \varphi_{1}(x_{l}+1+(-1)^li0)|\le
Cx_{l}^{-1}|z_{l}-x_{l}|\; ,
$$
and, for $0\le l<k-1$
$$
|\IM \varphi_{1}(z_{l})-\IM \varphi_{1}(x_{l}+(-1)^li0)|\le
C|x_{l}-1|^{-1}|z_{l}-x_{l}|\; .
$$
We have here
$|z_{l}-x_{l}|\le C|z_{k}|q_{l}^{2}q_{k}^{-2}$, $x_{l}^{-1}\le
Cm_{l+1}$ and $|x_{l}-1|^{-1}\le cx_{l+1}^{-1}\le cm_{l+2}$.
We thus get, taking (5.54) and (5.55) into account,
$$ \eqalign{|{\bf A}_l|&\le
C(q_{l-1}q_k^{-2}|z_k|\log{1\over x_l}+x_l^{-1}|z_l-x_l|q_l^{-1})
\cr&\le
C|z_k|q_k^{-2}(q_{l-1}\log m_{l+1}+q_lm_{l+1})
\cr&\le
C|z_k|q_k^{-2}q_{l+1}\; ,
\cr}$$
$$ \eqalign{|{\bf B}_l|&\le
C(q_{l-1}q_k^{-2}|z_k|\log{1\over x_{l+1}}
+|x_l-1|^{-1}|z_l-x_l|q_l^{-1})
\cr&\le
C|z_k|q_k^{-2}(q_{l-1}\log m_{l+2}+q_lm_{l+2})
\cr&\le
C|z_k|q_k^{-2}q_{l+2}\; .
\cr}$$
We thus have
$$
\sum_{l=0}^{k-1}|{\bf A}_l|
+\sum_{l=0}^{k-2}|{\bf B}_l|\le C |z_k|q_k^{-1}\; .
$$

Finally we note that when $m_{k}=1$, i.e. $\varepsilon_{k}=0$, one has
$z_{k-1}^{-1}=1+z_{k}$ thus
$$
\varphi_{1}(z_{k-1}) = (1+z_{k-1})({\pi\over 12}+{1\over \pi}\log 2)
-z_{k-1}\varphi_{1}(1+z_{k})\; ,
$$
which gives
$$\eqalign{ {\bf C}&=\IM\varphi_1(z_k+1)[(p_{k-1}-q_{k-1}\RE z_0)
+(p_{k-2}-q_{k-2}\RE z_0)\RE z_{k-1})]\cr
&+(p_{k-2}-q_{k-2}\RE z_{0})\IM z_{k-1}
\left({\pi\over 12}+{1\over \pi}\log 2-\RE\varphi_1(1+z_k)\right)
\; .\cr}$$
But we have the following inequalities
$$
\left|{\pi\over 12}+{1\over \pi}\log 2-\RE\varphi_1(1+z_k)\right|
\le C\; ,
$$
$$
\eqalign{
|(p_{k-2}-q_{k-2}\RE z_{0})\IM z_{k-1}|
&\le Cq_{k-1}^{-1}|\IM z_{k-1}|\cr
&\le Cq_{k}^{-1}|\IM z_{k}|\; ,\cr}
$$
$$\eqalign{
|(p_{k-1}-q_{k-1}\RE z_0)
+(p_{k-2}-q_{k-2}\RE z_0)\RE z_{k-1})|&\le q_{k-2}|\IM z_0||\IM z_{k-1}|\cr
&\le Cq_k^{-1}|\IM z_k| \; ,\cr}
$$
$$
|\IM\varphi_1(z_k+1)|\le C\log(1+|z_k|^{-1}) \; ,
$$
hence
$$
|{\bf C}|\le Cq_k^{-1}|\IM z_k|\, \log(1+|z_k|^{-1})\; ,
$$
and the proof of Theorem 5.10 is complete. \qed

\smallskip
\remark{5.18} We recall that the term involved in Theoorem 5.10
satisfies
$$
(p_{k-1}-q_{k-1}\RE z_0) \IM \varphi_{1}(z_{k}+1) =
\left(q_k^{-1} \log {1\over |z_{k}|}\right)(1+ \hbox{o}\, (1))\; ,
$$
as $z_k\to 0$.

\smallskip
{\bf 5.3.10} Here we consider the imaginary part of ${\cal B}$ near Brjuno
numbers. For $H>0$ let
$$
W_{H}= \{ w\in \H\; , \;\; \IM w \ge |\RE w|^{H}\}\eqno(5.57)
$$
and for $0<h<1/2$ let
$$
\tilde{W}_{h} = \{w\in \H\; , \;\; \IM w\ge \exp [-|\RE
w|^{-h}]\}\; . \eqno(5.58)
$$
Then we have
\smallskip
\Proc{Theorem 5.19}{
\item{1.} For any Brjuno number $\alpha$ and any $H>0$ we have }
$$
\lim_{w\rightarrow 0\, , \, w\in W_{H}}
\IM {\cal B}(w+\alpha) = B(\alpha )\; .
$$
\item{2.}{\sl Let  $\alpha$ be an irrational diophantine number  and
$0<h<1/2$ such
that }
$$
\liminf_{q\rightarrow\infty}\Vert q\alpha\Vert_{\Z} q^{1/h-1} =
+\infty\; ,
$$
{\sl where $\Vert\;\Vert_{\Z}$ denotes the distance from the nearest
integer. Then}
$$
\lim_{w\rightarrow 0\, , \, w\in \tilde{W}_{h}}
\IM {\cal B}(w+\alpha) = B(\alpha )\; .
$$

\smallskip
\proof
We begin by stating a useful Lemma (whose proof is an easy adaptation
of the arguments of 5.3.9 and is left to the reader).

\smallskip
\Proc{Lemma 5.20}{Let $k\ge 1$, $m_{1},\ldots ,m_{k}\ge 1$, $p_{k}/q_{k}$
be the point of intersection of $H(m_{1},\ldots ,m_{k})$ with $\R$.
For all $x\in  D(m_{1},\ldots ,m_{k})\cap\R$ one has}
$$
|B_{finite}(p_{k}/q_{k})-\sum_{0}^{k-1}\beta_{l-1}(x)\log{1\over
x_{l}}|\le Cx_{k}q_{k}^{-1}
$$
{\sl where $(x_{i})_{i\ge 0}$ is the continued fraction of $x$. }

\smallskip
Assume now that $\alpha\in (0,1)$ is irrational with continued
fraction
$$
\alpha = 1/m_{1}+1/m_{2}+\ldots +1/m_{k}+\ldots \; .
$$
Let $(p_{k}/q_{k})_{k\ge 0}$ denote the sequence of the partial
fractions. Let $w$ be a point close to $0$ and $z=\alpha +w$. For
$|\IM z|\le 1/2$, $z$ belongs to a domain $V(p/q)$
(defined  in  Section 5.3.1) and we distinguish two cases.
\item{(I)} Here we assume that $p/q=p_{k}/q_{k}$ is one of the partial
fractions
of $\alpha$. One has then by Proposition A1.1
$$
|\alpha - p_{k}/q_{k}|\ge (2q_{k}q_{k+1})^{-1}\; .
$$
\itemitem{(I.1)} If $w\in W_{H}$ one gets
$$
|z-p_{k}/q_{k}|\ge c^{-1} (q_{k}q_{k+1})^{-H}\; ,
$$
thus
$$
|z_{k}|^{-1}\le cq_{k}^{-2}(q_{k}q_{k+1})^{H}\; ,
$$
and
$$
{1\over q_{k}}\log |z_{k}|^{-1}\le q_{k}^{-1}[c+H\log q_{k+1}+
(H-2)\log q_{k}]\; ,
$$
which is small when $\alpha$ is a Brjuno number and $k$ is large.
Lemma 5.20, Remark 5.18, and
Theorem 5.10 lead to the desired conclusion.
\itemitem{(I.2)} If $w\in\tilde{W}_{h}$ we will have
$$
|z-p_{k}/q_{k}|\ge \exp [-c(q_{k}q_{k+1})^h]\; ,
$$
thus
$$
|z_{k}|^{-1}\le cq_{k}^{-2}\exp [c(q_{k}q_{k+1})^h]\; ,
$$
and
$$
{1\over q_{k}}\log |z_{k}|^{-1}\le c q_{k+1}^h q_{k}^{h-1}\; ,
$$
which is small if $\alpha$ satisfies the diophantine condition we
have assumed and $k$ is large. Once again the conclusion follows
from Lemma 5.20, Remark 5.18,  and Theorem 5.10.
\item{(II)} Here we assume that $p/q$ is {\it not} one of the partial
fractions of $\alpha$. We denote $(p_{l}'/q_{l}')_{0\le l\le L}$
the partial fractions of $p/q$ and $k$ the largest integer such that
$p'_{k}/q'_{k}=p_{k}/q_{k}$. Clearly one has $k<L$ and
$p'_{L}/q'_{L}=p/q$. By a classical result ([HW], Theorem 184, p. 153)
$$
|\alpha - p/q| \ge {1\over 2q^{2}}\; .
$$
For $w\in\tilde{W}_{h}$ one has
$$
|z-p/q|\ge \exp[-cq^{2h}]\; ,
$$
thus
$$
|z_{L}|^{-1}\le cq^{-2}\exp [cq^{2h}]\; ,
$$
and
$$
{1\over q}\log |z_{L}|^{-1}\le q^{-1}cq^{2h}\; ,
$$
which is small since $h<1/2$. Taking into account Theorem 5.10,
Lemma 5.20,  and Remark 5.18
we only need to check that $B_{finite}(p_{k}/q_{k})$ and
$B_{finite}(p/q)$ are close.

Let us introduce
$$
\rho = \max_{k\le l<L}q_{l}^{\prime -1}\log{q'_{l+1}\over q'_{l}}(l-k+1)^{2}\;
{}.
\eqno(5.59)
$$
By Lemma 5.20 we have
$$
|B_{\hbox{fin}}(p_{k}/q_{k}) - B_{\hbox{fin}}(p/q)|\le
cq_{k}^{-1}+c\rho\; ,
$$
and we must show that $\rho$ is small. Let $l$ be such that
$$
\rho=q_{l}^{\prime -1}\log{q'_{l+1}\over q'_{l}}(l-k+1)^{2}\; .
\eqno(5.60)
$$
We have
$$
|p/q-p'_{l}/q'_{l}|\le {1\over q'_{l}q'_{l+1}}\; .
$$
On the other hand, if $w\in\tilde{W}_{h}$ and $z=w+\alpha\in V(p/q)$
one has
$$
|z-p/q|\le cq^{-2}
$$
from which it follows that
$$
|p/q-\alpha|\le c(\log q)^{-1/h}\; ,
$$
thus
$$
|\alpha - p'_{l}/q'_{l}|\le {1\over q'_{l}q'_{l+1}}+ c(\log q)^{-1/h}\, .
\eqno(5.61)
$$
By the choice (5.60) of $l$ we have
$$
q'_{l+1}=q'_{l}\exp [{\rho q'_{l}\over (l-k+1)^{2}}]\; ,
$$
which implies
$$
(\log q)^{-1/h} \le (\log q'_{l+1})^{-1/h}=
[{\rho q'_{l}\over (l-k+1)^{2}}+ \log q'_{l}]^{-1/h}\; .
$$
When $w$ approaches zero, $q$ must be large so $k$ must be large too and
if $\rho$ is not small
one has by (5.59)
$$
[{\rho q'_{l}\over (l-k+1)^{2}}+ \log q'_{l}]^{-1/h}\le 1/{10}
q_{l}^{\prime -2}
$$
since $h<1/2$.

We will also have
${1\over q'_{l}q'_{l+1}}\le {1\over 10} q_{l}^{\prime -2}$, thus
by (5.61)
$$
|\alpha - p'_{l}/q'_{l}|\le {1\over 5}q_{l}^{\prime -2}\; .
$$
But  (once again thanks to [HW], Theorem 184, p. 153)
this implies that  $p'_{l}/q'_{l}$ is one of the partial fractions of $\alpha$
and it must be $p_{k}/q_{k}$ by definition  of $k$. So one has $l=k$
and
$$
|\alpha - p_{k}/q_{k}|\le {1\over q_{k}q'_{k+1}}+c(\log
q'_{k+1})^{-1/h}\; ,
$$
with $q'_{k+1}=q_{k}\exp (\rho q_{k})$. This leads to the conclusion
that $\rho$ is
small if $\alpha$ verifies the diophantine condition of the second
part of the theorem. If $\alpha$ is a Brjuno number and $w\in W_{H}$
the condition
$$
|z-p/q|\le cq^{-2}
$$
implies the stronger inequality
$$
|\alpha - p/q|\le cq^{-2/H}
$$
thus
$$
|\alpha - p_{k}/q_{k}|\le {1\over q_{k}q'_{k+1}}+cq_{k+1}^{\prime -2/H}
$$
where, once again, $q'_{k+1}=q_{k}\exp (\rho q_{k})$. Therefore $\rho$
must be small in this case too.  \qed

\smallskip
\remark{5.21} A careful examination of the previous proof
leads to a slightly stronger version of the first part of Theorem
5.19 (inspired by the work of Risler, [Ri]).

The set $B$ of Brjuno numbers $\alpha$ has an injective image
into $l^{1}(\N )$ as follows:
$$
\alpha \mapsto (\beta_{l-1}\log \alpha_{l}^{-1})_{l\ge 0}\; .
$$
Let $K$ be a subset of $B$ such that its image is {\it  relatively
compact } in $l^{1}(\N )$. Then the convergence
$$
\lim_{w\in W_{H}\, , \,w\rightarrow 0}\IM {\cal B}(w+\alpha ) =
B(\alpha )
$$
is uniform w.r.t. $\alpha \in K$.

We recall that a subset $K$ of $l^{1}(\N )$ is relatively compact if
and only if
\item{(i)} $\forall n\ge 0$, $\exists C_{n}$ such that $\forall
(u_{l})_{l\ge 0}\in K$ one has $|u_{n}|\le C_{n}$;
\item{(ii)} $\forall \varepsilon >0$ $\exists n_{0}$ such that
$\forall
(u_{l})_{l\ge 0}\in K$ one has $\sum_{l>n_{0}}|u_{l}|\le
\varepsilon$.

\vskip .5 truecm
%%%%%%%%%%%%%%   Appendix 1 : real continued fraction %%%%%%%%%%
\beginsection A1. Appendix 1: real continued fractions\par
\vskip .3 truecm
\noindent
In this appendix we recall some elementary facts on standard real continued
fractions (we refer to [MMY], and references therein, for more general
continued fractions).

We will consider the iteration of the Gauss map
$$
A : (0,1) \mapsto [0,1] \; , \eqno(A1.1)
$$
defined by
$$
A (x) = {\ 1\ \over x} -
\left[\ {\ 1\ \over x}\ \right]\;. \eqno(A1.2)
$$
Let
$$
G = {\sqrt{5}+1\over 2}\; , \;
g = G^{-1} = {\sqrt{5}-1\over 2} \; . \;\;\;
$$
To each $x \in \Bbb R \setminus \Bbb Q$ we associate a continued fraction
expansion by iterating $A$ as follows. Let
$$
x_0  = x - [x] \; , \quad
            a_0  = [x] \; ,  \eqno(A1.3)
$$
then
$x = a_0 +  x_0$.
We now define inductively for all $n \ge 0$
$$
x_{n+1}  = A(x_n) \; , \quad
a_{n+1}  = \left[ {1 \over x_n} \right]  \ge 1\; , \eqno(A1.4)
$$
thus
$$
x_{n}^{-1} = a_{n+1} + x_{n+1} \; . \eqno(A1.5)
$$
Therefore we have
$$
x=a_0 + x_0=a_0+{1\over a_1 +
     x_1}= \ldots =a_0 + \displaystyle{1 \over a_1
     + \displaystyle{1  \over a_2 + \ddots +
     \displaystyle{1 \over a_n + x_n}}}\; ,
     \eqno(A1.6)
$$
and we will write
$$
x=[a_0,a_1,\ldots ,a_n,
     \ldots] \;. \eqno(A1.7)
$$
The nth-convergent is defined by
$$
{p_n \over q_n} = [a_0,a_1,\ldots ,
                     a_n] =
                     a_0 + \displaystyle{1 \over a_1
     + \displaystyle{1  \over a_2 + \ddots +
     \displaystyle{1 \over a_n }}}
     \;. \eqno(A1.8)
$$
The numerators $p_n$ and denominators
$q_n$ are recursively determined by
$$
p_{-1}=q_{-2}=1 \;\;,\;\;\;p_{-2}=q_{-1}=0 \;\;,\eqno(A1.9)
$$
and for all $n \ge 0$
$$
p_n = a_n p_{n-1} +  p_{n-2} \; , \quad
q_n = a_n q_{n-1} + q_{n-2} \; .  \eqno(A1.10)
$$
Moreover
$$
\eqalignno{x &= {p_n + p_{n-1}  x_n \over q_n + q_{n-1}
       x_n } \; , &(A1.11) \cr
              x_n &= - {q_n x -p_n \over q_{n-1} x - p_{n-1}}
              \; , &(A1.12) \cr
              q_n &p_{n-1} - p_n q_{n-1} = (-1)^n \;\; .  &(A1.13) \cr}
$$
Let
$$
\beta_n = \Pi_{i=0}^n x_i = (-1)^n (q_n x - p_n)\quad\hbox{for\ }n\ge 0,\quad
  \hbox{and\ }\be_{-1}=1\;\; .\eqno(A1.14)
$$
{}From the definitions given one easily proves by induction
the following proposition (see [MMY])\par
\smallskip
\Proc{Proposition A1.1} {For all
$x \in \Bbb R \setminus \Bbb Q$ and for all $n \ge 1$ one has
\item{(i)}\qquad $ \left|q_n x - p_n\right|
={\dst 1\over\dst q_{n+1}+q_nx_{n+1}}$,
so that ${\dst 1\over\dst 2}<\beta_nq_{n+1}<1$\ ;
\item{(ii)}$\beta_n\le g^n$ and $q_n\ge{\dst1
\over\dst 2}G^{n-1}$\ .
}
\par
\smallskip
\noindent
Note that from {\it (ii)} it follows that
$\sum_{k=0}^\infty {\log q_{k}\over q_{k}}$ and
$\sum_{k=0}^\infty {1\over q_{k}}$
are always convergent and their sum is uniformly bounded.

With the notations of Section 2.1, equation (A1.5) can be written
$x_{n}=g(a_{n+1})x_{n+1}$, thus we have
$x_{0}=g({a_1})g({a_2})\cdots g({a_n})x_{n}$. The following
characterization of the monoid ${\cal M}$  defined in Section 2.1
is therefore relevant.

\smallskip
\Proc{Proposition A1.2}{Let
$g(m)=\left(\matrix{0 & 1\cr 1 & m\cr}\right)$, where $m\ge 1$.
${\cal M}$ is the {\it free monoid}, with unit,
generated by the elements $g(m)$, $m\ge 1$: each element $g$ of
${\cal M}$ can be written as
$$
g = g(m_1)
\cdots g(m_r)\; , \;\;r\ge 0\; , \; m_{i}\ge 1
$$
and this decomposition is unique.}
\smallskip
\proof
Let ${\cal M}^*$ be the monoid with unit generated by the $g(m), m\ge1$.
If $m\ge 1$ one has $g(m)\in {\cal M}$
and ${\cal M} g(m)\subset {\cal M}$, thus ${\cal M}^*\subset {\cal M}$.
\par
Conversely let $g\in {\cal M}$, $g\neq {\rm id}$. We now prove that
there exists a unique integer $m\ge 1$ such that
$
\left(\matrix{a' & b'\cr c' & d'\cr}\right) =g' =g(g(m))^{-1}\in
{\cal M} \; ,
$
which leads to the conditions
$
b'=a\; , \;\; d'=c\; , \;\; a'=b-ma\; , \;\; c'=d-mc\; .
$
We consider separately three cases.
\par
\item{1)} $a=0$, thus $b=c=1$ and $g=g(d)$ where $d\ge 1$. If there were
$m\ge 1$ such that $g'=g(g(m))^{-1}\in {\cal M}$, and $g'\neq {\rm id}$,
one should have $b'=0$, thus
$a'=0$ and $a'd'-b'c'=0$ which is impossible.
\item{2)}
$a=1$, thus $b,c\ge 1$ and $d=bc\pm 1$. We therefore
have $b'=1$ and
$a'=0$ or $a'=1$. If $a'=0$ then $b=m$ and $c'=d-bc$,
which is admissible if and only if
$d=bc+1$. If $a'=1$ then $b=m+1$ from which it follows that
$c'=d-mc=bc\pm 1-mc=c\pm 1$, which is admissible if and only if
$d=bc-1$, and then
$b,c\ge 2$.
\item{3)}
$a>1$. Since $a' \wedge b' =1$ the relation $0\le a'=b-ma\le b'=a$
determines uniquely $m\ge 1$ and one has $0<a' <b'$. But one also has
$b'=a\le c=d'$ and $|a'd'-b'c'|=1$, from which one easily gets
$d'\ge c'\ge a'$. \qed
%\vskip .5 truecm
%%%%%%%%%%%%%%%%% Appendix 2: Hyperfunction %%%%%%%%%%%%%%%%
\beginsection A2. Appendix 2: Hyperfunctions \par
%\vskip .3 truecm
\noindent
{\bf A2.1} We follow here [H], Chapter 9.
Let $K$ be a non empty compact subset of ${\Bbb R}$.
A {\it hyperfunction with support in} $K$ is a linear functional $u$
on the space ${\cal O}(K)$ of functions analytic in a neighborhood of $K$
such that for all neighborhood $V$ of $K$ there is a constant
$C_{V}>0$ such that
$$
|u(\varphi )| \le C_{V}\sup_{V}|\varphi|\; , \;\;\;
\forall \varphi\in {\cal O}(V)
\; .
$$
We denote by $A'(K)$ the space of hyperfunctions with support in $K$.
It is a Fr\'echet space: a seminorm is associated to each
neighborhood $V$ of $K$.  One has the following
proposition.
\vfill\eject
\smallskip
\Proc{Proposition A1.1}{ The spaces $A'(K)$ and
${\cal O}^{1}(\overline{\Bbb C}\setminus K)$ are canonically isomorphic.
To each $u\in A'(K)$ corresponds $\varphi\in
{\cal O}^{1}(\overline{\Bbb C}\setminus K)$ given by
$$
\varphi (z) = u(c_{z})\; , \; \forall z\in {\Bbb C}\setminus K\; ,
$$
where $c_{z}(x) = {1\over \pi}{1\over x-z}$. Conversely to each
$\varphi \in {\cal O}^{1}(\overline{\Bbb C}\setminus K)$
corresponds the hyperfunction
$$
u(\psi ) = {i\over 2\pi} \int_{\gamma} \varphi (z)\psi (z) dz\; , \;
\forall \psi \in A
$$
where $\gamma$ is any piecewise ${\cal C}^{1}$ path  winding around
$K$ in the positive direction. We will also use the notation
$
u(x) = {\dst 1\over\dst 2i}[\varphi (x+i0)-\varphi (x-i0)]
$
for short.}
\smallskip

\smallskip
\noindent
{\bf A2.2} Let ${\Bbb T}^{1}= {\Bbb R}/{\Bbb Z}\subset
\C/\Z$. A {\it hyperfunction on} $\T$ is a linear funtional $U$ on
the space ${\cal O}({\Bbb T}^{1})$ of
functions analytic in a complex
neighborhood of ${\Bbb T}^{1}$ such that for all neighborhood $V$
of $\T$ there exists $C_{V}>0$ such that
$$
|U(\Phi )| \le C_{V}\sup_{V}|\Phi|\; , \;\;\;
\forall \Phi\in {\cal O}(V)
\; .
$$
We will denote
$A'({\Bbb T}^{1})$.
the Fr\'echet space of  hyperfunctions with support in $\T$.
For $U\in A'(\T )$, let $\hat{U}(n) := U(e_{-n})$ with
$e_{n}(z)=e^{2\pi i n z}$. The doubly infinite sequence
$(\hat{U}(n))_{n\in {\Bbb Z}}$ satisfies
$$
|\hat{U}(n)|<C_{\varepsilon}e^{2\pi |n|\varepsilon}\; .
$$
for all $\varepsilon >0$ and for all $n\in \Z$ with a suitably chosen
$C_{\varepsilon}>0$. Conversely any such sequence is the Fourier
expansion of a unique hyperfunction with support in $\T$.

Let ${\cal O}_{\Sigma}$ denote the complex vector space of
holomorphic functions $\Phi\, : {\Bbb C}\setminus {\Bbb R}\rightarrow
{\Bbb C}$, $1$--periodic, bounded at $\pm i \infty$ and such that
$\Phi (\pm i\infty) := \lim_{\IM z\rightarrow \pm\infty }\Phi (z)$
exist and verify $\Phi (+i\infty) = -\Phi (-i\infty )$.

\smallskip
\Proc{Proposition A2.2}{ The spaces $A'({\Bbb T}^{1})$ and
${\cal O}_{\Sigma}$ are canonically isomorphic.
To each $U\in A'({\Bbb T}^{1})$ corresponds $\Phi\in
{\cal O}_{\Sigma}$ given by
$$
\Phi (z) = U(C_{z})\; , \; \forall z\in {\Bbb C}\setminus K\; ,
$$
where $C_{z}(x) = \cotg \pi ( x-z)$. Conversely to each
$\Phi\in {\cal O}_{\Sigma}$
corresponds the hyperfunction
$$
U(\Psi ) = {i\over 2} \int_{\Gamma} \Phi (z)\Psi (z) dz\; , \;
\forall \Psi \in A({\Bbb T}^{1})
$$
where $\Gamma$ is any piecewise ${\cal C}^{1}$ path  winding around a
closed interval $I\subset {\Bbb R}$ of length $1$
in the positive direction. We will also use the notation
$$
U(x) = {1\over 2i}[\Phi (x+i0)-\Phi (x-i0)]
$$
for short.}
\smallskip

The nice fact is that the following diagram commutes:
$$
\diagram{
A'([0,1]) &
\hfl{}{} & {\cal O}^{1}(\overline{\Bbb C}\setminus [0,1])
    \cr
\vfl{\sum_{Z}}{}&&\vfl{}{\sum_{Z}}
    \cr
A'({\Bbb T}^{1}) &
\hfl{}{} & {\cal O}_{\Sigma}
    \cr}
$$
the horizontal lines are the above mentioned isomorphisms, $\sum_{Z}$
is defined in 2.2.2.
%\vskip .5 truecm
%%%%%%%%%%%%% Appendix  3: Some properties of the dilogarithm %%%%%%%%
\beginsection A3. Appendix 3: Some properties of the dilogarithm  \par
%\vskip .3 truecm
\noindent
{\bf A3.1}
The classical dilogarithmic series
(see [Le], [O] for more information) is defined by
$$
\hbox{Li}_2(z)=\sum_{n=1}^{+\infty}{z^n\over n^2} \eqno(A3.1)
$$
and it is convergent for $|z|\le 1$. Since $-\log (1-z)=
\sum_{n=1}^{+\infty} {z^n\over n}$, dividing by $z$ and integrating
one obtains the analytic continuation of the dilogarithm to
${\Bbb C}\setminus [1,+\infty )$ by means of the integral
formula
$$
\hbox{Li}_2(z)=-\int_0^z{\log (1-t)\over t}dt
=\int_0^z\left(\int_0^t{d\zeta\over 1-\zeta}\right){dt\over t}
\eqno(A3.2)
$$
which we will use as a definition of the dilogarithm. Note that
$[1,+\infty )$ is a branch cut.
\par
Since
$$
\hbox{Li}_2 (z) = z \int_0^1{\log t\over tz-1}dt
$$
one obviously has that
$$
\hbox{Li}_2\left({1\over z}\right) = -\int_0^1 {\log t\over z-t}dt\; ,
\eqno(A3.3)
$$
which shows that $\hbox{Li}_2\left({1\over z}\right)$ is the
Cauchy--Hilbert transform of
the real function
$$
\varphi_0(t) = \cases{-\log t & if $t\in [0,1]$\cr
                      0 & elsewhere\cr}
					  \eqno(A3.4)
$$
Note also that
$$
\hbox{Im} \hbox{Li}_2(t\pm i0) = \pm \pi \log t\; ,  \eqno(A3.5)
$$
where $t\in [1,+\infty)$. Moreover
$$
|\hbox{Li}_2(z)|={\cal O}(\log^2|z|) \;\hbox{as} |z|\to +\infty\; .
\eqno(A3.6)
$$

\smallskip
\noindent
{\bf A3.2}\ {\it Euler's functional equations}
$$
\eqalignno{
\hbox{Li}_2(z)+\hbox{Li}_2\left({1\over z}\right) &= -{1\over 2}(\log (-z))^2
-{\pi^2\over 6} \; , &(A3.7)
\cr
\hbox{Li}_2(z)+\hbox{Li}_2(1-z) &= -\log z\log (1-z)+{\pi^2\over 6}
\; , &(A3.8)
\cr}
$$
where z varies in ${\Bbb C}\setminus [0,+\infty]$ and ${\Bbb C}\setminus
((-\infty , 0]\cup [1,+\infty))$ respectively.

\smallskip
\noindent
{\bf A3.3}\ {\it Special values}
$$
\eqalign{
\hbox{Li}_2 (1) &= {\pi^2\over 6}\; , \;\;\hbox{Li}_2(-1)=-{\pi^2\over 12}
\; , \cr
2\hbox{Li}_2(1/2)&={\pi^2\over 6}-(\log 2)^2\; , \;\;
\hbox{Li}_2(2\pm i0)={\pi^2\over 4}\pm \pi i\log 2\; . \cr
}
$$
%\vskip .5 truecm
%%%%%%%%%%%%%%%%% Appendix 4: Even Brjuno functions %%%%%%%%%%
\beginsection A4. Appendix 4: Even Brjuno functions \par
%\vskip .3 truecm
In [MMY] we also considered an even version of the Brjuno function
and we proved that
this differs from the one considered here by a $1/2$--H\"older
continuous function. In this appendix we explicit the relation
among the two associated complex Brjuno functons.

\smallskip
\noindent
{\bf A4.1}
Let $\sigma$ denote the matrix $\left(\matrix{-1 & 1\cr 0 &
1\cr}\right)$ which corresponds to $x\mapsto 1-x$.

At the {\it real level}, replacing periodic even functions with functions
on $[0,1/2]$ and null outside this interval,
the operator $T_{even}$ acting on $L^{2}([0,1/2])$ (for example)
can be written explicitely as
$$
T_{even}f(x) = \sum_{m\ge 2}xf\left({1\over x}-m\right)+\sum_{m\ge 3}
xf\left(m-{1\over x}\right)\; .
$$

At the {\it complex level} (i.e. after the identification of ${\cal
A}'([0,1/2])$ to ${\cal O}^{1}(\overline{\C}\setminus [0,1/2])$)
one gets
$$
T_{even}\varphi = \sum_{m\ge 2}L_{g(m)}\varphi + \sum_{m\ge 2}
L_{g'(m)}\varphi
$$
where
$
g'(m)= \left(\matrix{0 & 1\cr -1 &
m+1\cr}\right) = g(m)\sigma\; .
$

\smallskip
\noindent
{\bf A4.2} We want to consider
$$
(1-T_{even})^{-1}\, : \, {\cal O}^{1}(\overline{\C}\setminus [0,1/2])
\rightarrow {\cal O}^{1}(\overline{\C}\setminus [0,1/2])
$$
and then one will have to make the resulting function even and
periodic, thus one will take
$$
\sum_{Z}(1+L_{\sigma})(1-T_{even})^{-1}\, :
{\cal O}^{1}(\overline{\C}\setminus [0,1/2])
\rightarrow {\cal O}_{even}(\H/\Z)\; .
$$

Note that $Z\sqcup Z\sigma = Z\sqcup \sigma Z$.

When we expand $(1+L_{\sigma})(1-T_{even})^{-1}$ we obtain a sum
$\sum L_{g}$ where the matrices $g$ have the form
$$
g = \varepsilon_{0}g(i_{1})\varepsilon_{1}g(i_{2})\ldots g(i_{r})
\varepsilon_{r}
$$
with $r\ge 0$, $i_{k}\ge 2$ and $\varepsilon_{k}\in \{1,\sigma\}$.

Note that $\sigma g(i)=g(1)g(i-1)$ for all $i\ge 2$, thus all
matrices of the form $\varepsilon_{0}g(i_{1})\varepsilon_{1}\ldots
g(i_{r})$ belong to the monoid ${\cal M}$.

\smallskip
\noindent
{\bf A4.3} Let $r\ge 0$,
$$
{\cal M}^{(r)}=\{ g(j_{1})\ldots g(j_{r})\; , \; j_{k}\ge 1\}
\;
$$
and let $\widehat{\cal M}^{(r)}$ denote the part of ${\cal M}^{(r)}$ made of
matrices which can be written as a product
$\varepsilon_{0}g(i_{1})\varepsilon_{1}\ldots
g(i_{s})$. We have the following

\Proc{Lemma A4.1}{Each matrix $g\in \widehat{\cal M}^{(r)}$ can be uniquely
written as  a product $\varepsilon_{0}g(i_{1})\varepsilon_{1}\ldots
g(i_{s})$. Moreover one has ${\cal M}^{(0)}=\widehat{\cal M}^{(0)}=\{1\}$ and
for all $r>0$ }
$$
{\cal M}^{(r)}\setminus \widehat{\cal M}^{(r)} = \widehat{\cal M}^{(r-1)}
g(1)\; .
$$

\smallskip
\proof
Uniqueness is evident (just consider the first place at which
the product $\varepsilon_{0}g(i_{1})\varepsilon_{1}\ldots
g(i_{s})$ differs from $\varepsilon'_{0}g(i'_{1})\varepsilon'_{1}\ldots
g(i'_{s'})$). The second assertion follows easily from the remark that
$\widehat{\cal M}^{(r)}$ is indeed made of matrices $g=g(j_{1})\ldots
g(j_{r})$ which end with an even number of $g(1)$'s.
\qed

\smallskip
Let $\widehat{\cal M}= \sqcup_{r\ge 0} \widehat{\cal M}^{(r)}$ and
$\sum_{\widehat{\cal M}}= \sum_{g\in \widehat{\cal M}}L_{g}$. One clearly
has
$$
(1+L_{\sigma})\circ (1-T_{even})^{-1} =\bigg(\sum_{\widehat{\cal M}}\bigg)
\circ (1+L_{\sigma})
$$
and by the previous Lemma
$$
\eqalign{
\sum_{\cal M}=
\sum_{g\in{\cal M}} L_{g}&= \sum_{r\ge 0}\sum_{g\in{\cal M}^{(r)}}L_{g} =
\sum_{r\ge 0}\bigg(\sum_{g\in\widehat{\cal M}^{(r)}}L_{g}+
\sum_{g\in{\cal M}^{(r)}\setminus
\widehat{\cal M}^{(r)}}L_{g}\bigg)\cr
&= \sum_{r\ge 0}\sum_{g\in\widehat{\cal M}^{(r)}}L_{g} + \sum_{r\ge 1}
\sum_{g\in\widehat{\cal M}^{(r-1)}g(1)}L_{g}=
\bigg(\sum_{\widehat{\cal M}} \bigg)\circ (1+L_{g(1)})\; . \cr}
$$
We are therefore led to conclude that
$$
(1+L_{\sigma})\circ(1-T_{even})^{-1}=\bigg(\sum_{\cal M}\bigg)\circ
(1+L_{g(1)})^{-1}\circ (1+L_{\sigma})\; .
$$

\smallskip
\noindent
{\bf A4.4} It is not hard to check, as we did for the
monoid ${\cal M}$ in Proposition A1.2, that a matrix $g$ belongs to
the monoid $\widehat{\cal M}$
if and only if $d\ge 2b>0$, $c\ge 2a\ge 0$,
$d\ge Gc$, where $G={\sqrt{5}+1\over 2}$. Moreover the
decomposition $g=\varepsilon_{0}g(i_{1})\varepsilon_{1}\ldots
g(i_{r})$ is unique.
%\vskip .5 truecm
%%%%%%%%%%%% Appendix 5: The real Brjuno function as a cocycle %%%%%%%
\beginsection  A5. Appendix 5: The real Brjuno function as a
cocycle\par
%\vskip .3 truecm
In this Appendix we show how to interpret the real Brjuno function
as a cocycle under the action of $\hbox{PGL}(2,{\Bbb Z})$ on
${\Bbb R}\setminus {\Bbb Q}$. To this purpose we first
recall some basic definitions taken from the cohomology of
groups. We refer to [Ja] and [Se] for more information and the proofs.

\smallskip
\noindent
{\bf A5.1 Group cohomology}\par
Let $G$ be a group and $M$ an abelian group with a left $G$--action,
i.e. a structure of a left ${\Bbb Z}^{G}$--module. Recall that for
$n\ge 0$, one defines
\item{(i)} the abelian group of $n$--cochains $C^n(G,M)$, whose
elements are applications from $G^n$ to $M$.\par
\item{(ii)} the coboundary operator
$d^n\, :C^n(G,M)\rightarrow C^{n+1}(G,M)$:
$$
\eqalign{ (d^n f)(g_0,\ldots ,g_n) &= g_0f(g_1,\ldots ,g_n)
+ \sum_{i=0}^{n-1}
(-1)^{i+1} f(g_0,\ldots ,g_ig_{i+1},\ldots ,g_n) \cr
& + (-1)^{n+1}
f(g_0,\ldots ,g_{n-1})\;  ; \cr
}
$$ \par
\item{(iii)} the abelian subgroups of $n$--cocycles
$Z^n(G,M)= {\rm Ker}\ d^n$
and of $n$--coboun\-daries, $B^n(G,M)={\rm Im}\ d^{n-1}$.\par
\item{(iv)} the $n$-th cohomology group
$H^n=Z^n(G,M)/B^n(G,M)$.\par\smallskip
Identifying $C^0(G,M)$ with $M$, one has
$H^0(G,M) =Z^0(G,M)=\{m\in M\ ;\ gm=m\ \hbox{for all}\  g\in G\}$ with $M$.
An application $c\ : G\to M$ is a $1$--cocycle iff
$
c(g_0g_1)=g_0c(g_1)+c(g_0)\ ,
$
and a $1$--coboundary iff
$
c(g)=g\cdot m -m\ \hbox{for some $m$ and all $g\in G$}.
$
\smallskip
\noindent
{\bf A5.2 Automorphic factors, cocycles and coboundaries.}\par
Let $G$ be a group acting on the left on a set $X$. Let $A$ be an abelian
ring, $A^*$ the multiplicative group of invertible elements of $A$,
and $M$ a $A$--module.
A function $\chi\, : G\times X\rightarrow A$ is an {\it automorphic factor}
if the application
$ G\times M^X  \longrightarrow M^X$ given by
$$
(g,\varphi ) \ \longmapsto g\cdot \varphi\, : \; g\cdot \varphi (x)
= \chi (g^{-1},x)\varphi (g^{-1}\cdot x)\; \forall x\in X\; ,
$$
defines a left action of $G$ on $M^X$: one must have
$$
\chi (g_0g_1, x)= \chi (g_0, g_1x)\chi (g_1, x)\; .
$$
One has therefore given to $M^X$ the structure of a
${\Bbb Z}^{[G]}$--module.
The  coboundary of an element $\varphi\in M^X$ is given  by
$ d^0\varphi (g) = g\cdot\varphi - \varphi \; , $
i.e.
$$
d^0\varphi (g)(x) = \chi (g^{-1},x)\varphi (g^{-1}\cdot x)
- \varphi (x)\;\forall x\in X\; .
$$
A $1$--{\it cocycle} is an application $c\, : G\rightarrow M^X$ verifying
$ g_0\cdot c(g_1) - c(g_0g_1)+c(g_0)=0\; , $
i.e., letting $\check c (g) = c(g^{-1})$:
$$
\check c (g_0g_1) = c(g_1^{-1}g_0^{-1}) =
c(g_1^{-1})+g_1^{-1}c(g_0^{-1})=\check c (g_1)+g_1^{-1}\check c
(g_0)\; ,
$$
or, equivalently,
$
\check c(g_0g_1,x) = \chi (g_1,x)\check c(g_0,g_1\cdot x)+
\check c (g_1,x)\; \forall x\in X\; .
$

\smallskip
\noindent
{\bf A5.3 Action of $\hbox{PGL}(2,{\Bbb Z})$ on
${\Bbb R}\setminus {\Bbb Q}$.}\par
Let us consider $G=\hbox{PGL}(2,{\Bbb Z})$ and
$X={\Bbb R}\setminus {\Bbb Q}$, the action being
given by the homographies.
The transformations $T(x)=x+1$ and $S(x)=x^{-1}$
generate $\hbox{PGL}(2,{\Bbb Z})$. One has the following
more precise result:

\smallskip
\noindent
{\bf Proposition A5.1}\ {\it Let $g\in \hbox{PGL}(2,{\Bbb Z})$ and let
$x_0\in {\Bbb R}\setminus {\Bbb Q}$. There exist $r\ge 0$
and elements $g_1,\ldots ,g_r\in \{S,T,T^{-1}\}$ such
that
\item{(i)}Ú$g = g_r\ldots g_1$;\par
\item{(ii)}Ú let $x_i=g_ix_{i-1}$ for $1\le i\le r$,
then $x_{i-1}>0$  if $g_i=S$.\par

\noindent
Moreover one can require that $g_ig_{i-1}\not= 1$ for $0<i\le r$,
and in this case $r,g_1,\ldots , g_r$ are uniquely determined. }

\smallskip
\proof

First we prove the existence. Let $U(x)=-x$. We consider five cases:
\item{1.} if $g=T^{\pm 1}$, any $x_0$, one takes $r=1$, $g_1= T^{\pm 1}$.
\item{2.}Úif $g=U$ and $x_0\in (0,1)$ then $r=6$ and
$g_1=S$, $x_1=x_0^{-1}$; $g_2=T^{-1}$, $x_2={1-x_0\over x_0}$;
$g_3=S$, $x_3={x_0\over 1-x_0}$; $g_4=T$, $x_4={1\over 1-x_0}$;
$g_5=S$, $x_5=1-x_0$; $g_6=T^{-1}$, $x_6=-x_0$.
\item{3.}Úif $g=U$ and $x_0\in (n,n+1)$, $n\in {\Bbb Z}$, one is led
to consider the previous case by using $U=T^{-n}UT^{-n}$.
\item{4.} if $g=S$ it is immediate if $x_0>0$, and if $x_0<0$ one is
led to consider case 3. by using the relation $S=USU$.
\item{5.}ÚOne has the cases $g=S$ and $g=T^{\pm 1}$ for all $x_0$.
Since $S$ and $T$ generate $\hbox{PGL}(2,{\Bbb Z})$ this implies
the existence in all possible cases.

We can now prove uniqueness. It is sufficient to show that if
$r>0$ and $g_1,\ldots ,g_r\in \{S,T,T^{-1}\}$, $x_0\in {\Bbb R}\setminus
{\Bbb Q}$ verify
$$
g_r\ldots g_1 = 1\; ,\quad \hbox{and}\ \
x_{i-1} >0 \ \hbox{if}\  g_i=S \; (1\le i \le r)\ ,
$$
then there exists $1< i\le r$ such that $g_ig_{i-1}=1$.
We prove this by contradiction: let $r$ be minimal, $r>0$.

If $x_0<0$ one must have $r\ge 2$, $g_1=T$ and $g_r=T^{-1}$, thus
$r\ge 3$ and $g_{r-1}\ldots g_2=1$, which contradicts the minimality
of $r$.

One is led to assume $x_0>0$. Clearly one must have $x_i>0$ for all
$i\in [0,r]$. Let $i_1<i_2<\ldots <i_k$ denote the indices $i$ such that
$g_i=S$. The integer $k>0$, even (because of the determinant sign) and one
has $i_{l+1}\ldots i_1\ge 2$ for $1\le l<k$. Let us assume that
$x_{i_1-1}>1$. Then $x_{i_1}\in (0,1)$ thus $x_{i_2-1}=x_{i_1}+(i_2-i_1
-1)>1$. Therefore $x_{i_l-1}>1$ for all $1\le l \le k$. But then
$\prod_1^r{dg_i\over dx_{i-1}}(x_{i-1})<1$, in contradiction with the
assumption
$g_r\ldots
g_1=1$. If there exists $l$ such that $x_{i_l-1}<1$ one permutes
circularly all $g_i$ and $x_i$ ($\hbox{mod}\, r$)
until one is back to the  case previously considered. Finally
if $x_{i_l-1}<1$ for all $1\le l\le k$ then
$\prod_1^r{dg_i\over dx_{i-1}}(x_{i-1})>1$, which is again
in contradiction with the assumption $g_r\ldots
g_1=1$.
\qed
\smallskip

\smallskip
\Proc{Corollary A5.2} {Let $A$ be an abelian ring, and the maps $T$
and $S$ such that $t\, :{\Bbb R}\setminus
{\Bbb Q}\rightarrow A^*$, $s\,:(0,1)\cap ({\Bbb R}\setminus
{\Bbb Q})\rightarrow A^*$. There exists a unique automorphic factor
$\chi$ such that}
$$
\eqalign{
\chi (T,x) &= t(x) \; \hbox{for all}\,x\in {\Bbb R}\setminus
{\Bbb Q}= X \; , \cr
\chi (S,x) &= s(x)\; \hbox{for all}\, x\in X\cap (0,1)\; . \cr
}
$$

\smallskip
\proof
Let $s(x)=(s(x^{-1}))^{-1}$ for all $x\in X\cap (1,+\infty )$. The
map $s$ is therefore defined on $X\cap (0,+\infty )$ and one must have
$$
\chi (S,x)=s(x) \; \hbox{for all}\, x\in X, x>0\; .
$$
{a)}  The {\it uniqueness} of $\chi$ follows from the {\it existence}
in the previous proposition: if $g\in \hbox{PGL}(2,{\Bbb Z})$ and
$x_0\in X$ one must have
$$
\chi (g,x_0) = \prod_{i=1}^r\chi (g_i,x_{i-1}) \eqno(1)
$$
where $g_1,\ldots ,g_r$ and $x_1,\ldots ,x_r$ are defined in the
proposition and
$$
\eqalign{
\chi (T,x) &= t(x) \;\hbox{for all}\, x\in X\cr
\chi (T^{-1},x) &= (t(x-1))^{-1} \;\hbox{for all}\, x\in X\cr
\chi (S,x) &= s(x) \;\hbox{for all}\, x\in X, x>0.\cr
}
$$
{b)} The {\it existence} of $\chi$ follows from the {\it uniqueness}
in the previous proposition: here we use (1) with $r$ minimal
(i.e. $g_ig_{i-1}\not= 1$ for all $1<i\le r$) to define $\chi$ and one
must check that
$$
\chi (g'g,x_0) = \chi (g',gx_0)\chi (g,x_0)\; .
$$
Let
$
y_0=x_r=gx_0\; , \;\;\; g'=g_s'\ldots g_1'\; ,
$
following the previous proposition, and let
$$
g''=g'g\; , \;\;\; g_i'' = \cases{ g_i & if $1\le i\le r$\cr
g_{i-s}' & if $r<i\le r+s$\cr}\; \;\; x_i=y_{i-r}
$$
Then
$
g'' = g_{r+s}''\ldots g_1''
$
satisfies the conclusions of the proposition. The decomposition may
be not minimal (since one may have $g_1'g_r=1$) but one can obtain a
minimal decomposition by deleting $g_1'g_r$ if $g_1'g_r=1$ then
(if $g_1'g_r=1$) by deleting $g_2'g_{r-1}$ if $g_2'g_{r-1}=1$ and so
on. Given the definition (1) of $\chi$ the automorphic property
is now verified if
$$
\chi (g,g^{-1}x)\chi (g^{-1},x) = 1
$$
when $g=T,T^{-1}$ or $S$ and $x>0$ if $g=S$, which is immediate to check.
\qed
\smallskip

\smallskip
\Proc{Corollary A5.3}\ {Let $A$ be an abelian ring, $\chi$ an automorphic
factor, $M$ a $A$--module, $M^X$ with the structure of
${\Bbb Z}^{[G]}$--module defined by $\chi$. Let
$$
\eqalign{
\check c_T \, &:X\rightarrow M\cr
\check c_S \, &:X\cap (0,1)\rightarrow M\cr
}
$$
denote two maps. There exists a unique cocycle $\check c\,:G\times M
\rightarrow M$ such that }
$$
\eqalign{
\check c (T;x) &= \check c_T (x)\;\hbox{for all}\,x\in X\cr
\check c (S;x) &= \check c_S (x)\;\hbox{for all}\,x\in X\cap (0,1)\; .\cr
}
$$

\smallskip
\proof
One must have
$$
\check c (T^{-1};x) = -\chi (T^{-1},x)\check c_T(x-1)\; \hbox{for all }
x\in X
$$
and
$$
\check c (S;x) = -\chi (S,x)\check c_S(x^{-1})\; \hbox{for all }
x\in X\; , \; x>1\; .
$$
Moreover, if $g=g_r\ldots g_1$ and $x_0$ are given as in the proposition
$$
\check c (g;x_0) = \sum_{i=1}^r(\check c (g_i,x_{i-1})
\chi (g_{i-1}\ldots g_1,x_0)) \eqno(1')
$$
from which the uniqueness follows. The proof of existence is the
same as the one given for Corollary A5.2. \qed

\smallskip
\noindent
{\bf A5.4 The real Brjuno function as a cocycle.}\par
Let $A={\Bbb R}$, $t(x)=1$ and $s(x)=\varepsilon x^\nu$
with $\varepsilon\in\{-1,+1\}$, $\nu\in {\Bbb R}$ and apply Corollary
A5.2. Then
$$
\eqalign{
\chi (T^n,x) &= 1 \; , \;\;\hbox{for all}\,n\in {\Bbb Z}\;,\; x\in X\cr
\chi (S,x) &= \varepsilon x^\nu\; , \;\;
\hbox{for all}\, x\in X\; ,\; x>0\; .\cr
}
$$
If $x_0\in (0,1)$, one has seen that $U=T^{-1}STST^{-1}S$, thus
$$
\chi (U,x_0) = \varepsilon x_0^\nu\varepsilon\left({1-x_0\over x_0}
\right)^\nu\varepsilon\left({1\over 1-x_0}\right)^\nu =
\varepsilon
$$
{}From $U=T^nUT^n$ it follows that
$
\chi (U,x)=\varepsilon
$
for all $x\in X$, and from $S=USU$ follows that
$
\chi (S,x) = \varepsilon\varepsilon |x|^\nu\varepsilon
$
for $x<0$, or
$
\chi (S,x)=\varepsilon |x|^\nu
$
for all $x\in X$. One concludes that one must have
$$
\chi (g,x) = \cases{ |cx+d|^\nu & if $\varepsilon = +1$\cr
\det (g)|cx+d|^\nu & if $\varepsilon =-1$\cr}
$$
for all $g=\left(\matrix{ a & b \cr c & d\cr}\right)\in\hbox{PGL}(2,
{\Bbb Z})$.\par
Consider now the functional equations
$$
\eqalign{
B_f(x)&=xB_f(1/x)+f(x)
\quad, \quad x\in (0,1)\cap{\Bbb R}\setminus{\Bbb Q}\cr
B_f(x)&=B_f(x+1)\quad, \quad x\in{\Bbb R}\setminus{\Bbb Q}\ ,\cr}
$$
where $f\ :\
(0,1)\cap{\Bbb R}\setminus{\Bbb Q}\to {\Bbb C}$ is given.
Now we look for $B_f\ :\ {\Bbb R}\setminus {\Bbb Q}\to {\Bbb C}$,
and we easily see that the relevant automorphic factor is the case
$\epsilon=+1$, $\nu=+1$ above (other values of $\nu$ have also
been considered in [MMY]). By Corollary A5.3, there exists exactly one
$1$--cocycle $c_f$ such that
$$
\eqalign{
C_f(T,x)=&0\quad \forall x\in {\Bbb R}\setminus {\Bbb Q}\cr
C_f(S,x)=&f(x)\quad \forall x\in {\Bbb R}\setminus {\Bbb Q}
\cap (0,1)\ .\cr}
$$
The $1$--cocycle is a $1$--coboundary if and only if the functional equations
have a solution $B_f$, in which case we have $c_f=-d^0(B_f)$.
These considerations also apply and may become fruitful in  case we restrict
${\Bbb C}^{{\Bbb R}\setminus{\Bbb Q}}$ to one of its ${\Bbb
C}^{[G]}$--submodules,
for instance measurable functions.

%\vskip .5 truecm

%%%%% references %%%%
\beginsection References \par
%\vskip .3 truecm
\item{[BG1]} A. Berretti and G. Gentile, ``Scaling properties for the radius
of convergence of the Lindstedt series: the standard map'', University of
Roma (Italy), Preprint
(1998), to appear in Journal de Math\'ematiques pures et appliqu\'ees.
\item{[BG2]} A. Berretti and G. Gentile, ``Bryuno function and the standard
map'', University of Roma (Italy), Preprint (1998).
\item{[BPV]} N. Buric, I. Percival and F. Vivaldi
``Critical Function and Modular Smoothing'' {\it Nonlinearity}
{\bf 3} (1990), 21--37.
\item{[Br]} A. D. Brjuno ``Analytical form of differential equations''
{\it Trans. Moscow Math. Soc.} {\bf 25} (1971), 131--288; {\bf 26}
(1972), 199--239.
\item{[Da1]} A.M. Davie ``The critical function for the
semistandard map'' {\it Nonlinearity} {\bf 7} (1994), 219--229.
\item{[Da2]} A. M. Davie ``Renormalisation for analytic
area--preserving maps'' University of Edinburgh preprint, (1995).
\item{[Du]} P. L. Duren ``Theory of $H^p$ spaces'' Academic
Press, New York, (1970).
\item{[Ga]} J. B. Garnett ``Bounded Analytic Functions'' Academic Press,
New York, (1981).
\item{[GCRF]} J. Garcia--Cuerva and J.L. Rubio de Francia ``Weighted Norm
Inequalities and Related Topics'' North Holland Mathematical
Studies {\bf 116}, Amsterdam, (1985).
\item{[HW]} G.H. Hardy and E.M. Wright ``An introduction to the theory of
numbers'' Fifth Edition, Oxford Science Publications (1990).
\item{[H]} L. H\"ormander ``The Analysis of Linear Partial
Differential Operators I'' Grundlehren der mathematischen
Wissenschaften {\bf 256}, Springer--Verlag, Ber\-lin, Heidelberg,
New York, Tokyo (1983).
\item{[Ja]} N. Jacobson, ``Basic Algebra I and II$\,$'', Freeman, San
Francisco (1980).
\item{[KH]} A. Katok and B. Hasselblatt ``Introduction to the modern
theory of dynamical systems'' Encyclopedia of Mathematics and its
Applications {\bf 54}, Cambridge University Press, (1995).
\item{[Le]} L. Lewin ``Polylogarithms and Associated Functions''
Elsevier North--Holland, New York, (1981).
%\item{[Me1]} D. H. Meyer ``On a $\zeta$ function related to the
%continued fraction transformation'' {\it Bull. Soc. Math. France}
%{\bf 104} (1976), 195-203
%\item{[Me2]} D. H. Meyer ``On the Thermodynamic Formalism for
%the Gauss Map'' {\it Commun. Math. Phys.} {\bf 130} (1990), 311-333
\item{[Ma]} S. Marmi ``Critical Functions for Complex Analytic Maps''
{\it J. Phys. A: Math. Gen.} {\bf 23} (1990), 3447--3474.
\item{[MMY]} S. Marmi, P. Moussa and J.C. Yoccoz
``The Brjuno functions and their regularity properties''
{\it Commun. Math. Phys.} {\bf 186} (1997), 265--293.
\item{[O]} J. Oesterl\'e ``Polylogarithmes'' {\it S\'eminaire
Bourbaki n. 762} {\it Ast\'erisque} {\bf 216} (1993), 49--67.
\item{[Ri]} E. Risler ``Holomorphic perturbations of rotations and
rotation domains of rational functions'' preprint (March 10, 1997),
to appear in the Bulletin de la Soc. Math. France.
%\item{[Ra]} R.A. Rankin ``ÚModular forms and functions'' Cambridge
%University Press (1977)
\item{[S]} C.L. Siegel ``Iteration of analytic functions''
Annals of Mathematics $\bf{43}$ (1942), 807-812.
\item{[Se]} J.-P. Serre ``Cohomologie Galoisienne''
Lecture Notes in Mathematics,
{\bf 5}, Springer--Verlag (1973).
\item{[St]} E.M. Stein ``Singular integrals and differentiability properties
of functions'' Princeton University Press (1970).
%\item{[Ri]} G. J. Rieger ``Mischung und Ergodizit\"ata bei Kettenbruchen
%nach n\"achsten genzen'' {\it J. Reine Angew. Math.} {\bf 310} (1979),
%171-181.
\item{[Yo1]} J.-C. Yoccoz ``Th\'eor\`eme de Siegel, nombres de Bruno et
polyn\^omes quadratiques'' {\it Ast\'erisque} {\bf 231} (1995), 3--88.
\item{[Yo2]} J.-C. Yoccoz ``An introduction to small divisors problem'',
in ``From number theory to physics'', M. Waldschmidt, P. Moussa, J.-M. Luck
and C. Itzykson (editors) Springer--Verlag (1992) pp. 659--679.
\item{[Yo3]} J.-C. Yoccoz ``Analytic linearisation of analytic circle
diffeomorphisms'', in preparation.
\vfill\eject
%\input macbrc9.tex
%\vskip .5 truecm
\beginsection Figure Captions\par
\smk
%\vskip 1. truecm
\noindent
{\bf Figure 1.} The argument of the function $U$ close to the unit
circle: plot of $\arg U(0.999 e^{ix})$ for $-\pi <x<\pi $.
\mdk
%\vskip .5 truecm
\noindent
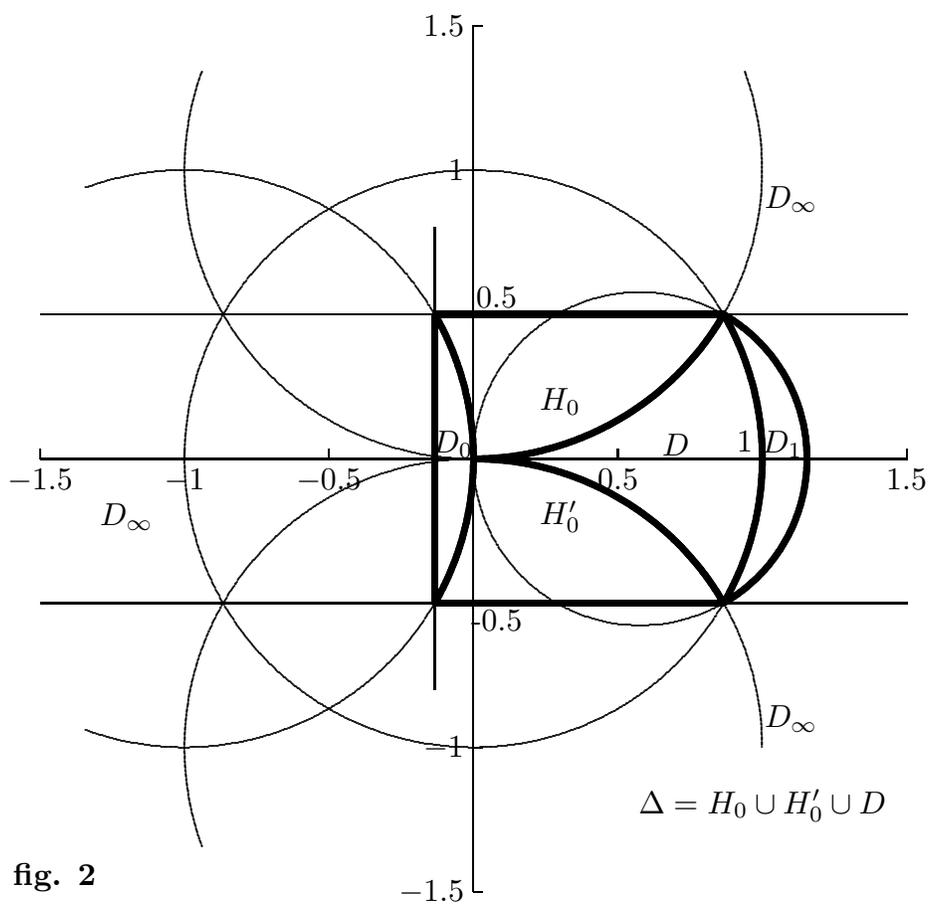
{\bf Figure 2.} The various domains $D_{i}$ used for the definition
of the complex continued fraction.
\mdk
%\vskip .5 truecm
\noindent
{\bf Figure 3.} Result of the action of $S\,:z\mapsto 1/z$ on the
domains of Figure 2.
\mdk
%\vskip .5 truecm
\noindent
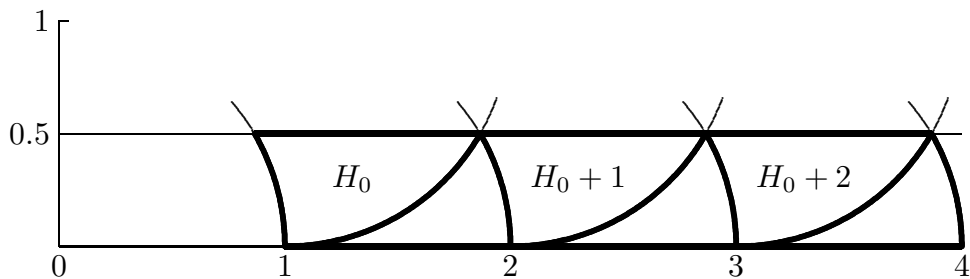
{\bf Figure 4.} The partition of the strip $0\le \IM z\le 1/2$
by the sets $D(m_{1},\ldots ,m_{r})$ and $H(m_{1},\ldots ,m_{r})$.
\bgk\bgk
%%%%%%% begin figure 1 %%%%%%
\hrule\vskip 1cm
\centerline{{\bf fig. 1}\hfil\hfil}\par\noindent
\includegraphics{brcpic1.epsi}
\vskip9.5cm
\hrule\vfill\eject
%%%%%%%%%%%% end figure 1 %%%%%%
\font\truc cmr10 scaled \magstep4
\font\struc cmr10 scaled \magstep3
\font\sstruc cmr10 scaled \magstep2
\input pictex.tex
\vglue 1pt\vfill\hrule\bigskip
%%%%%%% begin figure 2 %%%%%
$$
\beginpicture
\setcoordinatesystem units <3.2cm,3.2cm>
\setplotarea x from -1.5 to  +1.5 , y from -1.5 to 1.5
\axis bottom shiftedto y=0 ticks  short in 
withvalues $-1.5$ $-1$ $-0.5$ {} $0.5$ {} $1.5$ / quantity 7 / 
\axis left   shiftedto x=0  ticks short in 
withvalues $-1.5$ $-1$ {} {} {} $1$ $1.5$ / quantity 7 / 
\putrule from -1.5 0.5 to +1.5 0.5
\putrule from -1.5 -0.5 to +1.5 -0.5
\putrule from -0.133975 -0.8 to -0.133975 +0.8
\circulararc 110 degrees from 0 0 center at -1 0
\circulararc -110 degrees from 0 0 center at -1 0
\circulararc 110  degrees from 0 0 center at  0 1
\circulararc -110  degrees from 0 0 center at  0 1
\circulararc 110  degrees from 0 0 center at  0 -1
\circulararc -90  degrees from 0 0 center at  0 -1
\circulararc 360   degrees from 1 0 center at  0  0
\circulararc 360  degrees from 0 0 center at  0.57735 0 
\put {{\bf fig. 2}} at -1.45 -1.45
\put {1} at 0.94  +0.06
\put {0.5} at +0.08 0.56
\put {-0.5} at +0.08 -0.56
\put {$D$} at 0.7 0.05
\put {$D_0$} at -0.07 0.05
\put {$D_1$} at 1.07 0.05
\put {$D_{\infty}$} at 1.1 0.9
\put {$D_{\infty}$} at 1.1 -0.9
\put {$D_{\infty}$} at -1.2 -0.2
\put {$H_0$} at 0.3 0.2
\put {$H'_0$} at 0.3 -0.2
\put {$\Delta=H_0\cup H'_0\cup D$} at 1 -1.2 
\linethickness = 2pt
\putrule from  -0.133975 0.5 to 0.866025  0.5
\putrule from  -0.133975 -0.5 to 0.866025  -0.5
\putrule from -0.133975 -0.5 to -0.133975 0.5
\setplotsymbol({\truc .})
\circulararc 60  degrees from -0.133975 -0.5 center at  -1 0
\circulararc 60  degrees from 0.866025 -0.5 center at  0  0
\circulararc 60  degrees from 0 0 center at  0 1
\circulararc -60  degrees from 0 0 center at  0 -1
\circulararc 120  degrees from 0.866025 -0.5 center at  0.577350  0
\endpicture
$$
%% end figure2 %%%%
\bigskip\hrule
\vfill\eject\vglue1pt\vfill\hrule\bigskip
%%%%%%% begin figure3 %%%%%
$$
\beginpicture
\setcoordinatesystem units <1.6cm,1.6cm>
\setplotarea x from -3 to  +3 , y from -3 to 3
\axis bottom shiftedto y=0 ticks  short in 
withvalues $-3$ $-2$ $-1$ {} $1$ $2$ $3$ / quantity 7 / 
\axis left   shiftedto x=0  ticks short in 
withvalues $-3$ {} {}   {} {} {} $3$ / quantity 7 / 
\putrule from -3 0.5 to +3 0.5
\putrule from -3 -0.5 to +3  -0.5
\putrule from -0.5 -3 to -0.5 3 
\putrule from 0.866025 -1 to 0.866025 1
\circulararc 360 degrees from 1 0 center at 0 0
\circulararc 360 degrees from 0 0 center at 0 1 
\circulararc 360  degrees from 0 0 center at  0 -1
\circulararc -55  degrees from 0 0 center at  -3.73205 0
\circulararc +55  degrees from 0 0 center at  -3.73205 0 
\put {{\bf fig. 3}} at -2.9 -2.9
\put {1} at 1.12  0.12
\put {1} at 0.12 1.12  
\put {-1} at 0.12 -1.12 
\put {-2} at -0.12  -2.12
\put {2} at -0.12  +2.12
%\put {0.5} at +0.08 0.56
%\put {-0.5} at +0.08 -0.56
\put {$SD$} at 1.8 0.2 
\put {$SD_0$} at -0.75 2.8
\put {$SD_0$} at -0.75 -2.8
\put {$SD_1$} at 2.0 -1.39
\setlinear \plot  0.9 0.1 1.8 -1.51 2.2  -1.51 /
\put {$SD_{\infty}$} at -3 -1  
\put {$SH_0$} at  1.5 -2
\put {$SH'_0$} at 1.5 2
\linethickness = 2pt
\putrule from  0.866205 0.5 to 3  0.5
\putrule from  0.866205 -0.5 to 3  -0.5
\putrule from 0.866025 -0.5 to 0.866025 0.5
\putrule from -0.5 +1.866025 to -0.5 +3 
\putrule from -0.5 -1.866025 to -0.5 -3 
\setplotsymbol({\truc .})
\circulararc 30  degrees from 1 0 center at 0 0 
\circulararc -30  degrees from 1 0 center at 0 0 
\circulararc -120  degrees from 0 2 center at  0  1
\circulararc 30  degrees from 0 2 center at  0  1
\circulararc 120  degrees from 0 -2 center at  0 -1
\circulararc -30  degrees from 0 -2 center at  0 -1
\circulararc +25 degrees from -0.5 +1.866025 center at -3.73205 0 
\circulararc -25 degrees from -0.5 -1.866025 center at -3.73205 0 
\endpicture
$$
%% end figure3 %%%%
\bigskip \hrule
\vfill\eject
\vglue1pt
\hrule
\vskip 5 mm
%%%%%%% begin figure4 %%%%%
$$
\beginpicture
\setcoordinatesystem units <2.5cm,2.5cm>
\setplotarea x from  0  to  4 ,   y  from 0 to 1
\axis bottom shiftedto y=0 ticks  short in 
withvalues  $0$ $1$ $2$ $3$ $4$ / quantity 5 / 
\axis left   shiftedto x=0  ticks short in 
withvalues {}  $0.5$ $1$ / quantity 3 / 
\put{{\bf fig. 4}} at  2 1.2 
\putrule from  0 0 to 4 0 
\putrule from  0 0.5 to 4 0.5 
\circulararc 40 degrees from 1 0 center at 0 0
\circulararc 40 degrees from 2 0 center at 1 0
\circulararc 40 degrees from 3 0 center at 2 0
\circulararc 40 degrees from 4 0 center at 3 0
\circulararc 70 degrees from 1 0 center at 1 1
\circulararc 70 degrees from 2 0 center at 2 1
\circulararc 70 degrees from 3 0 center at 3 1
\put {$H_0$} at  1.3 0.3
\put {$H_0+1$} at  2.3 0.3 
\put {$H_0+2$} at  3.3 0.3
\linethickness = 1.66pt
\putrule from 1 0 to 4 0  
\putrule from 0.86625 0.5 to 3.86625 0.5  
\setplotsymbol({\struc  .})
\circulararc 30 degrees from 1 0 center at 0 0
\circulararc 30 degrees from 2 0 center at 1 0
\circulararc 30 degrees from 3 0 center at 2 0
\circulararc 30 degrees from 4 0 center at 3 0
\circulararc 60 degrees from 1 0 center at 1 1
\circulararc 60 degrees from 2 0 center at 2 1
\circulararc 60 degrees from 3 0 center at 3 1
%%%%%%%%%%%%%%%%%%%%
\setcoordinatesystem units <5cm,5cm> point at -0.8 1
\setplotarea x from  -0.5   to  1.5 ,   y  from 0 to 0.5 
\setplotsymbol({\fiverm .})
\linethickness = .4pt
\setlinear \plot 0.45 0.25 0.43 0.05 /
\setlinear \plot 0.25 0.23 0.31 0.025 /
\axis bottom shiftedto y=0 ticks  short in 
withvalues $-0.5$ $0$ $0.5$ $1$ $1.5$ / quantity 5 /
\put {After transformation $z\to 1/\overline{z}$\ :} at -0.3 0.4
\put {$(\sqrt{3}/2, 1/2)$} at 0.9 0.57
\put {$H(3)$} at 0.25 0.28
\put {$H(1)$} at 0.8  0.25
\put {$H(2)$} at 0.45 0.30
\setplotsymbol({\struc  .})
\circulararc 60 degrees from 0 0 center at 0 1
\circulararc -30 degrees from 1 0 center at 1 1
\circulararc -40.5 degrees from 0.5 0 center at 0.5 0.25
\circulararc -44.5 degrees from 0.33333 0 center at 0.33333 0.11111
\circulararc 30 degrees from 1 0 center at 0 0
\circulararc -10.5 degrees from 0.33333 0 center at 0.66666 0
\circulararc -14.5 degrees from 0.25 0 center at 0.375 0
\linethickness = 1.66pt
\putrule from 0 0 to 1 0  
\putrule from 0.5 0 to  0.5 0.133975  
%%%%%%%%%%%%%%%%%%%%%%%%
\setcoordinatesystem units <10cm,10cm> point at 0 1.3
\setplotarea x from  0   to  1 ,   y  from 0 to 0.5 
\put {After iteration of the process :} at 0.25 0.4
\put {$(\sqrt{3}/2, 1/2)$} at 0.9 0.54
\put {${1\over 2}$} at 0.5 -0.03
\put {${1\over 3}$} at 0.33333 -0.03
\put {${2\over 3}$} at 0.66666 -0.03
\put {${1\over 4}$} at 0.25 -0.03
\put {${3\over 4}$} at 0.75 -0.03
\put {$H(3)$} at 0.25 0.28
\put {$H(1)$} at 0.84  0.30
\put {$H(2)$} at 0.45 0.30
\put {$H(1,1)$} at 0.57 0.25
\put {$H(2,1)$} at 0.35 0.20
\put {$H(1,3)$} at 0.85 0.10
\put {$H(1,1,1)$} at 0.67 0.17
\put {$H(1,2)$} at 0.80 0.20
\setplotsymbol({\fiverm .})
\linethickness = .4pt
\setlinear \plot 0.45 0.26 0.43 0.05 /
\setlinear \plot 0.25 0.24 0.31 0.025 /
\setlinear \plot 0.57 0.21 0.55 0.05 /
\setlinear \plot 0.35 0.16 0.36 0.020 /
\setlinear \plot 0.85 0.08 0.77 0.01 /
\setlinear \plot 0.75 0.15 0.70 0.02 /
\setlinear \plot 0.63 0.14 0.65 0.02 /
\axis bottom shiftedto y=0 ticks  short in
withvalues  $0$   $1$  / quantity 2 /
\setplotsymbol({\sstruc  .})
\circulararc 60 degrees from 0 0 center at 0 1
\circulararc -30 degrees from 1 0 center at 1 1
\circulararc -40.5 degrees from 0.5 0 center at 0.5 0.25
\circulararc +40.5 degrees from 0.5 0 center at 0.5 0.25
\circulararc -44.5 degrees from 0.33333 0 center at 0.33333 0.11111
\circulararc +37   degrees from 0.33333 0 center at 0.33333 0.11111
\circulararc 30 degrees from 1 0 center at 0 0
\circulararc -10.5 degrees from 0.33333 0 center at 0.66666 0
\circulararc -14.5 degrees from 0.25 0 center at 0.375 0
\circulararc 11 degrees from 0.66666  0 center at 0.33333 0
\circulararc +44.5 degrees from 0.66666 0 center at 0.66666 0.11111
\circulararc -37   degrees from 0.66666 0 center at 0.66666 0.11111
\circulararc +14.5 degrees from 0.75 0 center at 0.625 0
\circulararc +34.5 degrees from 0.75 0 center at 0.75 0.10
\linethickness = 1.4 pt
\putrule from 0 0 to 1 0  
\putrule from 0.5 0 to  0.5 0.133975  
\endpicture
$$
%% end figure4 %%%%
\vglue 5mm
\hrule
\vfill\eject
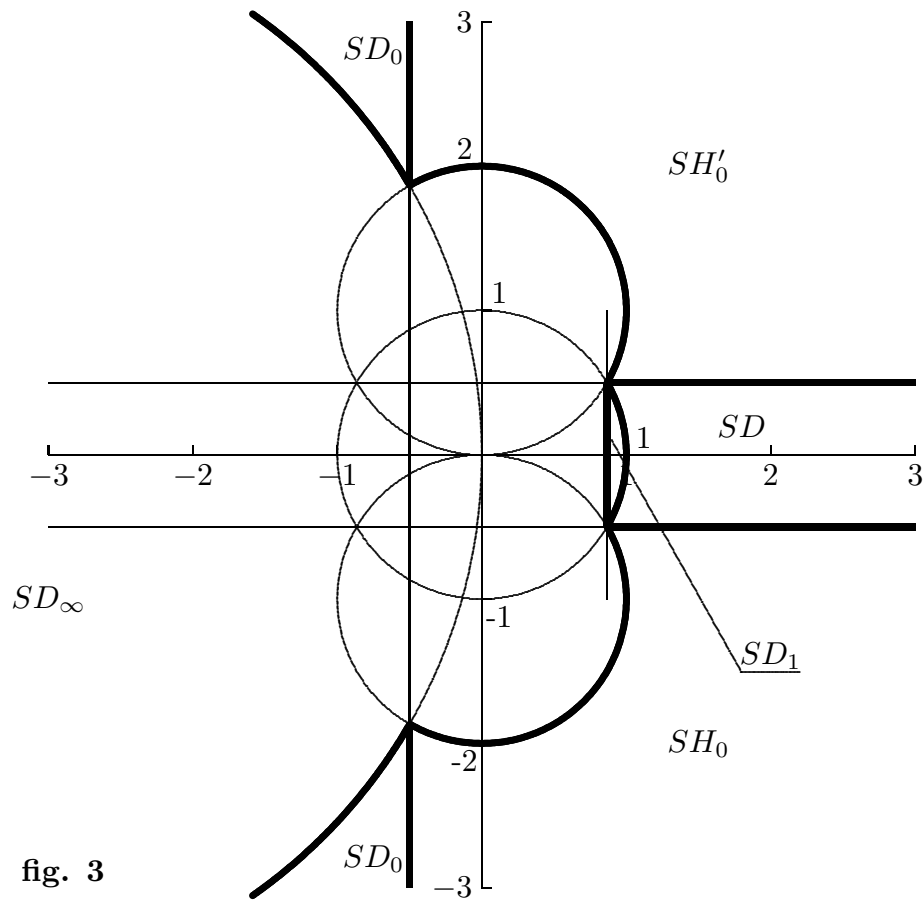
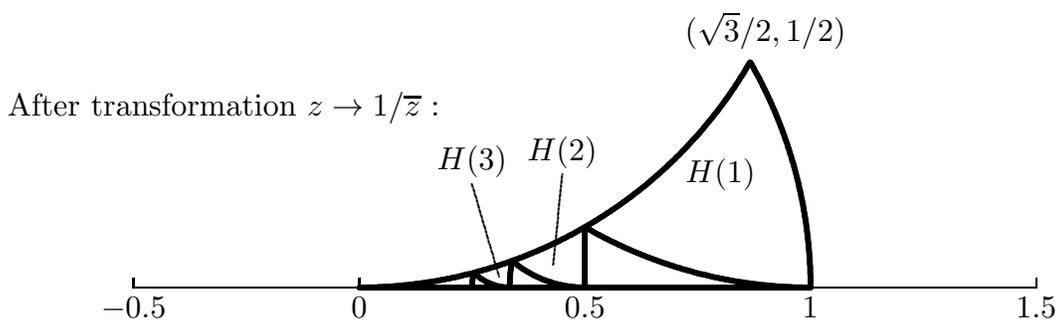
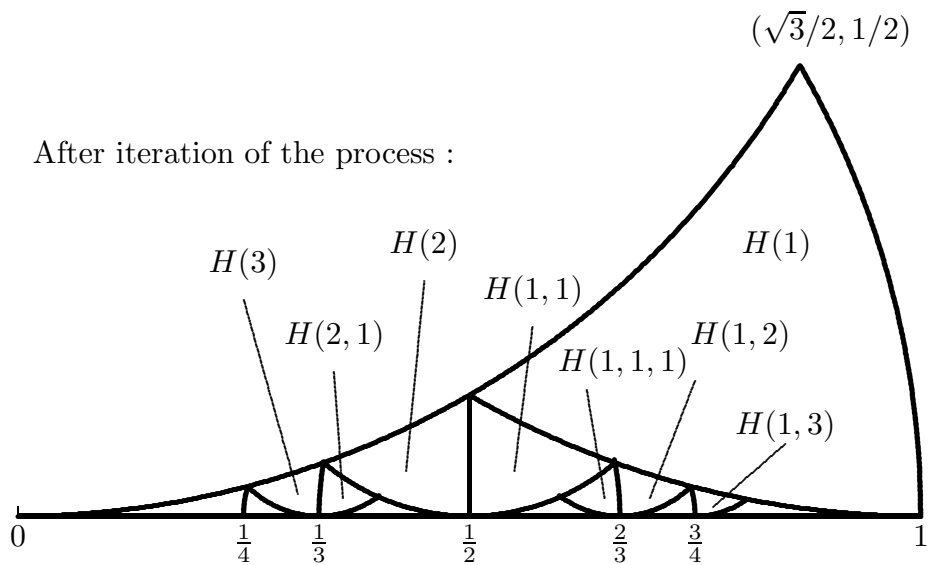
\bye